\documentclass[a4paper,10pt,twoside]{amsart}      
\usepackage{amsmath,amssymb,amsfonts,amsthm} 
\usepackage{graphics}                 
\usepackage{hyperref}                 
\usepackage[all]{xy}
\usepackage{enumerate}

\newtheorem{theorem}{Theorem}[section]
\newtheorem{lemma}[theorem]{Lemma}
\newtheorem{proposition}[theorem]{Proposition}

\newtheorem{corollary}[theorem]{Corollary}
\newtheorem{definition}[theorem]{Definition}
\newtheorem{example}[theorem]{Example}
\newtheorem{remark}[theorem]{Remark}
\numberwithin{equation}{section}

\DeclareMathOperator{\grad}{grad}
\DeclareMathOperator{\Ob}{Ob}
\DeclareMathOperator{\pd}{pd}
\DeclareMathOperator{\Ext}{Ext}
\DeclareMathOperator{\ext}{ext}
\DeclareMathOperator{\Hom}{Hom}
\DeclareMathOperator{\End}{End}
\DeclareMathOperator{\Aut}{Aut}
\DeclareMathOperator{\Stab}{Stab}

\title{A Generalized Koszul Theory and Its Application}
\author{Liping Li}
\address{School of Mathematics, University of Minnesota, MN, 55455, USA}
\email{lixxx480@math.umn.edu}
\thanks{The author wants to express great appreciation to his thesis advisor, Professor Peter Webb, for the proposal to develop a generalized Koszul theory, and the invaluable suggestions and contributions provided in numerous discussions.}
\begin{document}
\begin{abstract}
Let $A$ be a graded algebra. In this paper we develop a generalized Koszul theory by assuming that $A_0$ is self-injective instead of semisimple and generalize many classical results. The application of this generalized theory to directed categories and finite EI categories is described.
\end{abstract}
\maketitle

\section{Introduction}

The work described in this paper originated in the exploration of homological properties of finite EI categories. We want to apply Koszul theory, which has been proved to be very useful in the representation theory of algebras, to study the Ext groups of representations of finite EI categories. Examples of such applications can be found in \cite{Reiner,Woodcock}, where Koszul theory has been applied to incidence algebras of posets. In general this theory applies to graded algebras, and we do not assume that the degree 0 part of the algebra is semisimple, unlike the classical Koszul theory described in \cite{BGS,Green1,Green2,Martinez}. This generalization is necessary so that we can apply the theory to finite EI categories.

There do already exist several generalized Koszul theories where the degree 0 part $A_0$ of a graded algebra $A$ is not required to be semisimple, see \cite{Green3,Madsen1,Madsen2,Woodcock}. Each Koszul algebra $A$ defined by Woodcock in \cite{Woodcock} is supposed to satisfy that $A$ is both a left projective $A_0$-module and a right projective $A_0$-module. This requirement is too strong for us. Indeed, even the category algebra $k\mathcal{E}$ of a standardly stratified finite EI category $\mathcal{E}$ (studied in \cite{Webb}) does not satisfy this requirement: $k\mathcal{E}$ is a left projective $k\mathcal{E}_0$-module but in general not a right projective $k\mathcal{E}_0$-module. In Madsen's paper \cite{Madsen2}, $A_0$ is supposed to have finite global dimension. But for a finite EI category $\mathcal{E}$, this happens in our context if and only if $k\mathcal{E}_0$ is semisimple. The theory developed by Green, Reiten and Solberg in \cite{Green3} works in a very general framework, but some efforts are required to fill the gap between their theory and our applications.

Thus we want to develop a generalized Koszul theory which can inherit many useful results of the classical theory, and can be applied to structures with nice properties such as finite EI categories. Let $A$ be a graded $k$-algebra with $A_0$ self-injective instead of being semisimple. Then we define generalized \textit{Koszul modules} and \textit{Koszul algebras} in a similar way to the classical case. That is, a graded $A$-module $M$ is \textit{Koszul} if $M$ has a linear projective resolution, and $A$ is a \textit{Koszul algebra} if $A_0$ viewed as a graded $A$-module is Koszul. We also define \textit{quasi-Koszul modules} and \textit{quasi-Koszul algebras}: $M$ is \textit{quasi-Koszul} if the $\Ext _A^{\ast} (A_0, A_0)$-module $\Ext _A^{\ast} (M, A_0)$ is generated in degree 0, and $A$ is a \textit{quasi-Koszul algebra} if $A_0$ is a quasi-Koszul $A$-module. It turns out that this generalization works nicely for our goal. Many classical results described in \cite{BGS,Green1,Green2,Martinez} generalize to our context. In particular, we obtain the Koszul duality.

We then focus on the applications of this generalized Koszul theory. First we define \textit{directed categories}. A directed category $\mathcal{C}$ is a $k$-linear category equipped with a partial order $\leqslant$ on $\Ob \mathcal{C}$ such that for each pair of objects $x, y \in \Ob \mathcal{C}$, the space of morphisms $\mathcal{C} (x,y)$ is non-zero only if $x \leqslant y$. Directed categories include the $k$-linearizations of finite EI categories as special examples. This partial order determines a canonical pre-order $\preccurlyeq$ on the isomorphism classes of simple representations. Following the technique in \cite{Webb}, we develop a stratification theory for directed categories, describe the structures of standard modules and characterize every directed category $\mathcal{C}$ standardly stratified with respect to the canonical pre-order:

\begin{theorem}
Let $\mathcal{C}$ be a directed category with respect to a partial order $\leqslant$. Then $\mathcal{C}$ is standardly stratified for the canonical pre-order $\preccurlyeq$ if and only if the morphism space $\mathcal{C} (x,y)$ is a projective $\mathcal{C} (y,y)$-module for every pair of objects $x, y \in \Ob \mathcal{C}$.
\end{theorem}

By the correspondence between graded $k$-linear categories and graded algebras described in \cite{Mazorchuk1}, we can view a directed category as a directed algebra and vice-versa. Therefore, all of our results on graded algebras can be applied to graded directed categories. In particular, the following theorem relates Koszul theory to stratification theory:

\begin{theorem}
Let $\mathcal{C}$ be a graded directed category with $\mathcal{C}_0 = \bigoplus _{x \in \Ob \mathcal{C}} \mathcal{C} (x, x)$ self-injective. Then:
\begin{enumerate}
\item $\mathcal{C}$ is standardly stratified with respect to the canonical pre-order $\preccurlyeq$ if and only if $\mathcal{C}$ is a projective $\mathcal{C}_0$-module.
\item $\mathcal{C}$ is a Koszul category if and only if $\mathcal{C}$ is a quasi-Koszul category standardly stratified for $\preccurlyeq$.
\item If $\mathcal{C}$ is standardly stratified for $\preccurlyeq$, then a graded $\mathcal{C}$-module $M$ generated in degree 0 is Koszul if and only if it is a quasi-Koszul $\mathcal{C}$-module and a projective $\mathcal{C}_0$-module.
\end{enumerate}
\end{theorem}

Applying the homological dual functor $E = \Ext _{\mathcal{C}} ^{\ast} (-, \mathcal{C}_0)$ to a graded directed category $\mathcal{C}$, we construct the \textit{Yoneda category} $E(\mathcal{C}_0) = \Ext ^{\ast} _{\mathcal{C}} (\mathcal{C}_0, \mathcal{C}_0)$. We prove that if $\mathcal{C}$ is a Koszul directed category with $\mathcal{C}_0$ self-injective, then $E(\mathcal{C}_0)$ is also a Koszul directed category.

We acquire a very nice correspondence between the classical Koszul theory and our generalized Koszul theory for directed categories.

\begin{theorem}
Let $\mathcal{C}$ be a graded directed category with $\mathcal{C}_0$ self-injective. Define $\mathcal{D}$ to be the subcategory of $\mathcal{C}$ obtained by replacing each endomorphism ring by $k\cdot 1$, the span of the identity endomorphism. Then:\
\begin{enumerate}
\item $\mathcal{C}$ is a Koszul category in our sense if and only if $\mathcal{C}$ is standardly stratified for the canonical pre-order $\preccurlyeq$ and $\mathcal{D}$ is a Koszul category in the classical sense.
\item If $\mathcal{C}$ is a Koszul category, then a graded $\mathcal{C}$-module $M$ is Koszul if and only if $M$ is a projective $\mathcal{C}_0$-module, and $M \downarrow _{\mathcal{D}} ^{\mathcal{C}}$ is a Koszul $\mathcal{D}$-module.
\end{enumerate}
\end{theorem}

Finite EI categories have nice combinatorial properties. These properties can be used to define length gradings on the sets of morphisms. We discuss the possibility to put such a grading on an arbitrary finite EI category, and prove the following result for \textit{finite free EI categories} (defined in \cite{Li}):

\begin{theorem}
Let $\mathcal{E}$ be a finite free EI category. Then the following are equivalent:
\begin{enumerate}
\item $\pd _{k\mathcal{E}} k\mathcal{E} _0 < \infty$;
\item $\pd _{k\mathcal{E}} k\mathcal{E}_0 \leqslant 1$;
\item $\mathcal{E}$ is standardly stratified in a sense defined in \cite{Webb};
\item $k \mathcal{E}$ is a Koszul algebra.
\end{enumerate}
Moreover, $k \mathcal{E}$ is quasi-Koszul if and only if $\Ext _{k \mathcal{E}} ^i (k\mathcal{E} _0, k\mathcal{E} _0) =0$ for $i \geqslant 2$.
\end{theorem}

We then give a sufficient condition for the category algebra of a finite EI category $\mathcal{E}$ to be quasi-Koszul. An object $x \in \Ob \mathcal{E}$ is called \textit{left regular} if for every morphism $\alpha: z \rightarrow x$, the stabilizer of $\alpha$ in $\Aut _{\mathcal {E}} (y)$ has order invertible in $k$. Similarly we define \textit{right regular} objects. Then:

\begin{theorem}
Let $\mathcal{E}$ be a finite free EI category. If every object in $\mathcal{E}$ is either left regular or right regular, then $k \mathcal{E}$ is a quasi-Koszul algebra.
\end{theorem}

Motivated by the fact that the Yoneda category $E(\mathcal{C}_0)$ of a directed Koszul category $\mathcal{C}$ is still a directed Koszul category, and hence is standardly stratified, we ask whether the Koszul dual algebra $\Gamma = \Ext _A ^{\ast} (A_0, A_0)$ of a graded algebra $A$ standardly stratified for a partial order $\leqslant$ is still standardly stratified for $\leqslant$ (or $\leqslant ^{\textnormal{op}}$). This question has only been studied for the case that $A_0$ is semisimple, see \cite{Agoston1,Agoston2,Drozd,Mazorchuk2,Mazorchuk3}. By assuming that $A_0$ is a self-injective algebra, and supposing that all standard modules are concentrated in degree 0 and Koszul, we get a sufficient condition for $\Gamma$ to be standardly stratified with respect to $\leqslant ^{\textnormal{op}}$.

The layout of this paper is as follows. The generalized Koszul theory is developed in the first three sections. In Section 2 we define Koszul modules and quasi-Koszul modules, which generalize Koszul modules in the classical theory, and describe their basic properties. Since Koszul modules and quasi-Koszul modules do not coincide in our context, we also give a relation between these two concepts. Koszul algebras are studied in Section 3. The Koszul dualities are proved in Section 4. Most results in these three sections are generalized from works in \cite{BGS,Green1,Green2,Martinez}. Some results can be deduced from the paper \cite{Green3} of Green, Reiten and Solberg, who worked in a more general context, but we present full arguments for the sake of completeness.

The last three sections are on the application of the general theory developed before. Directed categories are defined in Section 5. Their stratification properties and Koszul properties are discussed in details in this section as well. The main content of Section 6 is to apply the Koszul theory and stratification theory to finite EI categories, which have nice combinatorial structures. In Section 7 we modify the technique of \cite{Agoston1} to study standardly stratified algebras with Koszul standard modules.

Here are the notation and conventions we use in this paper. All algebras are $k$-algebras with $k$ an algebraically closed field. Let $A$ be a graded algebra. A graded $A$-module $M$ is said to be \textit{locally finite} if $\dim _k M_i < \infty$ for each $i \in \mathbb{Z}$. In this paper we only consider locally finite graded modules, for which the category is denoted by $A$-gmod. For $M, N \in A$-gmod, by $\Hom_A (M,N)$ and $\hom_A (M, N)$ we denote the space of module homomorphisms and the space of \textit{graded module homomorphisms} respectively. The $s$-th \textit{shift} $M[s]$ of $M$ is defined in the following way: $M[s]_i = M_{i-s}$ for all $i \in \mathbb{Z}$. If $M$ is generated in degree $s$, then $\bigoplus _{i \geqslant s+1} M_i$ is a graded submodule of $M$, and $M_s \cong M/ \bigoplus _{i \geqslant s+1} M_i$ as vector spaces. We then view $M_s$ as an $A$-module by identifying it with this quotient module. We also regard the zero module 0 as a projective module since with this convention the expressions of many results can be simplified.

\section{Generalized Koszul Modules}

Throughout this section $A$ is a non-negatively graded and locally finite associative $k$-algebra with identity 1 generated in degrees 0 and 1, i.e., $A = \bigoplus _{i=0}^{\infty} A_i$ such that $A_i \cdot A_j = A_{i+j}$ for all $i, j \geqslant 0$; each $A_i$ is finite-dimensional. We also suppose that $A_0$ is a \textit{self-injective} algebra, i.e., every projective $A_0$-module is injective as well. Define $J = \bigoplus _{i=1}^{\infty} A_i$, which is a two-sided ideal of $A$. An $A$-module $M$ is called \textit{graded} if $M = \bigoplus_{i \in \mathbb{Z}} M_i$ such that $A_i \cdot M_j \subseteq M_{i+j}$. We say $M$ is \textit{generated in degree s} if $M = A \cdot M_s$. It is clear that $M$ is generated in degree $s$ if and only if $JM \cong \bigoplus _{i \geqslant s+1} M_i$, which is equivalent to $J^lM \cong \bigoplus _{i \geqslant s+l} M_i$ for all $l \geqslant 1$.

Most results in this section are generalized from \cite{Green1,Green2,Martinez}.

We collect some preliminary results in the following lemma.
\begin{lemma}
Let $A$ be as above and $M$ be a locally finite graded $A$-module. Then:
\begin{enumerate}
\item $J = \bigoplus _{i \geqslant 1} A_i$ is contained in the graded radical of $A$;
\item $M$ has a graded projective cover;
\item the graded syzygy $\Omega M$ is also locally finite.
\end{enumerate}
\end{lemma}

\begin{proof}
By definition, the graded radical $\grad A$ is the intersection of all maximal proper graded submodules of $A$. Let $L \subsetneqq A$ be a maximal proper graded submodule. Then $L_0$ is a proper subspace of $A_0$. We claim that $A_i = L_i$ for all $i \geqslant 1$. Otherwise, we can define $\tilde{L} \subseteq A$ in the following way: $\tilde{L_0} = L_0$ and $\tilde{L}_i = A_i$ for $i \geqslant 1$. Then $L \subsetneqq \tilde{L} \subsetneqq A$, so $L$ is not a maximal proper graded submodule of $A$. This contradiction tells us that $L_i = A_i$ for all $i \geqslant 1$. Therefore, $J \subseteq \grad A$, and the first statement is proved.

We use the following fact to prove the second statement: every primitive idempotent in the algebra $A_0$ can be lifted to a primitive idempotent of $A$. Consequently, a projective $A_0$-module concentrated in some degree $d$ can be lifted to a graded projective $A$-module generated in degree $d$.

Define $\bar{M} = M/JM$, which is also a locally finite graded $A$-module. Write $\bar{M} = \bigoplus _{i \geqslant 0} \bar{M}_i$. Then each $\bar{M}_i$ is a finite-dimensional graded $A$-module since $J \bar{M} =0$ and $A_0 \bar{M}_i = \bar{M} _i$ for all $i \geqslant 0$. Therefore, $\bar{M}$ can be decomposed as a direct sum of indecomposable graded $A$-modules each of which is concentrated in a certain degree. Moreover, for each $i \in \mathbb{Z}$, there are only finitely many summands concentrated in degree $i$. 

Take $L$ to be such an indecomposable summand and without loss of generality suppose that it is concentrated in degree 0. As an $A_0$-module, $L$ has a finitely generated projective cover $P_0$. By the lifting property, $P_0$ can be lifted to a finitely generated graded projective module $P$ generated in degree 0, which is a graded projective cover of $L$. Take the direct sum of these projective covers $P$ when $L$ ranges over all indecomposable summands of $\bar{M}$. In this way we obtain a graded projective cover $\tilde{P}$ of $\bar{M}$. Clearly, $\tilde{P}$ is also a graded projective cover of $M$. The second statement is proved.

Now we turn to the third statement. By the above proof, the graded projective cover $\tilde{P}$ of $M$ can be written as a direct sum $\bigoplus _{i \geqslant 0} P^i$ of graded projective modules, where $P^i$ is generated in degree $i$. For each fixed degree $i \geqslant 0$, there are only finitely many indecomposable summands $L$ of $\bar{M}$ concentrated in degree $i$, and the graded projective cover of each $L$ is finitely generated. Consequently, $P^i$ is finitely generated, and hence locally finite.

For a fixed $n \geqslant 0$, we have $\tilde{P}_n = \bigoplus _{i \geqslant 0} P^i_n = \bigoplus _{0 \leqslant i \leqslant n} P^i_n$. Since each $P^i$ is locally finite, $\dim_k P^i_n < \infty$. Therefore, $\dim_k \tilde{P}_n < \infty$, and $\tilde{P}$ is locally finite as well. As a submodule of $\tilde{P}$, the graded syzygy $\Omega M$ is also locally finite.
\end{proof}

These results will be used frequently.

\begin{lemma}
Let $0 \rightarrow L \rightarrow M \rightarrow N \rightarrow 0$ be an exact sequence of graded $A$-modules.
\begin{enumerate}
  \item If $M$ is generated in degree $s$, so is $N$.
  \item If $L$ and $N$ are generated in degree $s$, so is $M$.
  \item If $M$ is generated in degree $s$, then $L$ is generated in degree $s$ if and only if $JM \cap L = JL$.
\end{enumerate}
\end{lemma}

\begin{proof}
(1): This is obvious.

(2): Let $P$ and $Q$ be graded projective covers of $L$ and $N$ respectively. Then $P$ and $Q$, and hence $P \oplus Q$ are generated in degree $s$. In particular, each graded projective cover of $M$, which is isomorphic to a direct summand of $P \oplus Q$, is generated in degree $s$. Thus $M$ is also generated in degree $s$.

(3): We always have $JL \subseteq JM \cap L$. Let $x \in L \cap JM$ be a homogeneous element of degree $i$. Since $M$ is generated in degree $s$, we have $i \geqslant s+1$. If $L$ is generated in degree $s$, then $x \in J^{i-s}L \subseteq JL$. Thus $L \cap JM \subseteq JL$, so $JL = L \cap JM$.

Conversely, the identity $JL = L \cap JM$ gives us the following commutative diagram where all rows and columns are exact:
\begin{align*}
\xymatrix{
& 0 \ar[d] & 0 \ar[d] & 0 \ar[d] &  \\
0 \ar[r] & JL \ar[r] \ar[d] & JM \ar[r] \ar[d] & JN \ar[r] \ar[d] & 0\\
0 \ar[r] & L \ar[r] \ar[d] & M \ar[r] \ar[d] & N \ar[r] \ar[d] & 0 \\
0 \ar[r] & L/JL \ar[r] \ar[d] & M/JM \ar[r] \ar[d] & N/JN \ar[r] \ar[d] & 0 \\
& 0 & 0 & 0 &
}
\end{align*}

Consider the bottom sequence. Notice that $(M / JM) \cong M_s $ is concentrated in degree $s$. Thus $L/JL$ is also concentrated in degree $s$, i.e., $L/JL \cong L_s$. Let $I = A \cdot L_s$. Then $L \subseteq I + JL \subseteq L$, so $I + JL = L$. Notice that $J$ is contained in the graded Jacobson radical of $A$. Therefore, by the graded Nakayama lemma, $I = A \cdot L_s = L$, so $L$ is generated in degree $s$.
\end{proof}

\begin{corollary}
Suppose that each graded $A$-module in the short exact sequence $0 \rightarrow L \rightarrow M \rightarrow N \rightarrow 0$ is generated in degree $0$. Then $J^iM \cap L = J^iL$ for all $i \geqslant 0$.
\end{corollary}

\begin{proof}
Since all modules $L$, $M$ and $N$ are generated in degree 0, all $J^s L$, $J^s M$ and $J^s N$ are generated in degree $s$ for $s \geqslant 0$. The exactness of the above sequence implies $JL = L \cap JM$, which in turns gives the exactness of $0 \rightarrow JL \rightarrow JM \rightarrow JN \rightarrow 0$. By the above lemma, $J^2 M \cap JL = J^2 L$ and this implies the exactness of $0 \rightarrow J^2L \rightarrow J^2M \rightarrow J^2N \rightarrow 0$. The conclusion follows from induction.
\end{proof}

Now we introduce generalized \textit{Koszul modules} (or called \textit{linear modules}).

\begin{definition}
A graded $A$-module $M$ generated in degree 0 is called a Koszul module if it has a (minimal) projective resolution
\begin{displaymath}
\xymatrix{ \ldots \ar[r] & P^n \ar[r] & P^{n-1} \ar[r] & \ldots \ar[r] & P^1 \ar[r] & P^0 \ar[r] & M \ar[r] & 0}
\end{displaymath}
such that $P^i$ is generated in degree $i$ for all $i \geqslant 0$.
\end{definition}

A direct consequence of this definition and the previous lemma is:

\begin{corollary}
Let $M$ be a Koszul module. Then $\Omega^i(M) /J\Omega^i (M) \cong \Omega^i(M)_i$ is a projective $A_0$-module for each $i \geqslant 0$, or equivalently, $\Omega^i (M) \subseteq JP^{i-1}$ where $P^{i-1}$ is a graded projective cover of $\Omega^{i-1} (M)$, where $\Omega$ is the Heller operator.
\end{corollary}

\begin{proof}  Since $M$ is Koszul, $\Omega^i (M)$ is generated in degree $i$, and $\Omega^i(M) /J\Omega^i (M) \cong \Omega^i(M)_i$. Moreover, all $\Omega^i(M)[-i]$ are Koszul $A$-modules as well. By induction, it is sufficient to prove $\Omega M \subseteq JP^0$. But this is obvious since $\Omega M$ is generated in degree 1. From the following commutative diagram we deduce that $\Omega M \subseteq JP^0$ if and only if the bottom sequence is exact, or equivalently $M/JM \cong P^0/JP^0 \cong P^0_0$ is a projective $A_0$-module.
\begin{align*}
\xymatrix{ 0 \ar[r] & \Omega M \ar@{=}[d] \ar[r] & JP^0 \ar[d] \ar[r] & JM \ar[d] \ar[r] & 0 \\
0 \ar[r] & \Omega M \ar[r] \ar[d] & P^0 \ar[r] \ar[d] & M \ar[r] \ar[d] & 0\\
0 \ar[r] & 0 \ar[r] & P^0/JP^0 \ar[r] & M/JM \ar[r] & 0 }
\end{align*}
\end{proof}

There are several characterizations of Koszul modules.

\begin{proposition}
Let $M$ be a graded $A$-module generated in degree 0. Then the following are equivalent:
\begin{enumerate}
  \item $M$ is Koszul.
  \item The syzygy $\Omega^i(M)$ is generated in degree $i$ for every $i \geqslant 0$.
  \item For all $i>0$, $\Omega^i(M) \subseteq JP^{i-1}$ and $\Omega^i(M) \cap J^2P^{i-1} = J \Omega^i(M)$, where $P^{i-1}$ is a graded projective cover of $\Omega^{i-1} (M)$.
  \item $\Omega^i(M) \subseteq JP^{i-1}$ and $\Omega^i(M) \cap J^{s+1} P^{i-1} = J^s \Omega^i(M)$ for all $i>0, s \geqslant 0$.
\end{enumerate}
\end{proposition}

\begin{proof}
The equivalence of (1) and (2) is clear. It is also obvious that (3) is the special case of (4) for $s=1$. Now we show (1) implies (4). Indeed, if $M$ is a Koszul module, then both $JP^0$ and $\Omega M$ are generated in degree 1 and $\Omega M \subseteq JP^0$. Therefore we have the following exact sequence
\begin{equation*}
\xymatrix@1 {0 \ar[r] & \Omega M \ar[r] & JP^0 \ar[r] & JM \ar[r] & 0}
\end{equation*}
in which all modules are generated in degree 1. By Corollary 2.5 $J^{s+1} P^0 \cap \Omega M = J^s \Omega M$ for all $s >0$. Notice that all syzygies of $M$ are also Koszul with suitable grade shifts. Replacing $M$ by $\Omega^i(M)[-i]$ and using induction we get (4).

Finally we show (3) implies (2) to finish the proof. Since $\Omega M \subseteq JP^0$ we still have the above exact sequence. Notice that both $JM$ and $JP^0$ are generated in degree 1 and $J^2P^0 \cap \Omega M = J\Omega M$, by Lemma 2.2, $\Omega M$ is generated in degree 1 as well. Now the induction procedure gives us the required conclusion.
\end{proof}

The condition that $\Omega^i(M) \subseteq JP^{i-1}$ (or equivalently, $\Omega^i(M)/J\Omega^i(M) \cong \Omega^i (M)_i$ is a projective $A_0$-module) in (3) of the previous proposition is necessary, as shown by the following example:

\begin{example}
Let $G$ be a finite cyclic group of prime order $p$ and $k$ be an algebraically closed field of characteristic $p$. Let the group algebra $kG$ be concentrated on degree 0, so $J=0$. Consider the trivial $kG$-module $k$. Obviously, $k$ is not a Koszul module. But since $J=0$, the condition $J \Omega^i(k) = J^2P^{i-1} \cap \Omega^i(k)$ holds trivially.
\end{example}

\begin{remark}
We do not use the property that $A_0$ is a self-injective algebra up to now. Therefore, all results described before still hold for a non-negatively graded, locally finite graded algebra $A$ with $A_0$ an arbitrary finite-dimensional algebra.
\end{remark}

\begin{proposition}
Let $0 \rightarrow L \rightarrow M \rightarrow N \rightarrow 0$ be a short exact sequence of graded $A$-modules such that $L$ is a Koszul $A$-module. Then $M$ is Koszul if and only if $N$ is Koszul.
\end{proposition}

\begin{proof}
We verify the conclusion by using statement (2) in the last proposition. That is, given that $\Omega^i(L)$ is generated in degree $i$ for each $i \geq 0$, we want to show that $\Omega^i(M)$ is generated in degree $i$ if and only if so is $\Omega^i (N)$.

Consider the given exact sequence. By Lemma 2.2, $M$ is generated in degree 0 if and only if $N$ is generated in degree 0. Therefore we have the following diagram in which all rows and columns are exact:
\begin{align}
\xymatrix{
& 0 \ar[d] & 0 \ar[d] & 0 \ar[d] &  \\
0 \ar[r] & \Omega L \ar[r] \ar[d] & M' \ar[r] \ar[d] & \Omega N \ar[r] \ar[d] & 0\\
0 \ar[r] & P \ar[r] \ar[d] & P \oplus Q \ar[r] \ar[d] & Q \ar[r] \ar[d] & 0 \\
0 \ar[r] & L \ar[r] \ar[d] & M \ar[r] \ar[d] & N \ar[r] \ar[d] & 0 \\
& 0 & 0 & 0 &
}
\end{align}
where $P$ and $Q$ are graded projective covers of $L$ and $N$ respectively. In general $P \oplus Q$ is not a graded projective cover of $M$, and hence $M' \ncong \Omega M$. However, under the hypothesis of this proposition we claim that $P \oplus Q$ is indeed a graded projective cover of $M$ and $M' \cong \Omega M$. To see this, we point out that the given exact sequence induces a short exact sequence of $A_0$-modules:
\begin{align*}
\xymatrix { 0 \ar[r] & L_0 \ar[r] & M_0 \ar[r] & N_0 \ar[r] & 0. }
\end{align*}
Observe that $L_0 = \Omega^0 (L)_0$ is projective (by Corollary 2.5) and hence injective (since $A_0$ is self-injective) as an $A_0$-module. Therefore, the above sequence splits, and $M_0 \cong L_0 \oplus N_0$. On one hand, a graded projective cover $\tilde{P}$ of $M$ should be isomorphic to a direct summand of $P \oplus Q$. On the other hand, since $\tilde{P}$ induces a projective cover $\tilde{P}_0$ of $M_0 \cong L_0 \oplus N_0$, $\tilde{P}$ should contain a summand isomorphic to $Q \oplus P$. This forces $\tilde{P} \cong P \oplus Q$, and our claim is proved.

Now let us consider the top row of diagram 2.1. Since $M' \cong \Omega M$, and $\Omega L$ is generated in degree 1, $\Omega M$ is generated in degree 1 if and only if $\Omega N$ is generated in degree 1 by Lemma 2.2. Replace $L$, $M$ and $N$ by $(\Omega L)[-1]$, $(\Omega M)[-1]$ and $(\Omega N)[-1]$ (all of which are Koszul) respectively in the short exact sequence. Repeating the above procedure and using the fact that the Heller operator $\Omega$ and the grade shift functor $[ - ]$ commute, we conclude that $\Omega^2(M)[-1]$ is generated in degree $1$ if and only if $\Omega^2(N)[-1]$ is generated in degree $1$, i.e., $\Omega^2(M)$ is generated in degree $2$ if and only if $\Omega^2(N)$ is generated in degree $2$. By induction, $M$ is Koszul if and only if $N$ is Koszul.
\end{proof}

The condition that $L$ is Koszul in this proposition is necessary. Indeed, quotient modules of a Koszul module might not be Koszul.

\begin{lemma}
Let $M$ be a graded $A$-module generated in degree $s$. If $M_s$ is a projective $A_0$-module, then $\Ext _{A}^i (M, A_0) \cong \Ext _{A}^{i-1} (\Omega M, A_0)$ for all $i \geqslant 1$.
\end{lemma}

\begin{proof}
It is true for $i >1$. When $i=1$, consider the following exact sequence:
\begin{equation*}
0 \rightarrow \Hom_A (M, A_0) \rightarrow \Hom_A (P, A_0) \rightarrow \Hom_A (\Omega M, A_0) \rightarrow \Ext^1_A (M, A_0) \rightarrow 0.
\end{equation*}
As a graded projective cover of $M$, $P$ is also generated in degree $s$. Since $M_s$ is a projective $A_0$-module, $P_s \cong M_s$. So
\begin{equation*}
\Hom_A (M, A_0) \cong \Hom_{A_0} (M_s, A_0) \cong \Hom_{A_0} (P_s, A_0) \cong \Hom_A (P, A_0).
\end{equation*}
Thus $\Hom_A (\Omega M, A_0) \cong \Ext^1 _{A} (M, A_0)$.
\end{proof}

\begin{proposition}
Let $0 \rightarrow L \rightarrow M \rightarrow N \rightarrow 0$ be an exact sequence of Koszul modules. Then it induces the following short exact sequence:
\begin{equation*}
\xymatrix{ 0 \ar[r] & \Ext ^{\ast} _A (N, A_0) \ar[r] & \Ext ^{\ast} _A (M, A_0) \ar[r] & \Ext ^{\ast} _A (L, A_0) \ar[r] & 0.}
\end{equation*}
\end{proposition}

\begin{proof}
As in the proof of Proposition 2.9, the above exact sequence gives exact sequences $0 \rightarrow \Omega^i (L) \rightarrow \Omega^i (M) \rightarrow \Omega^i (N) \rightarrow 0$, $i \geqslant 0$. For a fixed $i$, the sequence $ 0 \rightarrow \Omega^i (L)_i \rightarrow \Omega^i(M)_i \rightarrow \Omega^i (N)_i \rightarrow 0$ splits since all terms are projective $A_0$-modules by Corollary 2.5. Applying the exact functor $\Hom _{A_0} (-, A_0)$ we get an exact sequence
\begin{equation*}
0 \rightarrow \Hom _{A_0} (\Omega^i(N)_i, A_0) \rightarrow \Hom _{A_0} (\Omega^i(M)_i, A_0) \rightarrow \Hom _{A_0} (\Omega^i(L)_i, A_0) \rightarrow 0
\end{equation*}
which is isomorphic to
\begin{equation*}
0 \rightarrow \Hom_A (\Omega^i(N), A_0) \rightarrow \Hom_A (\Omega^i(M), A_0) \rightarrow \Hom_A (\Omega^i(L), A_0) \rightarrow 0
\end{equation*}
since all modules are generated in degree $i$. By Lemma 1.9, it is isomorphic to
\begin{equation*}
\xymatrix{ 0 \ar[r] & \Ext^i _A (N, A_0) \ar[r] & \Ext^i _A (M, A_0) \ar[r] & \Ext^i _A (L, A_0) \ar[r] & 0.}
\end{equation*}
Putting them together we have:
\begin{equation*}
\xymatrix{ 0 \ar[r] & \Ext ^{\ast} _A (N, A_0) \ar[r] & \Ext ^{\ast} _A (M, A_0) \ar[r] & \Ext ^{\ast} _A (L, A_0) \ar[r] & 0.}
\end{equation*}
\end{proof}

If a graded $A$-module $M$ is Koszul, so are all syzygies $\Omega^i(M)[-i]$, $i \geqslant 0$. The next example shows that the Koszul property of $M$ in general does not imply the Koszul property of $J^i M [-i]$:

\begin{example}
Let $\mathcal{E}$ be a finite EI category with two objects $x$ and $y$ such that: $\Aut _{\mathcal{E}} (x) = \langle g \rangle \cong \Aut _{\mathcal{E}} (y) = \langle h\rangle$ are cyclic groups of order 2; $\Hom _{\mathcal{E}} (x,y)$ has one element $\alpha$ on which both $\Aut _{\mathcal{E}} (x)$ and $\Aut _{\mathcal{E}} (y)$ act trivially; and $\Hom _{\mathcal{E}} (y,x) = \emptyset$. Let $k$ be an algebraically closed field of characteristic 2. We put the following grading on the category algebra $A = k \mathcal{E}$: $A_0$ is spanned by $\{g, 1_x, h, 1_y \}$ and $A_1$ is spanned by $\alpha$. Consider the projective $k\mathcal{E}$-module $P_x = k\mathcal{E} 1_x$. Obviously, $P_x$ is Koszul, but $JP_x \cong k_y$ is not Koszul.
\end{example}

However, for some special cases, we can get a conclusion as follows.

\begin{proposition}
Suppose that $A_0$ is a Koszul $A$-module. If $M$ is a Koszul $A$-module, then $J^iM[-i]$ and $M_0$ are also Koszul $A$-modules. In particular, $M$ is projective viewed as an $A_0$-module (or equivalently, $M_i$ is a projective $A_0$-module for every $i \geqslant 0$).
\end{proposition}

\begin{proof}
Without loss of generality we assume that $M$ is indecomposable. Since $M$ is Koszul, $M_0$ is a projective $A_0$-module and is contained in add$ (A_0)$, the category of all $A$-modules each of which is isomorphic to a direct summand of $A_0 ^{\oplus m}$ for some $m \geqslant 0$. But $A_0$ is Koszul, so is $M_0$.

Notice that $M$ and $M_0$ have the same graded projective cover (up to isomorphism) as $A$-modules. Thus we have the following commutative diagram:
\begin{align*}
\xymatrix{
& 0 \ar[r] & \Omega M \ar[r] \ar[d] & \Omega(M_0) \ar[r] \ar[d] & JM \ar[r] & 0 \\
& 0 \ar[r] & P^0 \ar[r]^{id} \ar[d] & P^0 \ar[r] \ar[d] & 0 & \\
0 \ar[r] & JM \ar[r] & M \ar[r] & M_0 \ar[r] & 0 &}
\end{align*}
Consider the top sequence. The modules $M$ and $M_0$ are Koszul, so are $\Omega M[-1]$ and $\Omega(M_0)[-1]$. Therefore, $JM[-1]$ is also Koszul by Proposition 2.9. Now replacing $M$ by $JM[-1]$ and using recursion, we conclude that $J^iM[-i]$ is a Koszul $A$-module for every $i \geqslant 0$. This proves the first statement.

Since $J^iM[-i]$ is a Koszul $A$-module, $(J^iM[-i])_0 \cong (J^iM)_i \cong M_i$ is a projective $A_0$-module for each $i \geqslant 0$ by Corollary 2.5. Since $M = \bigoplus _{i \geqslant 0} M_i$ as an $A_0$-module, we deduce that $M$ is a projective $A_0$-module if and only if $M_i$ is a projective $A_0$-module for every $i \geqslant 0$.
\end{proof}

The following lemma is very useful.

\begin{lemma}
Let $M$ be a non-negatively graded $A$-module and suppose that $A$ is a projective $A_0$-module. Then the following are equivalent:
\begin{enumerate}
\item all $\Omega^i(M)_j$ are projective $A_0$-modules, $i, j \geqslant 0$;
\item all $\Omega^i(M)_i$ are projective $A_0$-modules, $i \geqslant 0$;
\item all $M_i$ are projective $A_0$-modules, $i \geqslant 0$.
\end{enumerate}
\end{lemma}

\begin{proof}
It is clear that (1) implies (2).

(3) implies (1): Suppose that $M_i$ is a projective $A_0$-module for all $i \geqslant 0$. Let $P$ be a graded projective cover of $M$. The surjective homomorphism $\varphi: P \rightarrow M$ gives a surjective homomorphism $\varphi_j: P_j \rightarrow M_j$ with kernel $(\Omega M)_j$. By the hypothesis, $P_j$ and $M_j$ are projective $A_0$-modules for all $j \geqslant 0$. Then $P_j \cong M_j \oplus (\Omega M)_j$, so all $(\Omega M)_j$ are projective $A_0$-modules for $j \geqslant 0$. Replacing $M$ by $\Omega M$ and using recursion, we conclude that all $\Omega^i(M)_j$ are projective $A_0$-modules for $i, j \geqslant 0$. In particular, all $\Omega^i(M)_i$ are projective $A_0$-modules.

(2) implies (3): Conversely, suppose that $\Omega^i(M)_i$ is a projective $A_0$-modules for every $i \geqslant 0$. We use contradiction to show that all $M_i$ are projective $A_0$-modules. If this not the case, we can find the minimal number $n \geqslant 0$ such that $M_n$ is not a projective $A_0$-module. As above, consider $\varphi_n: P_n \rightarrow M_n$ with kernel $(\Omega M)_n$. We claim that this kernel is not a projective $A_0$-module. Indeed, if it is a projective $A_0$-module, then it is injective as well, so $P_n \cong (\Omega M)_n \oplus M_n$. Consequently, $M_n$ is isomorphic to a summand of the projective $A_0$-module $P_n$ and must be a projective $A_0$-module, too. This is impossible. Therefore, $(\Omega M)_n$ is not a projective $A_0$-module. Now replacing $M$ by $\Omega M$ and using induction, we deduce that $\Omega^n (M)_n$ is not a projective $A_0$-module. This contradicts our assumption. Therefore, all $M_i$ are projective $A_0$-modules.
\end{proof}

Now we define \textit{quasi-Koszul modules} over the graded algebra $A$.

\begin{definition}
A non-negatively graded $A$-module $M$ is called quasi-Koszul if
\begin{equation*}
\Ext _{A}^1 (A_0, A_0) \cdot \Ext _{A}^i (M, A_0)=  \Ext _{A}^{i+1} (M, A_0)
\end{equation*}
for all $i \geq 0$. The algebra $A$ is called a quasi-Koszul algebra if $A_0$ as an $A$-module is quasi-Koszul.
\end{definition}

A graded $A$-module $M$ is quasi-Koszul if and only if as a graded $\Ext _{A}^{\ast} (A_0, A_0)$-module $\Ext _{A} ^{\ast} (M, A_0)$ is generated in degree $0$. The graded algebra $A$ is a quasi-Koszul algebra if and only if the cohomology ring $\Ext^{\ast}_{A} (A_0, A_0)$ is generated in degree 0 and degree 1.

The quasi-Koszul property is preserved by the Heller operator. Explicitly, if $M$ is a quasi-Koszul $A$-module with $M_0$ a projective $A_0$-module, then its syzygy $\Omega M$ is also quasi-Koszul. This is because for each $i \geqslant 1$, we have:
\begin{align*}
& \Ext _{A}^{i} (\Omega M, A_0) \cong \Ext _{A}^{i+1} (M, A_0) \\
& = \Ext _{A}^1 (A_0, A_0) \cdot \Ext _{A}^i (M, A_0)\\
& = \Ext _{A}^1 (A_0, A_0) \cdot \Ext _{A}^{i-1} (\Omega M, A_0).
\end{align*}
The identity $\Ext _{A}^i (M, A_0) \cong \Ext _{A}^{i-1} (\Omega M, A_0)$ is proved in Lemma 2.10.

If $A_0$ is a semisimple $k$-algebra, quasi-Koszul modules generated in degree 0 coincide with Koszul modules. This is not true if $A_0$ is only self-injective. Actually, by the following theorem, every Koszul module is quasi-Koszul, but the converse does not hold in general. For example, let $kG$ be the group algebra of a finite group concentrated in degree 0. The reader can check that every $kG$-module generated in degree 0 is quasi-Koszul, but only the projective $kG$-modules are Koszul. If $|G|$ is not invertible in $k$, then all non-projective $kG$-modules generated in degree 0 are quasi-Koszul but not Koszul.

The following theorem gives us a close relation between quasi-Koszul modules and Koszul modules.

\begin{theorem}
A graded $A$-module $M$ generated in degree 0 is Koszul if and only if it is quasi-Koszul and $\Omega^i (M)_i$ is a projective $A_0$-module for every $i \geqslant 0$.
\end{theorem}

The following lemma will be used in the proof of this theorem.

\begin{lemma}
Let $M$ be a graded $A$-module generated in degree 0 with $M_0$ a projective $A_0$-module. Then $\Omega M$ is generated in degree 1 if and only if every $A$-module homomorphism $\Omega M \rightarrow A_0$ extends to an $A$-module homomorphism $JP \rightarrow A_0$, where $P$ is a graded projective cover of $M$.
\end{lemma}

\begin{proof}
The short exact sequence $0 \rightarrow \Omega M \rightarrow P \rightarrow M \rightarrow 0$ induces an exact sequence $0 \rightarrow (\Omega M)_1 \rightarrow P_1 \rightarrow M_1 \rightarrow 0$. Applying the exact functor $\Hom _{A_0} (-, A_0)$ we get another exact sequence
\begin{equation*}
0 \rightarrow \Hom _{A_0} (M_1, A_0) \rightarrow \Hom _{A_0} (P_1, A_0) \rightarrow \Hom _{A_0} ((\Omega M)_1, A_0) \rightarrow 0.
\end{equation*}
Since $M_0$ is a projective $A_0$-module, $M_0 \cong P_0$, so $(\Omega M)_0 =0$. Therefore, $\Omega M$ is generated in degree 1 if and only if $\Omega M / J(\Omega M) \cong (\Omega M)_1$, if and only if the above sequence is isomorphic to
\begin{equation*}
0 \rightarrow \Hom _{A_0} (M_1, A_0) \rightarrow \Hom _{A_0} (P_1, A_0) \rightarrow \Hom _{A_0} (\Omega M / J \Omega M, A_0) \rightarrow 0.
\end{equation*}
Here we use the fact that $A_0$ is self-injective and the functor $\Hom _{A_0} (-, A_0)$ is a duality. But the above sequence is isomorphic to
\begin{equation*}
\xymatrix {0 \ar[r] & \Hom _A (JM, A_0) \ar[r] & \Hom _A (JP, A_0) \ar[r] & \Hom _A (\Omega M, A_0) \ar[r] & 0}
\end{equation*}
since $JM$ and $JP$ are generated in degree 1. Therefore, $\Omega M$ is generated in degree 1 if and only if every (non-graded) $A$-module homomorphism $\Omega M \rightarrow A_0$ extends to a (non-graded) $A$-module homomorphism $JP \rightarrow A_0$.
\end{proof}

Now let us prove the theorem.

\begin{proof}
\textbf{The only if part.} Let $M$ be a Koszul $A$-module. Without loss of generality we can suppose that $M$ is indecomposable. Notice that all syzygies $\Omega^i(M)[-i]$ are also Koszul. Therefore, $\Omega^i(M)_i \cong (\Omega^i(M) [-i])_0$ is a projective $A_0$-module for all $i \geqslant 0$.

Now we show that $M$ is quasi-Koszul, i.e.,
\begin{equation*}
\Ext_A^{i+1} (M, A_0) = \Ext_A^1 (A_0, A_0) \cdot \Ext _A^i (M, A_0)
\end{equation*}
for all $i \geqslant 0$. By Lemma 2.10, we have $\Ext _A^{i+1} (M, A_0) \cong \Ext_A^1 (\Omega^i (M), A_0)$ and $\Ext _A^i (M, A_0) \cong \Hom_A (\Omega^i(M), A_0)$. Therefore, it suffices to show $\Ext _{A}^1 (M, A_0) = \Ext _{A}^1 (A_0, A_0) \cdot \Hom_A (M, A_0)$ since the conclusion follows immediately if we replace $M$ by $(\Omega M) [-1]$ recursively.

To prove this identity, we first identify $\Ext _{A}^1 (M, A_0)$ with $\Hom_A (\Omega M, A_0)$ by Lemma 2.10. Take an element $x \in \Ext _{A}^1 (M, A_0)$ and let $g: \Omega M \rightarrow A_0$ be the corresponding homomorphism. Since $M$ is Koszul, $M_0$ is a projective $A_0$-module, and $\Omega M$ is generated in degree 1. Thus by the previous lemma, $g$ extends to $JP^0$, and hence there is a homomorphism $\tilde{g}: JP^0 \rightarrow A_0$ such that $g = \tilde{g} \iota$, where $P^0$ is a graded projective cover of $M$ and $\iota: \Omega M \rightarrow JP^0$ is the inclusion.
\begin{align*}
\xymatrix { \Omega M \ar[r]^{\iota} \ar[d]^g & JP^0 \ar[dl]^{\tilde{g}}\\
A_0 & }
\end{align*}

We have the following commutative diagram:
\begin{align*}
\xymatrix{ 0 \ar[r] & \Omega M \ar[d]^{\iota} \ar[r] & P^0 \ar@{=}[d] \ar[r] & M \ar[d]^p \ar[r] & 0\\
0 \ar[r] & JP^0 \ar[r] & P^0 \ar[r] & P^0_0 \ar[r] & 0
}
\end{align*}
where the map $p$ is defined to be the projection of $M$ onto $M_0 \cong P_0^0$.

The map $\tilde{g}: JP^0 \rightarrow A_0$ gives a push-out of the bottom sequence. Consequently, we have the following commutative diagram:
\begin{align*}
\xymatrix{ 0 \ar[r] & \Omega M \ar[d]^{\iota} \ar[r] & P^0 \ar@{=}[d] \ar[r] & M \ar[d]^p \ar[r] & 0\\
0 \ar[r] & JP^0 \ar[r] \ar[d] ^{\tilde{g}} & P^0 \ar[d] \ar[r] & P^0_0 \ar[r] \ar@{=}[d] & 0\\
0 \ar[r] & A_0 \ar[r] & E \ar[r] & P^0_0 \ar[r] & 0.}
\end{align*}

Since $P_0 \in \text{add} (A_0)$, we can find some $m$ such $P_0$ can be embedded into $A_0 ^{\oplus m}$. Thus the bottom sequence $y \in \Ext _A^1 (P_0, A_0) \subseteq \bigoplus _{i=1}^m \Ext _A^1 (A_0, A_0)$ and we can write $y = y_1 + \ldots + y_m$ where $y_i \in \Ext _A^1 (A_0, A_0)$ is represented by the sequence
\begin{equation*}
\xymatrix{ 0 \ar[r] & A_0 \ar[r] & E_i \ar[r] & A_0 \ar[r] & 0}.
\end{equation*}
Composed with the inclusion $\epsilon: P_0 \rightarrow A_0 ^{\oplus m}$, the map $\epsilon \circ p = (p_1, \ldots, p_m)$ where each component $p_i$ is defined in an obvious way. Consider the pull-backs:
\begin{equation*}
\xymatrix {0 \ar[r] & A_0 \ar[r] \ar@{=}[d] & F_i \ar[r] \ar[d] & M \ar[r] \ar[d]^{p_i} & 0 \\
0 \ar[r] & A_0 \ar[r] & E_i \ar[r] & A_0 \ar[r] & 0.}
\end{equation*}
Let $x_i$ be the top sequence. Then $x = \sum_{i=1}^m x_i = \sum_{i=1}^m y_i p_i \in \Ext _A^1 (A_0, A_0) \cdot \Hom_A (M, A_0)$ and hence $\Ext _A^1 (M, A_0) \subseteq \Ext _A^1 (A_0, A_0) \cdot \Hom_A (M, A_0)$. The other inclusion is obvious.\\

\textbf{The if part.} By Proposition 2.6, it suffices to show that $\Omega^i (M)$ is generated in degree $i$, $i>0$. But we observe that if $M$ is quasi-Koszul and $\Omega^i (M)_i$ are projective $A_0$-modules for all $i \geqslant 0$, then each $\Omega^i(M)$ has these properties as well. Thus we only need to show that $\Omega M$ is generated in degree 1 since the conclusion follows if we replace $M$ by $\Omega M$ recursively. By the previous lemma, it suffices to show that each (non-graded) $A$-module homomorphism $g: \Omega M \rightarrow A_0$ extends to $JP^0$.

The map $g$ gives a push-out $x \in \Ext _{A}^1 (M, A_0)$ as follows:
\begin{align*}
\xymatrix{ 0 \ar[r] & \Omega M \ar[d]^g \ar[r] & P^0 \ar[d] \ar[r] & M \ar@{=}[d] \ar[r] & 0\\
0 \ar[r] & A_0 \ar[r] & E \ar[r] & M \ar[r] & 0
}
\end{align*}
Since $M$ is quasi-Koszul, $x$ is contained in $\Ext _{A}^1 (A_0, A_0) \cdot \Hom_A (M, A_0)$. Thus $x = \sum_{i} y_i h_i$ with $y_i \in \Ext_{A}^1 (A_0, A_0)$ and $h_i \in \Hom_A (M, A_0)$, and each $y_i h_i$ gives the following commutative diagram, where the bottom sequence corresponds to $y_i$:
\begin{align}
\xymatrix{ 0 \ar[r] & A_0 \ar@{=}[d] \ar[r] & E_i \ar[d] \ar[r] & M \ar[d]^{h_i} \ar[r] & 0\\
0 \ar[r] & A_0 \ar[r] & F_i \ar[r] & A_0 \ar[r] & 0
}
\end{align}
By the natural isomorphism $\Ext _{A}^1 (M, A_0) \cong \Hom_A (\Omega M, A_0)$ (see Lemma 2.10), each $y_ih_i$ corresponds an $A$-homomorphism $g_i: \Omega M \rightarrow A_0$ such that the following diagram commutes:
\begin{align}
\xymatrix{ 0 \ar[r] & \Omega M \ar[d]^{g_i} \ar[r] & P^0 \ar[d] \ar[r] & M \ar@{=}[d] \ar[r] & 0\\
0 \ar[r] & A_0 \ar[r] & E_i \ar[r] & M \ar[r] & 0
}
\end{align}

Diagrams 2.2 and 2.3 give us:
\begin{align*}
\xymatrix{ 0 \ar[r] & \Omega M \ar@{=}[d] \ar[r]^{\iota} & JP^0 \ar[d]^{\tilde{j}} \ar[r] & JM \ar[d]^j \ar[r] & 0\\
0 \ar[r] & \Omega M \ar[r] \ar[d]^{g_i} & P^0 \ar[d]^{\tilde{h_i}} \ar[r] & M \ar[r] \ar[d]^{h_i} & 0\\
0 \ar[r] & A_0 \ar[r]^{\rho} & F_i \ar[r] & A_0 \ar[r] & 0
}
\end{align*}
Since $JM$ is sent to 0 by $h_ij$, there is a homomorphism $\varphi_i$ from $JP_0$ to the first term $A_0$ of the bottom sequence such that $\rho \varphi_i = \tilde{h_i} \tilde{j}$. Then $g_i$ factors through $\varphi_i$, i.e., $g_i = \varphi_i \iota$. Since $g = \sum_i g_i$, we know that $g$ extends to $JP^0$. This finishes the proof.
\end{proof}

An easy corollary of the above theorem is:

\begin{corollary}
Suppose that $A$ is projective viewed as an $A_0$-module. Then a graded $A$-module $M$ is Koszul if and only if it is quasi-Koszul as an $A$-module and projective as an $A_0$-module (or equivalently all $M_i$ are projective $A_0$-modules).
\end{corollary}

\begin{proof}
By Lemma 2.14, all $\Omega^i (M)_i$ are projective $A_0$-modules for $i \geqslant 0$ if and only if $M_s$ is a projective $A_0$-module for every $s \geqslant 0$. The conclusion follows from the previous theorem.
\end{proof}

In particular, if $A_0$ is is a Koszul $A$-module, then by letting $M = A$ in Proposition 2.13, we deduce that all $A_i$ are projective $A_0$-modules for $i \geqslant 0$.

\section{Generalized Koszul Algebras}
In this section we generalize to our context some useful results on classical Koszul algebras which appear in \cite{BGS}. As before, throughout this section $A$ is a non-negatively graded, locally finite associative $k$-algebra with $A_0$ self-injective. For two graded $A$-modules $M$ and $N$, we use $\Hom_A (M, N)$ and $\hom_A (M, N)$ to denote the space of all module homomorphisms and the space of graded module homomorphisms respectively. The derived functors Ext and ext correspond to Hom and hom respectively.

Recall that $A$ a \textit{quasi-Koszul algebra} if $A_0$ is quasi-Koszul as an $A$-module. In particular, if $A_0$ is a Koszul $A$-module, then $A$ is a quasi-Koszul algebra.

\begin{theorem}
The graded algebra $A$ is quasi-Koszul if and only if the opposite algebra $A^{\textnormal{op}}$ is quasi-Koszul.
\end{theorem}

\begin{proof}
Since the quasi-Koszul property is invariant under the Morita equivalence, without loss of generality we can suppose that $A$ is a basic algebra. Therefore, $A_0$ is also a basic algebra. Let $M$ and $N$ be two graded $A$-modules. We claim $\ext^i_{A}(M, N) \cong \ext _{A^{\textnormal{op}}}^i (DN, DM)$ for all $i \geqslant 0$, where $D$ is the graded duality functor $\hom_k (-,k)$. Indeed, Let
\begin{align*}
\xymatrix{ \ldots \ar[r] & P^2 \ar[r] & P^1 \ar[r] & P^0 \ar[r] & M \ar[r] & 0 }
\end{align*}
be a projective resolution of $M$. Applying the graded functor $\hom_A (-,N)$ we get the following chain complex $C^{\ast}$:
\begin{align*}
\xymatrix{ 0 \ar[r] & \hom_A (P^0, N) \ar[r] & \hom_A (P^1, N) \ar[r] & \ldots. }
\end{align*}
Using the natural isomorphism $\hom_A (P^i, N) \cong \hom _{A^{\textnormal{op}}} (DN, DP^i)$, we get another chain complex $E^{\ast}$ isomorphic to the above one:
\begin{align*}
\xymatrix{ 0 \ar[r] & \hom_{A^{\textnormal{op}}} (DN, DP^0) \ar[r] & \hom_{A^{\textnormal{op}}} (DN, DP^1) \ar[r] & \ldots. }
\end{align*}
Notice that all $DP^i$ are graded injective $A^{\textnormal{op}}$-modules. Thus
\begin{equation*}
\ext ^i_{A}(M, N) \cong H^i(C^{\ast}) \cong H^i(E^{\ast}) \cong \ext _{A^{\textnormal{op}}}^i (DN, DM)
\end{equation*}
which is exactly our claim.

Now let $M = N = A_0$. Then $\ext ^i_{A} (A_0, A_0) \cong \ext _{A^{\textnormal{op}}}^i (DA_0, DA_0)$. Since $A_0$ is self-injective and basic, it is a Frobenius algebra. Therefore, $DA_0$ is isomorphic to $A_0^{\textnormal{op}}$ as a left $A_0^{\textnormal{op}}$-module (and hence as a left $A^{\textnormal{op}}$-module). By the next proposition, $A_0$ is a quasi-Koszul $A$-module if and only if $A_0^{\textnormal{op}}$ is a quasi-Koszul $A^{\textnormal{op}}$-module.
\end{proof}

However, if $A_0$ is a Koszul $A$-module, $A^{\textnormal{op}}_0$ need not be a Koszul $A^{\textnormal{op}}$-module, as shown by the following example.

\begin{example}
Let $\mathcal{E}$ be a finite EI category with two objects $x$ and $y$ such that: $\Aut _{\mathcal {C}} (x) = \langle g \rangle$ is a cyclic group of order 2; $\Aut _{\mathcal {C}} (y)$ is a trivial group; $\Hom _{\mathcal{E}} (x, y) = \{ \alpha \}$ (thus $\alpha \circ g = \alpha$) and $\Hom _{\mathcal{E}} (y, x) = \emptyset$. Let $k$ be an algebraically closed field of characteristic 2. The category algebra $A = k \mathcal{E}$ is of dimension 4. Let $A_0$ be the space spanned by $1_x, g$ and $1_y$, and let $A_1$ be the one-dimensional space spanned by $\alpha$. The reader can check that $A_0$ is a Koszul $A$-module. The opposite algebra $A^{\textnormal{op}} = k \mathcal{E} ^{\textnormal{op}}$ can also be graded in a similar way, but $A^{\textnormal{op}}_0$ is not a Koszul $A^{\textnormal{op}}$-module. However, we will show in Section 6 that $A^{\textnormal{op}}$ is a quasi-Koszul algebra.
\end{example}

\begin{proposition}
The graded algebra $A$ is Koszul if and only if $\Omega^i (A_0)_i$ are projective $A_0$-modules for all $i \geqslant 0$, and whenever $\ext^i _{A} (A_0, A_0 [n]) \neq 0$ we have $n=i$.
\end{proposition}

\begin{proof}
If $A$ is a Koszul algebra, then $A_0$ is a Koszul $A$-module, and by Corollary 2.5 all $\Omega^i (A_0)_i$ are projective $A_0$-modules. Moreover, there is a linear projective resolution
\begin{equation*}
\xymatrix { \ldots \ar[r] & P^2 \ar[r] & P^1 \ar[r] & P^0 \ar[r] & A_0 \ar[r] & 0}
\end{equation*}
where $P^i$ is generated in degree $i$. Applying $\hom_A (-, A_0 [n])$ we find that all terms in this complex except $\hom_A (P^{n}, A_0 [n])$ are 0. Consequently, $\ext^i _{A} (A_0, A [n]) = 0$ unless $i=n$.

Conversely, suppose that $\Omega^i (A_0)_i$ is a projective $A_0$-module for each $i \geqslant 0$ and $\ext^i _{A} (A_0, A_0 [n]) = 0$ unless $n=i$, we want to show that $\Omega^{i} (A_0)$ is generated in degree $i$ by induction. Obviously, $\Omega^0 (A_0) = A_0$ is generated in degree 0. Suppose that $\Omega^j(A_0)$ is generated in degree $j$ for $0 \leqslant j \leqslant i$. Now consider $\Omega^{i+1} (A_0)$. By applying the graded version of Lemma 2.10 recursively, we have
\begin{equation*}
\hom_A (\Omega^{i+1} (A_0), A_0 [n]) = \ext _{A}^{i+1} (A_0, A_0 [n]).
\end{equation*}
The right-hand side is 0 unless $n = i + 1$, so $\Omega^{i+1} (A_0)$ is generated in degree $i+1$. By induction we are done.
\end{proof}

The reader can check that the conclusion of this proposition is also true for Koszul modules. i.e., $M$ is a Koszul $A$-module if and only if $\ext _A ^i (M, A_0[n]) \neq 0$ implies $n = i$.

We can define a tensor algebra $T(A)$ generated by $A_1$, which is a $(A_0, A_0)$-bimodule. Explicitly,
\begin{equation*}
T(A) = A_0 \oplus A_1 \oplus (A_1 \otimes A_1) \oplus (A_1 \otimes A_1 \otimes A_1) \oplus \ldots,
\end{equation*}
where all tensors are over $A_0$ and we use $\otimes$ rather than $\otimes _{A_0}$ to simplify the notation. This tensor algebra has a natural grading. Clearly, $A$ is a quotient algebra of $T(A)$. Let $I$ be the kernel of the quotient map $q: T(A) \rightarrow A$. We say that $A$ is a \textit{quadratic algebra} if the ideal $I$ has a set of generators contained in $A_1 \otimes A_1$.

\begin{theorem}
If $A$ is a Koszul algebra, then it is a quadratic algebra.
\end{theorem}

\begin{proof}
This proof is a modification of the proofs of Theorem 2.3.2 and Corollary 2.3.3 in \cite{BGS}. First, consider the exact sequence
\begin{equation*}
\xymatrix{ 0 \ar[r] & W \ar[r] & A \otimes A_1 \ar[r] & A \ar[r] & A_0 \ar[r] & 0}
\end{equation*}
where $W$ is the kernel of the multiplication. Clearly, $\Omega(A_0) \cong J = \bigoplus _{i \geqslant 1} A_i$. Since the image of $(A \otimes A_1)_1 = A_0 \otimes A_1$ under the multiplication is exactly $A_1 = \Omega(A_0)_1$, $A \otimes A_1$ is a projective cover of $\Omega(A_0)$ and $\Omega^2 (A_0) = W \subseteq J \otimes A_1$. Therefore, $W$ is generated in degree 2, and hence $W/JW \cong W_2$ is concentrated in degree 2. Observe that $A$ is a quotient algebra of $T(A)$ with kernel $I$. Let $R_n$ be the kernel of the quotient map $A_1 ^{\otimes n} \rightarrow A_n$.

If $A$ is not quadratic, we can find some $x \in R_n$ with $n > 2$ such that $x$ is not contained in the two-sided ideal generated by $\sum _{i=2}^{n-1} R_i$. Consider the following composite of maps:
\begin{equation*}
\xymatrix {A_1 ^{\otimes n} = A_1 ^{\otimes n-1} \otimes A_1 \ar[r] ^-p & A_{n-1} \otimes A_1 \ar[r]^-m & A_n}.
\end{equation*}
Clearly $p(x) \in W$ since $m(p(x)) =0$. We show $p(x) \notin JW$ by contradiction.

Indeed, if $p(x) \in JW$, then $p(x) \in A_1 W$ since $JW \cong \bigoplus _{i \geqslant 3} W_i = A_1 W$ (notice that $W$ is generated in degree 2). Therefore, we can express $p(x)$ as a linear combination of vectors of the form $\lambda \cdot w$ with $\lambda \in A_1$ and $w \in W$. But $W \subseteq J \otimes A_1$, so each $w$ can be expressed as $\sum_i w_i' \otimes \lambda_i'$ with $w_i' \in A_{n-2}$, $\lambda_i' \in A_1$ such that $\sum_{i=1}^s w'_i \cdot \lambda'_i =0$ by the definition of $W$.

Since there is a surjective product map $\varphi: A_1 ^{\otimes n-2} \twoheadrightarrow A_{n-2}$, we can choose a pre-image $v_i^1 \otimes \ldots \otimes v_i^{n-2} \in \varphi^{-1} (w_i')$ for each $i$ and define
\begin{equation*}
\tilde{w} = \sum_{i=1}^s v_i^1 \otimes \ldots \otimes v_i^{n-2} \otimes \lambda'_i
\end{equation*}
which is contained in $R_{n-1}$ clearly. Observe that $p(\lambda \otimes \tilde{w}) = \lambda \cdot w$. Since $p(x)$ is a linear combination of vectors of the form $\lambda \cdot w$, by the above process we can get some $y$ which is a linear combination of vectors of the form $\lambda \otimes \tilde{w}$ such that $p(y) = p(x)$. Clearly, $p(x -y) = 0$ and $y \in A_1 \otimes R_{n-1}$.

Consider the following short exact sequence
\begin{equation*}
\xymatrix{ 0 \ar[r] & R_{n-1} \ar[r] & A_1 ^{\otimes n-1} \ar[r] & A_{n-1} \ar[r] & 0.}
\end{equation*}
Since $A_0$ is a Koszul $A$-module, so is $J[-1] = \Omega(A_0)[-1]$. Therefore, by Corollary 2.5, $A_1 \cong J[-1]_0$ is a projective $A_0$-module. Therefore, the following sequence is also exact:
\begin{equation*}
\xymatrix{ 0 \ar[r] & R_{n-1} \otimes A_1 \ar[r] & A_1 ^{\otimes n} \ar[r]^-p & A_{n-1} \otimes A_1 \ar[r] & 0.}
\end{equation*}
Thus $x - y \in R_{n-1} \otimes A_1$ since $p(x-y) =0$. It follows $x \in A_1 \otimes R_{n-1} + R_{n-1} \otimes A_1$, which contradicts our choice of $x$.

We proved $x \notin JW$. Then $p(x) \in W/JW \cong W_2$ is of degree 2. But this is impossible since $p$ as a graded homomorphism sends $x \in R_n$ with $n >2$ to an element of degree $n$.
\end{proof}

We can define the \textit{Koszul complex} for $A$ in a similar way to the classical situation.

Let $A \cong T(A) /(R)$ be quadratic where $R \subseteq A_1 \otimes A_1$ is a set of relations. Define $P_n^n = \bigcap _{i=0} ^{n-2} A_1 ^{\otimes i} \otimes R \otimes A_1 ^{\otimes n-i-2} \subseteq A_1 ^{\otimes n}$. In particular, $P_0^0 = A_0$, $P_1^1 = A_1$ and $P_2^2 = R$. Let $P^n  = A \otimes P_n^n$ such that $A_0 \otimes P_n^n \cong P_n^n$ is in degree $n$. Define $d^n: P^n \rightarrow P^{n-1}$ to be the restriction of $A \otimes A_1 ^{\otimes n} \rightarrow A \otimes A_1 ^{\otimes n-1}$ by $a \otimes v_1 \otimes \ldots \otimes v_n \mapsto av_1 \otimes v_2 \otimes \ldots \otimes v_n$. The reader can check $d^{n-1} d^{n} =0$ for $n \geqslant 1$. Therefore we get the following Koszul complex $K^{\ast}$:
\begin{equation*}
\xymatrix {\ldots \ar[r] & P^3 \ar[r]^{d^3} & A \otimes R \ar[r]^{d^2} & A \otimes A_1 \ar[r]^{d^1} & A \ar[r] & 0.}
\end{equation*}

\begin{theorem}
Let $A \cong T(A) /(R)$ be a quadratic algebra. Then $A$ is a Koszul algebra if and only if the Koszul complex is a projective resolution of $A_0$.
\end{theorem}

\begin{proof}
One direction is trivial. Now suppose that $A_0$ is a Koszul $A$-module. The Koszul complex $K^{\ast}$ of $A$ has the following properties:

\noindent (1). Let $Z^n$ be the kernel of $d^n: P^n \rightarrow P^{n-1}$. The restricted map $d_n^n:$
\begin{equation*}
P^n_n = A_0 \otimes \big{(}\bigcap _{i=0} ^{n-2} A_1 ^{\otimes i} \otimes R \otimes A_1 ^{\otimes n-i-2} \big{)} \rightarrow P^{n-1} _n = A_1 \otimes (\bigcap _{i=0} ^{n-3} A_1 ^{\otimes i} \otimes R \otimes A_1 ^{\otimes n-i-3})
\end{equation*}
is injective. Therefore $Z^n_i = 0$ for every $i \leqslant n$.

\noindent (2). $Z_{n+1}^n$, the kernel of the map $d_n^{n+1}:$
\begin{equation*}
P^n _{n+1} = A_1 \otimes \big{(} \bigcap_{i=0} ^{n-2} A_1^{\otimes i} \otimes R \otimes A_1 ^{\otimes n-i -2} \big{)} \rightarrow P^{n-1} _{n+1} = A_2 \otimes \big{(} \bigcap_{i=0} ^{n-3} A_1^{\otimes i} \otimes R \otimes A_1 ^{\otimes n-i -3} \big{)}
\end{equation*}
is
\begin{equation*}
A_1 \otimes \big{(} \bigcap_{i=0} ^{n-2} A_1^{\otimes i} \otimes R \otimes A_1 ^{\otimes n-i -2} \big{)} \cap (R \otimes A_1 ^{\otimes n-1}) = \bigcap_{i=0} ^{n-1} A_1^{\otimes i} \otimes R \otimes A_1 ^{\otimes n-i -1}
\end{equation*}
which is exactly $P_{n+1} ^{n+1}$ (or $d^{n+1}_{n+1} (P_{n+1} ^{n+1})$ since $d^{n+1}_{n+1}$ is injective by the last property).

We claim that each $P^n = A \otimes P^n_n$ is a projective $A$-module. Clearly, it is enough to show that each $P_n^n = Z_n ^{n-1}$ is a projective $A_0$-module. We prove the following stronger conclusion. That is, $Z^n_i$ are projective $A_0$-modules for $i \in \mathbb{Z}$ and $n \geqslant 0$. We use induction on $n$.

Since $A_0$ is a Koszul $A$-module, by Proposition 2.13, $A_i$ are projective $A_0$-modules for all $i \geqslant 0$. The conclusion is true for $Z^0 \cong J$ since $J_0 = 0$ and $J_m = A_m$ for $m \geqslant 1$. Suppose that it is true for $l \leqslant n$. That is, all $Z^l _i$ are projective $A_0$-modules for $l \leqslant n$ and $i \in \mathbb{Z}$. Consider $l = n+1$. By the second property described above, $P_{n+1} ^{n+1} = Z_{n+1}^n$, which is a projective $A_0$-module by the induction hypothesis. Therefore, $P^{n+1} = A \otimes P^{n+1}_{n+1}$ is a projective $A$-module, so $P^{n+1} _i$ are all projective $A_0$-modules for $i \in \mathbb{Z}$. But the following short exact sequence of $A_0$-modules splits
\begin{equation*}
\xymatrix{ 0 \ar[r] & Z^{n+1} _i \ar[r] & P^{n+1} _i \ar[r] & Z^n_i \ar[r] & 0}
\end{equation*}
since $Z^n_i$ is a projective $A_0$-module by the induction hypothesis. Now as a direct summand of $P^{n+1}_i$ which is a projective $A_0$-module, $Z^{n+1} _i$ is a projective $A_0$-module as well. Our claim is proved by induction.

We claim that this complex is acyclic. First, the sequence
\begin{equation*}
\xymatrix{ P^1 = A \otimes A_1 \ar[r] & P^0 = A \otimes A_0 \ar[r] & A_0 \ar[r] & 0}
\end{equation*}
is right exact. By induction on $n>1$
\begin{equation*}
\ext _{A} ^{n+1} (A_0, A_0 [m]) = \text{coker} \big{(} \hom_A (P^n, A_0 [m]) \rightarrow \hom_A (Z^n, A_0 [m]) \big{)}.
\end{equation*}
By Property (1), $Z^n_m = 0$ if $m < n+1$. If $m>n+1$, $\hom_A (P^n, A_0 [m]) =0$ since $P^n$ is generated in degree $n$, so $\ext _{A} ^{n+1} (A_0, A_0 [m]) = \hom_A (Z^n, A_0 [m])$ by the above identity. But the left-hand side of this identity is non-zero only if $m = n+1$ since $A_0$ is Koszul. Therefore, $\hom_A (Z^n, A_0 [m]) = 0$ for $m > n+1$. Consequently, $Z^n$ is generated in degree $n+1$. By property (2), $Z^n _{n+1} = d^{n+1}_{n+1} (P^{n+1} _{n+1})$, so $Z^n = d^{n+1} (P^{n+1})$ since both modules are generated in degree $n+1$. Therefore, the Koszul complex is acyclic, and hence is a projective resolution of $A_0$.
\end{proof}

\section{Generalized Koszul Duality}
In this section we prove the Koszul duality. As before, $A$ is a non-negatively graded, locally finite algebra with $A_0$ self-injective. Define $\Gamma = \Ext ^{\ast}_{A} (A_0, A_0)$ which has a natural grading. Notice that $\Gamma_0 \cong A_0^{\textnormal{op}}$ is also a self-injective algebra. Let $M$ be a graded $A$-module. Then $\Ext ^{\ast}_{A} (M, A_0)$ is a graded $\Gamma$-module. Moreover, if $M$ is Koszul, then $\Ext _A^{\ast} (M, A_0)$ is generated in degree 0, so it is a finitely generated $\Gamma$-module. Thus $E = \Ext ^{\ast}_{A} (-, A_0)$ gives rise to a functor from the category of Koszul $A$-modules to $\Gamma$-gmod.

\begin{theorem}
If $A$ is a Koszul algebra, then $E = \Ext ^{\ast}_{A} (-, A_0)$ gives a duality between the category of Koszul $A$-modules and the category of Koszul $\Gamma$-modules. That is, if $M$ is a Koszul $A$-module, then $E(M)$ is a Koszul $\Gamma$-module, and $E_{\Gamma}EM = \Ext ^{\ast} _{\Gamma} (EM, \Gamma_0) \cong M$ as graded $A$-modules.
\end{theorem}

\begin{proof}
Since $M$ and $A_0$ both are Koszul, by Proposition 2.13 $M_0$ and $JM[-1]$ are Koszul, where $J= \bigoplus _{i \geqslant 1} A_i$. Furthermore, we have the following short exact sequence of Koszul modules:
\begin{align*}
\xymatrix{ 0 \ar[r] & \Omega M [-1] \ar[r] & \Omega (M_0)[-1] \ar[r] & JM[-1] \ar[r] & 0.}
\end{align*}
As in the proof of Proposition 2.9, this sequence induces exact sequences recursively (see diagram 2.1):
\begin{align*}
\xymatrix{ 0 \ar[r] & \Omega^i(M)[-i] \ar[r] & \Omega^i(M_0)[-i] \ar[r] & \Omega^{i-1}(JM[-1]) [1-i] \ar[r] & 0,}
\end{align*}
and gives exact sequences of $A_0$-modules:
\begin{align*}
\xymatrix{ 0 \ar[r] & \Omega^i(M)_i \ar[r] & \Omega^i(M_0)_i \ar[r] & \Omega^{i-1}(JM [-1])_{i-1} \ar[r] & 0.}
\end{align*}
Applying the exact functor $\Hom_{A_0} (-, A_0)$ and using the following isomorphism for a graded $A$-module $N$ generated in degree $i$
\begin{equation*}
\Hom_A (N, A_0) \cong \Hom_A (N_i, A_0) \cong \Hom _{A_0} (N_i, A_0,)
\end{equation*}
we get:
\begin{align*}
0 \rightarrow \Hom_A (\Omega^{i-1}(JM[-1]), A_0) \rightarrow \Hom_A (\Omega^i(M_0), A_0) \rightarrow \Hom_A (\Omega^iM, A_0) \rightarrow 0.
\end{align*}
By Lemma 2.10, this sequence is isomorphic to
\begin{align*}
0 \rightarrow \Ext^{i-1} _{A} (JM[-1], A_0) \rightarrow \Ext^i _A (M_0, A_0) \rightarrow \Ext^i _A (M, A_0) \rightarrow 0.
\end{align*}
Now let the index $i$ vary and put these sequences together. We have:
\begin{align*}
\xymatrix{ 0 \ar[r] & E(JM[-1])[1] \ar[r] & E(M_0) \ar[r] & EM \ar[r] & 0.}
\end{align*}

Let us focus on this sequence. We claim $\Omega(EM) \cong E(JM[-1])[1]$. Indeed, since $M_0$ is a projective $A_0$-module and the functor $E$ is additive, $E(M_0)$ is a projective $\Gamma$-module. Since $JM[-1]$ is Koszul, $JM[-1]$ is quasi-Koszul and hence $E(JM[-1])$ as a $\Gamma$-module is generated in degree 0. Thus $E(JM[-1])[1]$ is generated in degree 1, and $E(M_0)$ is minimal. This proves the claim. Consequently, $\Omega(EM)$ is generated in degree 1 as a $\Gamma$-module. Moreover, replacing $M$ by $JM[-1]$ (which is also Koszul) and using the claimed identity, we have
\begin{equation*}
\Omega^2(EM) = \Omega(E(JM[-1])[1]) = \Omega( E(JM[-1]) [1] = E(J^2M[-2])[2],
\end{equation*}
which is generated in degree 2. By recursion, we know that $\Omega^i (EM) \cong E(J^iM [-i])[i]$ is generated in degree $i$ for all $i \geqslant 0$. Thus $EM$ is a Koszul $\Gamma$-module. In particular for $M = A$,
\begin{equation*}
EA = \Ext _{A} ^{\ast} (A, A_0) = \Hom_A (A, A_0) = \Gamma_0
\end{equation*}
is a Koszul $\Gamma$-module.

Since $\Omega ^i (EM)$ is generated in degree $i$ and
\begin{align*}
\Omega^i (EM)_i & \cong E(J^iM [-i])[i]_i \cong E(J^iM[-i])_0 \\
& = \Hom_A (J^iM[-i], A_0) \cong \Hom_A (M_i, A_0),
\end{align*}
we have
\begin{align*}
\Hom_{\Gamma} (\Omega^i(EM), \Gamma_0) & \cong \Hom _{\Gamma _0} (\Omega^i(EM)_i, \Gamma_0) \nonumber \\
& \cong \Hom _{\Gamma _0} (\Hom_{A} (M_i, A_0), \Gamma_0) \nonumber \\
& \cong \Hom _{\Gamma_0} (\Hom_{A_0} (M_i, A_0), \Gamma_0) \nonumber \\
& \cong M_i.
\end{align*}
The last isomorphism holds because $A_0$ is self-injective and $\Gamma_0 \cong A_0^{\textnormal{op}}$.

We have proved that $EM$ is a Koszul $\Gamma$-module. Therefore, $(\Omega ^i (EM))_i$ is a projective $\Gamma_0$-module for every $i \geqslant 0$. Applying Lemma 2.10 recursively, $\Ext _{\Gamma}^i (EM, \Gamma_0) \cong \Hom_{\Gamma} (\Omega^i(EM), \Gamma_0) \cong M_i$ for every $i \geqslant 0$. Adding them together, $E_{\Gamma}E(M) \cong \bigoplus_{i=0}^{\infty} M_i \cong M$.

Now we have $E_{\Gamma} ( E(A)) = E_{\Gamma} (\Gamma_0) \cong A$. Moreover, $\Gamma$ is a graded algebra such that $\Gamma_0 \cong A_0^{\textnormal{op}}$ is self-injective as an algebra and Koszul as a $\Gamma$-module. Using this duality, we can exchange $A$ and $\Gamma$ in the above reasoning and get $EE_{\Gamma}(N) \cong N$ for an arbitrary Koszul $\Gamma$-module $N$. Thus $E$ is a dense functor.

Let $L$ be another Koszul $A$-module. Since $L, M, EL, EM$ are all generated in degree 0, we have
\begin{align*}
\hom _{\Gamma} (EL, EM) & \cong \Hom_{\Gamma_0} ((EL)_0, (EM)_0)\\
& \cong \Hom_{\Gamma_0} (\Hom_A (L, A_0), \Hom_A (M, A_0))\\
& \cong \Hom_{A_0 ^{\textnormal{op}}} (\Hom_{A_0} (L_0, A_0), \Hom_{A_0} (M_0, A_0))\\
& \cong \Hom_{A_0} (L_0, M_0) \cong \hom_A (L, M).
\end{align*}
Consequently, $E$ is a duality between the category of Koszul $A$-modules and the category of Koszul $\Gamma$-modules.
\end{proof}

\begin{remark}
We can also use $\hom_A (-, A_0)$ to define the functor $E$ on the category of Koszul $A$-modules, namely $E := \bigoplus _{i \geqslant 0} \ext _A^i (-, A_0[i])$. Indeed, Let $M$ be a Koszul $A$-module. Note that all syzygies of $M$ are finitely generated. Thus
\begin{align*}
\Ext _A^{\ast} (M, A_0) & = \bigoplus _{i \geqslant 0} \Ext _A^i (M, A_0) \\
& = \bigoplus _{i \geqslant 0} \bigoplus _{j \in \mathbb{Z}} \ext _A^i (M, A_0 [j]) \\
& = \bigoplus _{i \geqslant 0} \ext _A^i (M, A_0[i])
\end{align*}
since $\ext _A^i (M, A_0 [j]) = 0$ for $i \neq j$, see the paragraph after Proposition 3.3.
\end{remark}

\section{Application to Directed Categories}
In this section we will apply the general theory to a type of structures called \textit{directed categories}, for which there exist very nice relations between stratification theory and Koszul theory and a nice correspondence between our generalized Koszul theory and the classical theory. All categories $\mathcal{C}$ we consider in this section are locally finite $k$-linear categories with finitely many objects, that is, for $x, y \in \Ob \mathcal{C}$, the set of morphisms $\mathcal{C} (x,y)$ is a finite-dimensional $k$-vector space. To simplify the technical part, we suppose furthermore that $\mathcal{C}$ is \textit{skeletal}, i.e., $x \cong y$ implies $x = y$ for $x,y \in \Ob \mathcal{C}$.

\begin{definition}
A locally finite $k$-linear category $\mathcal{C}$ is a directed category if there is a partial order $\leqslant$ on $\Ob \mathcal{C}$ such that $\mathcal{C} (x, y) \neq 0$ only if $x \leqslant y$.
\end{definition}

Correspondingly, we define \textit{directed algebras}.

\begin{definition}
A finite-dimensional algebra $A$ is called a directed algebra with respect to a partially ordered set of orthogonal idempotents $\{ e_i; \leqslant \} _{i=1}^n$ if $\sum _{i=1}^n  e_i = 1$ and $\Hom_A (A e_i, A e_j) \cong e_i A e_j \neq 0$ implies $e_j \leqslant e_i$.
\end{definition}

Notice that in the above definition we do not require the idempotents $e_i$ to be primitive. Clearly, every algebra $A$ is directed with respect to the trivial set $\{1\}$.

There is a bijective correspondence between directed categories and directed algebras. Let $A$ be a directed algebra with respect to a poset of orthogonal idempotents $(\{e_i\} _{i=1}^n, \leqslant)$. Then we can construct a directed category $\mathcal{A}$ in the following way: $\Ob \mathcal{A} = \{e_i\} _{i=1}^n$ with the same partial order, and $\mathcal{A} (e_i, e_j) = e_j A e_i \cong \Hom_A (A e_j, A e_i)$. The reader can check that $\mathcal{A}$ is indeed a directed category. We call $\mathcal{A}$ the \textit{associated category} of $A$.

Conversely, given a directed category $\mathcal{A}$ with the poset $(\Ob \mathcal{A}, \leqslant)$, we obtain an algebra $A$ which is directed with respect to the poset of orthogonal idempotents $( \{ 1_x \} _{x \in \Ob \mathcal{A}}, \leqslant)$, namely, $1_x \leqslant 1_y$ if and only if $x \leqslant y$. As a $k$-vector space, $A = \bigoplus _{x, y \in \Ob \mathcal{A}} \mathcal{A} (x, y)$. For two morphisms $\alpha: x \rightarrow y$ and $\beta: z \rightarrow w$, the product $\beta \cdot \alpha = 0$ if $y \neq z$, otherwise it is the composite morphism $\beta \alpha$. Since every vector in $A$ is a linear combination of morphisms in $\mathcal{A}$, the multiplication of morphisms can be extended linearly to a well defined product in $A$. The reader can check that the algebra $A$ we get in this way is indeed a directed algebra, which is called the \textit{associated algebra} of $\mathcal{A}$.

It is well known that $A$-mod is equivalent to the category of finite-dimensional $k$-linear representations of $\mathcal{A}$. If one of $A$ and $\mathcal{A}$ is graded, then the other one can be graded as well. Moreover, $A$-gmod is equivalent to the category of finite-dimensional graded $k$-linear representations of $\mathcal{A}$. For more details, see \cite{Mazorchuk1}. Because of these facts, we may view a directed category $\mathcal{A}$ as a directed algebra with respect to the set of idempotents $\{ 1_x \mid x \in \Ob \mathcal{A} \}$ and abuse notation and terminologies. For example, we may say idempotents in $\mathcal{A}$, ideals of $\mathcal{A}$ and so on. We hope this would not cause confusions to the reader and point out that all results in previous sections can be applied to directed categories.

Directed categories generalize $k$-linearizations of finite EI categories. Explicitly, let $\mathcal{E}$ be a skeletal finite EI category. Consider the category algebra $k\mathcal{E}$ with a set of idempotents $\{ 1_x \} _{x \in \Ob \mathcal{E}}$ on which there is a partial order $\leq$ such that $1_x \leq 1_y$ if and only if $\mathcal{E} (x, y) \neq \emptyset$. Then the category algebra $k\mathcal{E}$ is directed with respect to this poset of idempotents, so we can construct a direct category $\tilde {\mathcal{E}}$ by the above correspondence. Actually, $\tilde{ \mathcal{E}}$ is precisely the $k$-linearization of $\mathcal{E}$. In \cite{Li2, Li3}, we proved that for a finite-dimensional algebra $R$ standardly stratified with respect to a partial order, the associated category of the extension algebra of standard modules is a directed category. Moreover, if $R$ is standardly stratified for all linear orders, then the associated category of $R$ is a directed category as well.

Let $\mathcal{C}$ be a directed category. A \textit{$\mathcal{C}$-module} (or a \textit{representation} of $\mathcal{C}$) is defined to be a $k$-linear functor from $\mathcal{C}$ to the category of finite-dimensional $k$-vector spaces. The morphism space of $\mathcal{C}$ can be decomposed as the direct sum of $\mathcal{C} 1_x$ with $x$ ranging over all objects, where by $\mathcal{C} 1_x$ we denote the vector space formed by all morphisms with source $x$. Therefore, each $\mathcal{C} 1_x$ is a projective $\mathcal{C}$-module, and every indecomposable projective $\mathcal{C}$-module is isomorphic to a direct summand of a certain $\mathcal{C} 1_x$. The isomorphism classes of simple $\mathcal{C} (x, x)$-modules with $x$ varying within $\Ob \mathcal{C}$ give rise to isomorphism classes of simple $\mathcal{C}$-modules. Explicitly, let $V$ be a simple $\mathcal{C} (x, x)$-module for some object $x$, we can construct a simple $\mathcal{C}$-module $S$: $S(x) = V$ and $S(y) = 0$ for $y \neq x$. These results are well known for finite EI categories, see \cite{Webb}.

Our next task is to translate some results on finite EI categories in Section 2 of \cite{Webb} to directed categories. First, we want to show that every directed category is stratified with respect to the given partial order.

\begin{proposition}
Let $\mathcal{D}$ and $\mathcal{E}$ be full subcategories of a directed category $\mathcal{C}$ such that $\Ob \mathcal{D} \cup \Ob \mathcal{E} = \Ob \mathcal{C}$, $\Ob \mathcal{D} \cap \Ob \mathcal{E} = \emptyset$, and $\mathcal{C} (x,y) = 0$ for $x \in \Ob \mathcal{D}$ and $y \in \Ob \mathcal{E}$. Let $e = \sum _{x \in \Ob \mathcal{D}} 1_x$ and $I = \mathcal{C} e \mathcal{C}$. Then $I$ is a stratified ideal of $\mathcal{C}$.
\end{proposition}

\begin{proof}
The proof is similar to that of Proposition 2.2 in \cite{Webb}. Clearly $I$ is idempotent. Notice that $\mathcal{C} e$ is the space spanned by all morphisms with sources contained in $\Ob \mathcal{D}$ and $e \mathcal{C} e$ is the space spanned by all morphisms with both sources and targets contained in $\Ob \mathcal{D}$. Since $\mathcal{C} (x, y) = 0$ for $x \in \Ob \mathcal{D}$ and $y \in \Ob \mathcal{E}$, these two spaces coincide, i.e., $\mathcal{C} e = e \mathcal{C} e$. In particular, $\mathcal{C} e$ is projective $e \mathcal{C} e$-module, here $e \mathcal{C} e$ is an algebra for which the associated directed category is precisely $\mathcal{D}$. Therefore, Tor$ _n ^{e \mathcal{C} e} (\mathcal{C}e, e\mathcal{C} ) = 0$ for $n \geqslant 1$. Furthermore,
\begin{equation*}
\mathcal{C} e \otimes _{e \mathcal{C} e} e \mathcal{C} = e \mathcal{C} e \otimes _{e \mathcal{C} e} e \mathcal{C} \cong e \mathcal{C}.
\end{equation*}
We claim $e \mathcal{C} = \mathcal{C} e \mathcal{C}$. Clearly, $e \mathcal{C} \subseteq \mathcal{C} e \mathcal{C}$. On the other hand, since we just proved $\mathcal{C} e = e \mathcal{C} e$, we have $\mathcal{C} e \mathcal{C} = e \mathcal{C} e \mathcal{C} \subseteq e \mathcal{C}$. Therefore, $e \mathcal{C} = \mathcal{C} e \mathcal{C}$ as we claimed. In conclusion, $I$ is indeed a stratified ideal of $\mathcal{C}$.
\end{proof}

\begin{corollary}
Every directed category $\mathcal{C}$ is stratified with respect to the give partial order on $\Ob \mathcal{C}$.
\end{corollary}

\begin{proof}
The partial order $\leqslant$ on $\Ob \mathcal{C}$ gives a filtration on $\Ob \mathcal{C}$ in the following way: let $S_1$ be a set containing a maximal object in $\Ob \mathcal{C}$, $S_2$ is formed by adding a maximal object in $\Ob \mathcal{C} \setminus S_1$ into $S_1$, $S_3$ is formed by adding a maximal object in $\Ob \mathcal{C} \setminus S_2$ into $S_2$, and so on. Consider the full subcategories $\mathcal{D}_i$ formed by $S_i$ and let $e_i = \sum _{x \in S_i} 1_x$. Then the ideals $\mathcal{C} e_i \mathcal{C}$ give a stratification of $\mathcal{C}$ by the previous proposition.
\end{proof}

Now we want to describe standard modules and give a characterization of standardly stratified directed categories with respect to a particular pre-order. Before doing that, we need to define this pre-order on a complete set of primitive idempotents of $\mathcal{C}$ (or precisely, primitive idempotents of the assciated algebra). For every object $x$, $\mathcal{C} (x, x) = 1_x \mathcal{C} 1_x$ is a finite-dimensional $k$-algebra, so we can choose a complete set of orthogonal primitive idempotents $E_x = \{e_{\lambda}\} _{\lambda \in \Lambda_x}$ with $\sum_{\lambda \in \Lambda_x} e_{\lambda} = 1_x$. In this way we get a complete set of orthogonal primitive idempotents $\bigsqcup _{x \in \Ob \mathcal{C}} E_x$. The partial order $\leqslant$ on $\Ob \mathcal{C}$ can be applied to define a pre-ordered set $(\Lambda, \preceq)$ to index all these primitive idempotents, namely for $e_{\lambda} \in E_x$ and $e_{\mu} \in E_y$, $e_{\lambda} \preceq e_{\mu}$ if and only if $x \leqslant y$. We can check that $\preceq$ is indeed a pre-order. We denote $e_{\lambda} \prec e_{\mu}$ if $e_{\lambda} \preceq e_{\mu}$ but $e_{\mu} \npreceq e_{\lambda}$ for $\lambda, \mu \in \Lambda$. Notice that indecomposable summands of $\mathcal{C}$ (viewed as an algebra) can be indexed by these primitive idempotents in a obvious way, namely $P_{\lambda} = \mathcal{C} e_{\lambda}$. Therefore, the pre-ordered set $(\Lambda, \preceq)$ can also be used to index all indecomposable summands of $\mathcal{C}$.

We define standard $\mathcal{C}$-modules in the following way:
\begin{equation*}
\Delta_{\lambda} = P_{\lambda} / \sum_{\mu \npreceq \lambda, \text{ } \mu \in \Lambda} \text{tr} _{P_{\mu}} (P _{\lambda}),
\end{equation*}
where tr$_{P_{\mu}} (P_{\lambda})$ is the trace of $P_{\mu}$ in $P_{\lambda}$. The following proposition gives a description of standard $\mathcal{C}$-modules with respect to the above pre-order.

\begin{proposition}
The standard $\mathcal{C}$-module $\Delta_{\lambda}$ is only supported on $x$ with value $\Delta_{\lambda} (x) \cong 1_x \mathcal{C} e_{\lambda}$, where $x \in \Ob \mathcal{C}$ and $e_{\lambda} \in E_x$.
\end{proposition}

\begin{proof}
Let us first analyze the structure of $P_{\lambda} = \mathcal{C} e_{\lambda}$. Since $e_{\lambda} \in E_x$, $P_{\lambda}$ is a direct summand of $\mathcal{C} 1_x$. The value of $\mathcal{C} 1_x$ on an arbitrary object $y$ is $1_y \mathcal{C} 1_x$, the space of all morphisms from $x$ to $y$. Therefore, the value of $P_{\lambda}$ on $y$ is $1_y \mathcal{C} e_{\lambda}$. By our definition of the partial order on $\Ob \mathcal{C}$, if $x \nleqslant y$, then there is no nontrivial morphisms from $x$ to $y$. Therefore, $1_y \mathcal{C} 1_x$ and hence $1_y \mathcal{C} e_{\lambda}$ are 0. We deduce immediately that $\Delta_{\lambda}$ is only supported on objects $y$ satisfying $x \leqslant y$.

Let $y$ be an object such that $y > x$. Then every $e_{\mu} \in E_y$ satisfies $e_{\mu} \succ e_{\lambda}$. Since $\sum _{e_{\mu} \in E_y} e_{\mu} = 1_y$, by taking the sum we find that tr$ _{\mathcal{C} 1_y} (P_{\lambda})$ is contained in $\sum_{\mu \succ \lambda, \text{ } \mu \in \Lambda} \text{tr} _{P_{\mu}} (P_{\lambda})$. The value on $y$ of tr$ _{\mathcal{C} 1_y} (\mathcal{C} 1_x)$ is $1_y \mathcal{C} 1_x$. Since $P_{\lambda} = \mathcal{C} e_{\lambda}$ is a direct summand of $\mathcal{C} 1_x$, the value on $y$ of tr$ _{\mathcal{C} 1_y} (P_{\lambda})$ is $1_y \mathcal{C} e_{\lambda}$. Consequently, the value of $\sum_{\mu \succ \lambda, \text{ } \mu \in \Lambda} \text{tr} _{P_{\mu}} (P_{\lambda})$ on $y$ contains $1_y \mathcal{C} e_{\lambda}$, which equals the value of $P_{\lambda}$ on $y$. Therefore, the value of $\sum_{\mu \succ \lambda, \text{ } \mu \in \Lambda} \text{tr} _{P_{\mu}} (P_{\lambda})$ on $y$ is precisely $1_y \mathcal{C} e_{\lambda}$, so the value of $\Delta_{\lambda}$ on $y$ is 0.

We have proved that $\Delta_{\lambda}$ is only supported on $x$. Clearly, its value on $x$ is $1_x \mathcal{C} e_{\lambda}$.
\end{proof}

This proposition tells us that standard modules are exactly indecomposable direct summands of $\bigoplus _{x \in \Ob \mathcal{C}} \mathcal{C} (x, x)$ (viewed as a $\mathcal{C}$-module by identifying it with the quotient module $\bigoplus _{x, y \in \Ob \mathcal{C}} \mathcal{C} (x, y) / \bigoplus _{x \neq y \in \Ob \mathcal{C}} \mathcal{C} (x, y)$).

\begin{definition}
A directed category $\mathcal{C}$ is said to be standardly stratified if every (indecomposable) projective module $P_{\lambda}$ has a $\Delta$-filtration by standard modules.
\end{definition}

To simplify the expression, we stick to the following convention frow now on:\\

\textbf{Convention:} When we say a directed category is standardly stratified, we always refer to the preorder $\preceq$ induced by the given partial order $\leqslant$ on the set of objects.\\

This definition is very simple compared to the definition of standardly stratified algebras (for example, the definition in \cite{Cline}). However, from the previous proposition we find that every standard module $\Delta_{\lambda}$ of $\mathcal{C}$ satisfies the following condition: if $S_{\mu}$ and $S_{\nu}$ are two different composition factors of $\Delta_{\lambda}$, then both $S_{\mu} \preceq S_{\nu}$ and $S_{\nu} \preceq S_{\mu}$ (but in general $S_{\mu} \ncong S_{\nu}$). Moreover, If $S_{\mu}$ is a composition factor of the kernel $K_{\lambda}$ of the surjection $P_{\lambda} \rightarrow \Delta_{\lambda}$, then $S_{\mu} \succ S_{\lambda}$ since $K_{\lambda}$ is only supported on objects $y > x$, here $x$ is the object where $S_{\lambda}$ is supported. Therefore, if $\mathcal{C}$ satisfies the requirement in the above definition, then the associated algebra is standardly stratified as well.

The next theorem characterizes standardly stratified directed categories.

\begin{theorem}
Let $\mathcal{C}$ be a directed category. Then $\mathcal{C}$ is standardly stratified if and only if the morphism space $\mathcal{C} (x,y)$ is a projective $\mathcal{C} (y,y)$-module for every pair of objects $x, y \in \Ob \mathcal{C}$. Moreover, if $\mathcal{C}$ is standardly stratified, then $\bigoplus _{x \in \Ob \mathcal{C}}  \mathcal{C} (x, x)$ as a $\mathcal{C}$-module has finite projective dimension.
\end{theorem}

\begin{proof}
Suppose that $\mathcal{C}$ is standardly stratified and take two arbitrary objects $x$ and $y$ in $\mathcal{C}$. Since 0 is regarded as a projective module, we can assume $\mathcal{C} (x, y) \neq 0$ and want to show that it is a projective $\mathcal{C} (y, y)$-module. Consider the projective $\mathcal{C}$-module $\mathcal{C} 1_x$, which has a filtration with factors standard modules. Since each standard module is only supported on one object, the value of $\mathcal{C} 1_x$ on $y$ is exactly the sum of these standard modules with non-zero values on $y$. This sum is direct since standard modules supported on $y$ are non-comparable with respect to the pre-order and therefore there is no extension between them (or because by the previous proposition each of these standard modules is projective viewed as a $\mathcal{C} (y,y)$-module). Therefore, the value of $\mathcal{C} 1_x$ on $y$ is a projective $\mathcal{C} (y,y)$-module. But the value of $\mathcal{C} 1_x$ on $y$ is precisely $\mathcal{C} (x,y)$, so the only if part is proved.

Conversely, let $P_{\lambda} = \mathcal{C} e_{\lambda}$ be an indecomposable projective $\mathcal{C}$-module. Its value on an arbitrary object $y$ is $1_y \mathcal{C} e_{\lambda} \cong 1_y \mathcal{C} 1_x$ which is either 0 or isomorphic to a direct summand of $\mathcal{C} (x, y)$. If $\mathcal{C} (x, y)$ is a projective $\mathcal{C} (y, y)$-module, then the value of $P_{\lambda}$ on $y$ is a projective $\mathcal{C} (y, y)$-module as well. This value can be expressed as a direct sum of standard modules supported on $y$ since standard modules are exactly indecomposable direct summands of $\bigoplus _{x \in \Ob \mathcal{C}} \mathcal{C} (x, x)$. Therefore we can get a filtration of $P_{\lambda}$ by standard modules.

It is well known that the projective dimension of a standard module is finite if the algebra is standardly stratified. Since $\bigoplus _{x \in \Ob \mathcal{C}}  \mathcal{C} (x, x)$ as a $\mathcal{C}$-module is a direct sum of standard modules, the last statement follows from this fact immediately.
\end{proof}

If the directed category $\mathcal{C}$ is standardly stratified, then all standard modules have finite projective dimensions. But the converse is not true in general. However, we will prove later that for a finite EI category, all standard modules have finite projective dimension if and only if this category is standardly stratified with respect to the canonical pre-order.

From now on we suppose that $\mathcal{C}$ is a \textit{graded} category, that is, there is a grading on the morphisms in $\mathcal{C}$ such that $\mathcal{C}_i \cdot \mathcal{C}_j \subseteq \mathcal{C} _{i+j}$, where we denote the subspace spanned by all morphisms with grade $i$ by $\mathcal{C}_i$. Furthermore, $\mathcal{C}$ is supposed to satisfy the following condition: $\mathcal{C}_i \cdot \mathcal{C}_j = \mathcal{C} _{i+j}$. Every vector in $\mathcal{C}_i$ is a linear combination of morphisms with degree $i$. Clearly, $\mathcal{C}_i = \bigoplus _{x, y \in \Ob \mathcal{C}} \mathcal{C} (x, y)_i$. We always suppose $\mathcal{C}_i = 0$ for $i <0$ and $\mathcal{C}_0 = \bigoplus _{x \in \Ob \mathcal{C}} \mathcal{C} (x, x)$. This is equivalent to saying that $\mathcal{C}_0$ is the direct sum of all standard $\mathcal{C}$-modules by Proposition 5.5.

Given a graded directed category $\mathcal{C}$, we can apply the functor $E = \Ext _{\mathcal{C}} ^{\ast} (-, \mathcal{C}_0)$ to construct the \textit{Yoneda category} $E (\mathcal{C}_0)$: $\Ob E(\mathcal{C} _0) = \Ob \mathcal{C}$ and $E (\mathcal{C} _0) (x, y) _n = \Ext ^n _{\mathcal{C}} (\mathcal{C}(x, x), \mathcal{C} (y,y))$. This is precisely the categorical version of Yoneda algebras. By the correspondence between graded algebras and graded categories, we can define \textit{Koszul categories, quasi-Koszul categories, quadratic categories, Koszul modules, quasi-Koszul modules} for graded categories as well. We do not repeat these definitions here but emphasize that all results described in the previous sections can be applied to graded categories.

A corollary of Theorem 5.7 and Corollary 2.18 relates stratification theory to Koszul theory in the context of directed categories.

\begin{theorem}
Let $\mathcal{C}$ be a graded directed category with $\mathcal{C}_0 = \bigoplus _{x \in \Ob \mathcal{C}} \mathcal{C} (x, x)$ self-injective. Then:
\begin{enumerate}
\item $\mathcal{C}$ is standardly stratified if and only if $\mathcal{C}$ is a projective $\mathcal{C}_0$-module.
\item $\mathcal{C}$ is a Koszul category if and only if $\mathcal{C}$ is standardly stratified and quasi-Koszul.
\item If $\mathcal{C}$ is standardly stratified, then a graded $\mathcal{C}$-module $M$ generated in degree 0 is Koszul if and only if it is a quasi-Koszul $\mathcal{C}$-module and a projective $\mathcal{C}_0$-module.
\end{enumerate}
\end{theorem}

\begin{proof}
Take an arbitrary pair of objects $x, y \in \Ob \mathcal{C}$. If $\mathcal{C}$ is standardly stratified, then $\mathcal{C} (x,y)$ is either 0 (a zero projective module) or a projective $\mathcal{C} (y,y)$-module by the previous theorem. Notice that each $\mathcal{C} (x, y)_i$ is a $\mathcal{C} (y,y)$-module since $\mathcal{C} (y,y) \subseteq \mathcal{C}_0$, and we have the decomposition $\mathcal{C} (x,y) = \bigoplus _{i \geqslant 0} \mathcal{C} (x,y)_i$. Therefore, $\mathcal{C}(x,y)_i$ is a projective $\mathcal{C} (y,y)$-module, and hence a projective $\mathcal{C}_0$-module since only the block $\mathcal{C} (y,y)$ of $\mathcal{C}_0$ acts on $\mathcal{C}(x,y)_i$ nontrivially. In conclusion, $\mathcal{C}_i = \bigoplus _{ x, y \in \Ob \mathcal{C}} \mathcal{C} (x,y)_i$ is a projective $\mathcal{C}_0$-module. Conversely, if $\mathcal{C}_i$ are projective $\mathcal{C}_0$-modules for all $i \geqslant 0$, then each $\mathcal{C} (x, y)_i$, and hence $\mathcal{C} (x,y)$ are projective $\mathcal{C} (y,y)$-modules, so $\mathcal{C}$ is standardly stratified again by the previous theorem. The first statement is proved.

We know that $\mathcal{C}$ is Koszul if and only if it is quasi-Koszul and is a projective $\mathcal{C}_0$-module. Then (2) follows from (1) immediately.

The last part is an immediate result of Corollary 2.18.
\end{proof}

To each graded category $\mathcal{C}$ we can associate an \textit{associated quiver} $Q$ in the following way: the vertices of $Q$ are exactly the objects in $\mathcal{C}$; if $\mathcal{C}(x, y)_1 \neq 0$, then we put an arrow from $x$ to $y$ with $x, y$ ranging over all objects in $\mathcal{C}$. Clearly, the associated quiver of $\mathcal{C}$ is completely determined by $\mathcal{C}_0$ and $\mathcal{C}_1$. There is no loop in $Q$ since $\mathcal{C} (x, x)_1 =0$ for each $x \in \Ob \mathcal{C}$.

\begin{proposition}
Let $\mathcal{C}$ be a graded category with $\mathcal{C}_0 = \bigoplus _{x \in \Ob \mathcal{C}} \mathcal{C} (x, x)$ and $Q$ be its associated quiver. Then $\mathcal{C}$ is a directed category if and only if $Q$ is an acyclic quiver.
\end{proposition}

\begin{proof}
Assume that $\mathcal{C}$ is directed. By the definition, there is a partial order $\leqslant$ on $\Ob \mathcal{C}$ such that $\mathcal{C} (x,y) \neq 0$ only if $x \leqslant y$ for $x, y \in \Ob \mathcal{C}$. In particular, $\mathcal{C}(x, y)_1 \neq 0$ only if $x < y$. Therefore, an arrow $x \rightarrow y$ exists in $Q$ only if $x < y$. If there is an oriented cycle
\begin{equation*}
x_1 \rightarrow x_2 \rightarrow \ldots \rightarrow x_n \rightarrow x_1
\end{equation*}
in $Q$, then $x_1 < x_2 < \ldots < x_n < x_1$, which is impossible. Hence $Q$ must be acyclic.

Conversely, if $Q$ is acyclic, we then define $x \leqslant y$ if and only if there is a directed path (including trivial path with the same source and target) from $x$ to $y$ in $Q$ for $x, y \in \Ob \mathcal{C}$. This gives rise to a well defined partial order on $\Ob \mathcal{C}$. We claim that $\mathcal{C}$ is a directed category with respect to this partial order, i.e., $\mathcal{C} (x, y) \neq 0$ implies $x \leqslant y$. Since it holds trivially for $x = y$, we assume that $x \neq y$. Take a morphism $0 \neq \alpha \in \mathcal{C} (x, y)$ with a degree $n$ (this is possible since $\mathcal{C} (x,y)$ is a non-zero graded space). Since $\mathcal{C}_n = \mathcal{C}_1 \cdot \ldots \cdot \mathcal{C}_1$, we can express $\alpha$ as a linear combination of composite morphisms
\begin{equation*}
\xymatrix {x=x_0 \ar[r] ^-{\alpha_1} & x_1 \ar[r] ^-{\alpha_2} & \ldots \ar[r] ^-{\alpha_n} & x_n = y}
\end{equation*}
where each $\alpha_i \in \mathcal{C}_1$ and all $x_i$ are distinct (since endomorphisms in $\mathcal{C}$ are contained in $\mathcal{C}_0$). Therefore, there is a nontrivial directed path
\begin{equation*}
x=x_0 \rightarrow  x_1 \rightarrow x_2 \rightarrow \ldots \rightarrow x_n=y
\end{equation*}
in $Q$, and we have $x < x_1 < x_2 < \ldots < y$, which proves our claim.
\end{proof}

Let $\mathcal{C}$ be a graded directed category. We define the \textit{free cover} $\hat {\mathcal{C}}$ of $\mathcal{C}$ by using the associated quiver $Q$. Explicitly, $\hat {\mathcal{C}}$ has the same objects and endomorphisms as $\mathcal{C}$. For each pair of objects $x \neq y$ we construct $\hat{ \mathcal {C}} (x, y)$ as follows. let $\Gamma_{x,y}$ be the set of all paths from $x$ to $y$ in $Q$. In the case that $\Gamma_{x,y} = \emptyset$ we let $\hat{ \mathcal {C}} (x, y) = 0$. Otherwise, take an arbitrary path $\gamma \in \Gamma_{x,y}$ pictured as below
\begin{equation*}
x \rightarrow x_1 \rightarrow x_2 \rightarrow \ldots \rightarrow x_{n-1} \rightarrow y,
\end{equation*}
and define $(x, y)_{\gamma}$ to be
\begin{equation*}
\mathcal{C}(x_{n-1}, y)_1 \otimes _{\mathcal{C} (x_{n-1}, x_{n-1})} \mathcal{C}(x_{n-2}, x_{n-1})_1 \otimes _{\mathcal{C} (x_{n-2}, x_{n-2})} \ldots \otimes _{\mathcal{C} (x_1, x_1)} \mathcal{C}(x, x_1)_1.
\end{equation*}
Finally, we define
\begin{equation*}
\hat {\mathcal{C}} (x, y) = \bigoplus _{\gamma \in \Gamma_{x, y}} (x, y)_{\gamma}.
\end{equation*}
It is clear that $\hat{\mathcal{C}}$ is also a graded category with $\hat{\mathcal{C}}_0 = \mathcal{C}_0$ and $\hat{\mathcal{C}}_1 = \mathcal{C}_1$. Therefore, $\hat{\mathcal{C}}$ has the same associated quiver as that of $\mathcal{C}$ and is also a directed category by the above lemma. Actually, if two grade categories $\mathcal{C}$ and $\mathcal{D}$ have the same degree 0 and degree 1 components, then one is a directed category if and only if so is the other.

\begin{theorem}
Let $\mathcal{C}$ be a directed Koszul category with $\mathcal{C}_0$ self-injective, then the Yoneda category $\mathcal{E} = E(\mathcal{C}_0)$ is also directed and Koszul.
\end{theorem}

\begin{proof}
Applying the Koszul duality (Theorem 4.1) we know that $\mathcal{E}$ is a Koszul category. What we need to show is that $\mathcal{E}$ is a directed category as well. Since $\mathcal{C}$ is standardly stratified, $\pd _{\mathcal{C}} \mathcal{C} _0 < \infty$. Therefore, all morphisms in $\mathcal{E}$ spans a finite-dimensional space by Lemma 4.3. In particular, for $x, y \in \Ob \mathcal{E}$, $ \dim_k \mathcal{E} (x, y) < \infty$. Therefore, $\mathcal{E}$ is a locally finite $k$-linear category.

Let $\leqslant$ be the partial order on $\Ob \mathcal{C}$ with respect to which $\mathcal{C}$ is directed. This partial order gives a partial order on $\Ob \mathcal{E}$ as well because $\Ob \mathcal{E} = \Ob \mathcal{C}$. We claim that $\mathcal{E}$ is directed with respect to this partially ordered set, i.e., if $x \nleq y$ are two distinct objects in $\mathcal{E}$, then $\mathcal{E} (x, y) = 0$.

Since $\mathcal{E}$ is the Yoneda category of $\mathcal{C}$, $\mathcal{E} (x, y) = 1_y \mathcal{E} 1_x \cong \Ext _{\mathcal{C}} ^{\ast} (\mathcal{C}_0 1_x, \mathcal{C}_0 1_y)$. But $\mathcal{C}$ is a Koszul category, so $\mathcal{C}_i$ is a projective $\mathcal{C}_0$-modules for each $i \geqslant 0$ by Theorem 5.8. Therefore, all $\Omega^i (\mathcal{C}_0 1_x)_i$ are projective $\mathcal{C}_0$-modules (Lemma 2.14) and we have (Lemma 2.10)
\begin{equation*}
\mathcal{E} (x, y)_i \cong \Ext _{\mathcal{C}} ^i (\mathcal{C}_0 1_x, \mathcal{C}_0 1_y) \cong \Hom _{\mathcal{C}} (\Omega^i (\mathcal{C}_0 1_x), \mathcal{C}_0 1_y).
\end{equation*}
Observe that $\mathcal{C}_0 1_y$ is only supported on the object $y$ and $y \ngeq x$. If we can prove the statement that each $\Omega^i (\mathcal{C}_0 1_x)$ is only supported on objects $z$ with $z \geqslant x$, then our claim is proved.

Clearly, $\Omega^0 (\mathcal{C}_0 1_x) = \mathcal{C}_0 1_x = \mathcal{C} (x,x)$ is only supported on $x$, so the statement is true for $i =0$. Now suppose that $\Omega^n (\mathcal{C}_0 1_x)$ is only supported on objects $z \geqslant x$ and consider $\Omega^{n+1} (\mathcal{C}_0 1_x)$. Let $S$ be the set of objects $z$ such that the value $1_z \Omega^n (\mathcal{C}_0 1_x)$ of $\Omega^n (\mathcal{C}_0 1_x)$ on $z$ is non-zero. Then we can find a short exact sequence:
\begin{equation*}
\xymatrix {0 \ar[r] & N \ar[r] & \bigoplus _{z \in S} (\mathcal{C} 1_z) ^{m_z} \ar[r]^p & \Omega^n (\mathcal{C}_0 1_x) \ar[r] & 0}
\end{equation*}
such that the map $p$ gives a surjection $p_z: (1_z \mathcal{C} 1_z) ^{m_z} \rightarrow 1_z \Omega^n (\mathcal{C}_0 1_x)$ for $z \in S$. Thus $p$ is a surjection and $\Omega^{n+1} (\mathcal{C}_0 1_x)$ is a direct summand of $N$. Notice that all $\mathcal{C} 1_z$ are supported only on objects $w \geqslant z$, and $z \geqslant x$ by the induction hypothesis. Therefore, the submodule $\Omega^{n+1} (\mathcal{C}_0 1_x) \subseteq N \subseteq \bigoplus _{z \in S} (\mathcal{C} 1_z) ^{m_z} $ is only supported on objects $w \geqslant x$. Our statement is proved by induction. This finishes the proof.
\end{proof}

The condition that $\mathcal{C}_0$ is exactly the space spanned by endomorphisms in $\mathcal{C}$ is crucial since it implies that standard modules are precisely indecomposable summands of $\mathcal{C}_0$. The following example tells us that without this assumption, the Yoneda category $E (\mathcal{C} _0)$ might not be directed even if $\mathcal{C}$ is a Koszul directed category.

\begin{example}
Let $\mathcal{C}$ be the following category. Put an order $x < y$ on the objects and the following grading on morphisms: $\mathcal{C}_0 = \langle 1_x, 1_y, \beta \rangle$, $\mathcal{C}_1 = \langle \alpha \rangle$.
\begin{equation*}
\xymatrix{ x \ar@/^/ [rr]^{\alpha} \ar@ /_/ [rr]_{\beta} & & y.}
\end{equation*}
This category is directed obviously. It is standardly stratified (actually hereditary) with $\Delta_x \cong k_x$ and $\Delta_y \cong k_y$. By the exact sequence
\begin{equation*}
\xymatrix{ 0 \ar[r] & k_y [1] \ar[r] & \mathcal{C} \ar[r] & \mathcal{C}_0 \ar[r] & 0,}
\end{equation*}
$\mathcal{C} _0$ is a Koszul module. But $\Delta_x \oplus \Delta_y \ncong \mathcal{C} _0$. Furthermore, $\Delta_x \oplus \Delta_y$ is not Koszul since from the short exact sequence
\begin{equation*}
\xymatrix{ 0 \ar[r] & k_y [1] \oplus k_y \ar[r] & \mathcal{C} \ar[r] & \Delta_x \oplus \Delta_y \ar[r] & 0}
\end{equation*}
we find that $\Omega(\Delta_x \oplus \Delta_y)$ is not generated in degree 1.

By computation we get the Yoneda category $\mathcal{D} = E (\mathcal{C} _0)$ pictured as below, with relation $\alpha \cdot \beta = 0$.
\begin{equation*}
\xymatrix{ x \ar@/^/ [rr]^{\alpha} & & y \ar@ /^/ [ll]^{\beta}.}
\end{equation*}
This is not a directed category with respect to the order $x <y$. However, we check: $P_y = \Delta'_y = \mathcal{D}_0 1_y$ and $\Delta'_x = k_x \cong P_x/P_y$. Therefore, $\Delta'_x \oplus \Delta'_y \cong \mathcal{D}_0$, and $\mathcal{D}$ is standardly stratified. The exact sequence
\begin{equation*}
\xymatrix{ 0 \ar[r] & P_y [1] \ar[r] & \mathcal{D} \ar[r] & \mathcal{D}_0 \ar[r] & 0}
\end{equation*}
tells us that $\mathcal{D}_0$ is a Koszul $\mathcal{D}$-module. Therefore, $\mathcal{D}$ is standardly stratified and Koszul, but not directed.
\end{example}

There is a close relation between the classical Koszul theory and our generalized Koszul theory in the context of graded directed categories. Let $\mathcal{C}$ be a graded directed category. We then define a subcategory $\mathcal{D}$ of $\mathcal{C}$ by replacing all endomorphism rings in $\mathcal{C}$ by $k \cdot 1$, the span of the identity endomorphism. Explicitly, $\Ob \mathcal{D} = \Ob \mathcal{C}$; for $x, y \in \Ob \mathcal{D}$, $\mathcal{D} (x, y) = k \langle 1_x \rangle $ if $x = y$ and $ \mathcal{D} (x, y) = \mathcal{C} (x,y)$ otherwise. Clearly, $\mathcal{D}$ is also a graded directed category with $\mathcal{D}_i = \mathcal{C}_i$ for every $i \geqslant 1$. Observe that the degree 0 component $\mathcal{D}_0$ is semisimple.

\begin{theorem}
Let $\mathcal{C}$ be a graded directed Koszul category and define the subcategory $\mathcal{D}$ as above. If $M$ is a Koszul $\mathcal{C}$-module, then $M \downarrow _{\mathcal{D}} ^{\mathcal{C}}$ is a Koszul $\mathcal{D}$-module. In particular, $\mathcal{D}$ is a Koszul category in the classical sense.
\end{theorem}

\begin{proof}
We prove the conclusion by induction on the size of $\Ob \mathcal{C}$. If the size of $\Ob \mathcal{C}$ is 1, the conclusion holds trivially. Now suppose that the conclusion is true for categories with at most $n$ objects and let $\mathcal{C}$ be a graded directed category with $n+1$ objects. Take $x$ to be a minimal object in $\mathcal{C}$ and define $\mathcal{C}_x$ ($\mathcal{D}_x$, resp.) to be the full subcategory of $\mathcal{C}$ ($\mathcal{D}$, resp.) formed by removing $x$ from it. Clearly $\mathcal{C}_x$ and $\mathcal{D}_x$ have $n$ objects.

The following fact, which is well known in the context of finite EI categories (see \cite{Xu}), is essential in the proof.\\

\textbf{Fact:} Every graded $\mathcal{C}_x$-module $N$ can be viewed as a $\mathcal{C}$-module with $N(x) =0$ by induction. Conversely, every graded $\mathcal{C}$-module $M$ with $M(x) = 0$ can also be viewed as a $\mathcal{C}_x$-module by restriction. Furthermore, if $M(x) = 0$, then $\Omega^i (M)(x) = 0$ for all $i \geqslant 0$. The above induction and restriction preserves projective modules: a projective $\mathcal{C}_x$-module is still projective when viewed as a $\mathcal{C}$-module; conversely, a projective $\mathcal{C}$-module $P$ with $P(x) =0$ is still projective viewed as a $\mathcal{C}_x$-module. Therefore, a graded $\mathcal{C}$-module $M$ with $M(x)=0$ is Koszul if and only if it is Koszul as a $\mathcal{C}_x$-module. All these results hold for the pair $(\mathcal{D}, \mathcal{D}_x)$ similarly.\\

By this fact, we only need to handle Koszul $\mathcal{C}$-modules $M$ with $M(x) \neq 0$. Indeed, if $M(x) = 0$, then $M$ is also Koszul regarded as a $\mathcal{C}_x$-module. By the induction hypothesis, $M \downarrow _{\mathcal{D}_x} ^{\mathcal{C}_x}$ is a Koszul $\mathcal{D}_x$-module. By the above fact, $M \downarrow _{\mathcal{D}} ^{\mathcal{C}}$ is a Koszul $\mathcal{D}$-module. Thus the conclusion is true for Koszul $\mathcal{C}$-modules $M$ with $M(x) =0$.

Firstly we consider the special case $M = \mathcal{C}_0 1_x = \mathcal{C} (x,x)$ which is concentrated on $x$ when viewed as a $\mathcal{C}$-module. It is clear that
\begin{equation*}
\Omega (\mathcal{C}_0 1_x) \downarrow _{\mathcal{D}} ^{\mathcal{C}} = \Omega (\mathcal{D}_0 1_x)
\end{equation*}
as vector spaces since for each pair $u \neq v \in \Ob \mathcal{C}$, $\mathcal{C} (u, v) = \mathcal{D} (u, v)$, and
\begin{equation*}
\mathcal{C}_0 1_x \downarrow _{\mathcal{D}} ^{\mathcal{C}} \cong (\mathcal{D}_0 1_x)^m \cong k_x^m,
\end{equation*}
where $m = \dim_k \mathcal{C} (x, x)$. Since $\mathcal{C}_01_x$ is a Koszul $\mathcal{C}$-module, $\Omega (\mathcal{C}_0 1_x) [-1]$ is a Koszul $\mathcal{C}$-module supported on $\Ob \mathcal{C}_x$. By the induction hypothesis, $\Omega(\mathcal{D}_0 1_x) [-1] = \Omega (\mathcal{C}_0 1_x) \downarrow _{\mathcal{D}} ^{\mathcal{C}} [-1]$ is a Koszul $\mathcal{D}_x$-module, and hence a Koszul $\mathcal{D}$-module. Therefore, $\mathcal{D}_0 1_x$, and hence $\mathcal{C}_0 1_x \downarrow _{\mathcal{D}} ^{\mathcal{C}} \cong (\mathcal{D}_0 1_x)^m$ are Koszul $\mathcal{D}$-modules. In the case that $y \neq x$, $\mathcal{D}_0 1_y$ is a direct summand of $\mathcal{C}_0 1_y \downarrow _{\mathcal{D}} ^{\mathcal{C}}$. It is Koszul viewed as a $\mathcal{D}_x$-module by the induction hypothesis, and hence is a Koszul $\mathcal{D}$-module. Consequently, $\mathcal{D}_0$ is a Koszul $\mathcal{D}$-module, so $\mathcal{D}$ is a Koszul category in the classical sense.

Now let $M$ be an arbitrary Koszul $\mathcal{C}$-module with $M(x) \neq 0$. Consider the exact sequence
\begin{equation}
\xymatrix {0 \ar[r] & \Omega M \downarrow _{\mathcal{D}} ^{\mathcal{C}} \ar[r] & P \downarrow _{\mathcal{D}} ^{\mathcal{C}} \ar[r]^p & M \downarrow _{\mathcal{D}} ^{\mathcal{C}} \ar[r] & 0}
\end{equation}
induced by
\begin{equation*}
\xymatrix {0 \ar[r] & \Omega M \ar[r] & P \ar[r]^p & M \ar[r] & 0.}
\end{equation*}
The structures of $\mathcal{C}$ and $\mathcal{D}$ give the following exact sequence:
\begin{equation*}
\xymatrix {0 \ar[r] & _{\mathcal{D}} \mathcal{D} \ar[r] & _{\mathcal{D}} \mathcal{C} \ar[r] & \bigoplus _{x \in \Ob \mathcal{C}} k_x^{m_x} \ar[r] & 0}
\end{equation*}
where $k_x \cong \mathcal{D}_01_x$ and $m_x = \dim_k \mathcal{C} (x, x) - 1$. Since $P \downarrow _{\mathcal{D}} ^{\mathcal{C}} \in \text{add} (_{\mathcal{D}} \mathcal{C})$, the above sequence gives us a corresponding sequence for $P \downarrow _{\mathcal{D}} ^{\mathcal{C}}$:
\begin{equation}
\xymatrix {0 \ar[r] & P' \ar[r]^{\iota} & P \downarrow _{\mathcal{D}} ^{\mathcal{C}} \ar[r] & T \ar[r] & 0,}
\end{equation}
here $P'$ is a projective $\mathcal{D}$-module and $T \in \text{add} (\mathcal{D}_0)$. Both of them are generated in degree 0.

Putting sequences 5.1 and 5.2 together we get
\begin{equation*}
\xymatrix{ 0 \ar[r] & \Omega (M') \ar[r] \ar[d]^{\varphi} & P' \ar[r]^{p \circ \iota} \ar[d]^{\iota} & M' \ar[d] \ar[r] & 0 \\
0 \ar[r] & \Omega (M) \downarrow _{\mathcal{D}} ^{\mathcal{C}} \ar[r] & P \downarrow _{\mathcal{D}} ^{\mathcal{C}} \ar[r]^p \ar[d] & M \downarrow _{\mathcal{D}} ^{\mathcal{C}} \ar[r] \ar[d] & 0 \\
& & T \ar[r] &  M \downarrow _{\mathcal{D}} ^{\mathcal{C}} /M'
}
\end{equation*}
where $M' = (p \circ \iota) (P')$. Notice that $p$ and $\iota$ both are injective restricted to degree 0 components. Therefore $p \circ \iota$ is also injective restricted to the degree 0 component of $P'$ (actually it is an isomorphism restricted to the degree 0 component). Consequently, $P'$ is a projective cover of $M'$ and the kernel of $p \circ \iota$ is indeed $\Omega (M')$.

We claim that $\varphi$ is an isomorphism and hence $\Omega(M') \cong (\Omega M) \downarrow _{\mathcal{D}} ^{\mathcal{C}}$. It suffices to show $T \cong M \downarrow _{\mathcal{D}} ^{\mathcal{C}} /M'$ by the snake lemma. First, since $T$ is concentrated in degree 0 in sequence 5.2, $(P \downarrow _{\mathcal{D}} ^{\mathcal{C}} )_i = P'_i$ and
\begin{equation*}
M'_i = (p \circ \iota) (P'_i) = p ((P \downarrow _{\mathcal{D}} ^{\mathcal{C}} )_i) = (M \downarrow _{\mathcal{D}} ^{\mathcal{C}}) _i, i >0.
\end{equation*}
Therefore, $ M \downarrow _{\mathcal{D}} ^{\mathcal{C}} /M'$ is concentrated in degree 0. Furthermore, since $M$ is Koszul, $(P \downarrow _\mathcal{D} ^{\mathcal{C}} )_0 = P_0 = M_0 = (M \downarrow _{\mathcal{D}} ^{\mathcal{C}})_0$, and $P'_0 = M'_0$ because $p\circ \iota$ restricted to $P_0'$ is an isomorphism as well. We deduce that
\begin{equation*}
T \cong (P \downarrow _\mathcal{D} ^{\mathcal{C}} )_0 / P'_0 \cong (M \downarrow _\mathcal{D} ^{\mathcal{C}} )_0 / M'_0 = M \downarrow _{\mathcal{D}} ^{\mathcal{C}} /M',
\end{equation*}
exactly as we claimed.

Now consider the rightmost column of the above diagram. Clearly, the bottom term $M \downarrow _{\mathcal{D}} ^{\mathcal{C}} /M' \cong T \in \text{add} (\mathcal{D}_0)$ is Koszul since we just proved that $\mathcal{D}_0$ is Koszul. The $\mathcal{C}$-module $(\Omega M) [-1]$ is Koszul since $M$ is supposed to be Koszul. Moreover, because $x$ is minimal and $M$ is generated in degree 0, $M(x) \subseteq M_0$ and hence $(\Omega M) [-1] (x) = (\Omega M)(x) = 0$. Therefore, $(\Omega M) [-1]$ is a Koszul $\mathcal{C}$-module supported on $\Ob \mathcal{C}_x$, so it is also a Koszul $\mathcal{C}_x$-module. By the induction hypothesis, $(\Omega M') [-1] \cong (\Omega M) [-1] \downarrow _{\mathcal{D}} ^{\mathcal{C}}$ is Koszul viewed as a $\mathcal{D}_x$-module, and hence Koszul as a $\mathcal{D}$-module. Thus the top term $M'$ is a Koszul $\mathcal{D}$-module since as a homomorphic image of $P'$ (which is generated in degree 0) it is generated in degree 0 as well. By Proposition 2.9, $M \downarrow _{\mathcal{D}} ^{\mathcal{C}}$ is also a Koszul $\mathcal{D}$-module since $\mathcal{D}_0$ is semisimple by our construction. The conclusion follows from induction.
\end{proof}

The converse of the above theorem is also true.

\begin{theorem}
Let $\mathcal{C}$ be a graded directed category and construct the subcategory $\mathcal{D}$ as before. Suppose that $\mathcal{D}$ is Koszul in the classical sense. Let $M$ be a graded $\mathcal{C}$-module generated in degree 0 such that $\Omega^i(M)_i$ are projective $\mathcal{C}_0$-modules for all $i \geqslant 0$. Then $M$ is a Koszul $\mathcal{C}$-module whenever $M \downarrow _{\mathcal{D}} ^{\mathcal{C}}$ is a Koszul $\mathcal{D}$-module.
\end{theorem}

\begin{proof}
We use the similar technique to prove the conclusion. Notice that we always assume that $\mathcal{C}_0 = \bigoplus _{x \in \Ob \mathcal{C}} \mathcal{C} (x,x)$. If $\mathcal{C}$ has only one object, then Koszul modules are exactly projective modules generated in degree 0 and the conclusion holds. Suppose that it is true for categories with at most $n$ objects. Let $\mathcal{C}$ be a category of $n+1$ objects and take a minimal object $x$. Define $\mathcal{C}_x$ and $\mathcal{D}_x$ as before. As in the proof of last theorem, a graded $\mathcal{C}$-module $M$ with $M(x) =0$ is Koszul if and only if it is Koszul viewed as a $\mathcal{C}_x$-module by restriction. and the same result holds for the pair $(\mathcal{D}, \mathcal{D}_x)$. In particular, $\mathcal{D}_x$ is a Koszul category. Therefore, we only need to show that an arbitrary graded $\mathcal{C}$-module $M$ which is generated in degree 0 and satisfies the following conditions is Koszul: $\Omega^i(M)_i$ is a projective $\mathcal{C}_0$-module for each $i \geqslant 0$; $M(x) \neq 0$; and $M \downarrow _{\mathcal{D}} ^{\mathcal{C}}$ is a Koszul $\mathcal{D}$-module.

Let $M$ be such a $\mathcal{C}$-module and consider the commutative diagram:

\begin{equation}
\xymatrix{ & K \ar@{=}[r] \ar[d] & K \ar[d] \ar[r] & 0 \ar[d] \\
0 \ar[r] & \Omega(M \downarrow _{\mathcal{D}} ^{\mathcal{C}}) \ar[d] \ar[r] & P \ar[r] \ar[d]^{\varphi} & M \downarrow _{\mathcal{D}} ^{\mathcal{C}} \ar[r] \ar@{=}[d] & 0 \\
0 \ar[r] & (\Omega M) \downarrow _{\mathcal{D}} ^{\mathcal{C}} \ar[r] & \tilde{P} \downarrow _{\mathcal{D}} ^{\mathcal{C}} \ar[r] & M \downarrow _{\mathcal{D}} ^{\mathcal{C}} \ar[r] & 0}
\end{equation}
where $P$ and $\tilde{P}$ are projective covers of $M \downarrow _{\mathcal{D}} ^{\mathcal{C}}$ and $M$ respectively. Since $M_0$ is a projective $\mathcal{C}_0$-module, $P_0 = M_0 = (M \downarrow _{\mathcal{D}} ^{\mathcal{C}} )_0 =(\tilde{P} \downarrow _{\mathcal{D}} ^{\mathcal{C}} )_0$ as vector spaces, and the induced map $\varphi$ restricted to $P_0$ is an isomorphism. Therefore $\varphi$ is surjective since both $P$ and $\tilde{P} \downarrow _{\mathcal{D}} ^{\mathcal{C}}$ are generated in degree 0. Let $K$ be the kernel of $\varphi$.

We have the following exact sequence similar to sequence 5.2:
\begin{equation*}
\xymatrix{ 0 \ar[r] & P' \ar[r] & \tilde{P} \downarrow _{\mathcal{D}} ^{\mathcal{C}} \ar[r] ^{\tilde{p}} & T \ar[r] & 0, }
\end{equation*}
where $P'$ is a projective $\mathcal{D}$-module such that $P'_i = (\tilde{P} \downarrow _{\mathcal{D}} ^{\mathcal{C}} )_i$ for every $i \geqslant 1$, and $T \in \text{add} (\mathcal{D}_0)$.

Let $P''$ be a projective cover of $T$ (as a $\mathcal{D}$-module). Then we obtain:

\begin{equation}
\xymatrix{ & & K \ar@{=}[r] \ar[d] & K \ar[d] \\
0 \ar[r] & P' \ar@{=}[d] ^{\alpha} \ar[r] & P \ar[r]^p \ar[d]^{\varphi} & P'' \ar[d]^{p''} \ar[r] & 0 \\
0 \ar[r] & P' \ar[r] & \tilde{P} \downarrow _{\mathcal{D}} ^{\mathcal{C}} \ar[r] ^{\tilde{p}} & T \ar[r] & 0}.
\end{equation}
We give some explanations here. Since $P$ is a projective $\mathcal{D}$-module and the map $p''$ is surjective, the map $\tilde{p} \circ \varphi$ factors through $p''$ and gives a map $p: P \rightarrow P''$. Restricted to degree 0 components, $p''$ and $\varphi$ (see diagram 5.3) are isomorphisms and $\tilde{p}$ is surjective. Thus $p$ restricted to the degree 0 components is also surjective. But $P''$ is generated in degree 0, so $p$ is surjective. Since $P_0 = (\tilde{P} \downarrow _{\mathcal{D}} ^{\mathcal{C}} )_0$ and $P''_0 = T_0 = T$, $\alpha$ restricted to the degree 0 components is an isomorphism, and hence an isomorphism of projective $k\mathcal{D}$-modules (notice that the middle row splits since $P''$ is a projective $\mathcal{D}$-module, so the kernel should be a projective $\mathcal{D}$-module generated in degree 0). By the snake Lemma, the kernel of $p''$ is also $K$ up to isomorphism.

Let $J = \bigoplus _{i \geqslant 1} \mathcal{D}_i$. Since $\mathcal{D}_0$ is supposed to be a Koszul $\mathcal{D}$-module, $J [-1] \cong \Omega (\mathcal{D} _0) [-1]$ is a Koszul $\mathcal{D}$-module, too. Consider the leftmost column in diagram 5.3. The top term $K[-1]$ is a Koszul $\mathcal{D}$-module since $K \cong P'' / T \cong P''/P''_0 \in \text{add} (J)$. The middle term $\Omega (M \downarrow _{\mathcal{D}} ^{\mathcal{C}}) [-1]$ is Koszul as well since $M \downarrow _{\mathcal{D}} ^{\mathcal{C}}$ is supposed to be Koszul. By Proposition 2.9, the bottom term $(\Omega M)[-1] \downarrow _{\mathcal{D}} ^{\mathcal{C}}$ must be Koszul.

Since $M$ is generated in degree 0 and $x$ is a minimal object, $M(x) \subseteq M_0$, so $M(x) = (M \downarrow _{\mathcal{D}} ^{\mathcal{C}}) (x) \subseteq (M \downarrow _{\mathcal{D}} ^{\mathcal{C}}) _0$ as vector spaces. Similarly, $P(x) \subseteq P_0$ and $P(x) \cong (M \downarrow _{\mathcal{D}} ^{\mathcal{C}}) (x)$, so $\Omega (M \downarrow _{\mathcal{D}} ^{\mathcal{C}}) (x) = 0$. Consequently, $(\Omega M)[-1]$ is supported on $\Ob \mathcal{C}_x$ by observing the leftmost column of diagram 5.3. Moreover, we can show as in the proof of Theorem 5.10 that all of its syzygies are supported on $\Ob \mathcal{C}_x$, and $\Omega^i ((\Omega M)[-1])_i = \Omega^{i+1} (M)_{i+1}$ are projective $(\mathcal{C}_x)_0$-modules. Therefore, applying the induction hypothesis to $(\Omega M)[-1]$ supported on $\Ob \mathcal{C}_x$ and $(\Omega M) [-1] \downarrow _{\mathcal{D}} ^{\mathcal{C}}$ supported on $\Ob \mathcal{D}_x$, we conclude that $(\Omega M) [-1]$ is a Koszul $\mathcal{C} _x$-module, and hence a Koszul $\mathcal{C}$-module. Clearly, $M$ is a Koszul $\mathcal{C}$-module since it is generated in degree 0. The conclusion follows from induction.
\end{proof}

\begin{remark}
We remind the reader that in the previous two theorems we do not require $\mathcal{C}_0 = \bigoplus _{x \in \Ob \mathcal{C}} \mathcal{C} (x,x)$ to be a self-injective algebra. By our construction, $\mathcal{D}_0 \cong \bigoplus _{x \in \Ob \mathcal{D}} k_x$ is a semisimple algebra.
\end{remark}

Assuming that $\mathcal{C}_0$ is self-injective, we get the following nice correspondence.

\begin{theorem}
Let $\mathcal{C}$ be a graded directed category with $\mathcal{C}_0$ self-injective and construct $\mathcal{D}$ as before. Then:
\begin{enumerate}
\item $\mathcal{C}$ is a Koszul category in our sense if and only if $\mathcal{C}$ is standardly stratified and $\mathcal{D}$ is a Koszul category in the classical sense.
\item If $\mathcal{C}$ is a Koszul category, then a graded $\mathcal{C}$-module $M$ generated in degree 0 is Koszul if and only if $M \downarrow _{\mathcal{D}} ^{\mathcal{C}}$ is a Koszul $\mathcal{D}$-module and $M$ is a projective $\mathcal{C}_0$-module.
\end{enumerate}
\end{theorem}

\begin{proof}
If $\mathcal{C}$ is Koszul in our sense, then it is standardly stratified by (2) of Theorem 5.8, and $\mathcal{D}$ is Koszul in the classical sense by Theorem 5.12. Conversely, if $\mathcal{D}$ is Koszul in the classical sense, then $\mathcal{C}_0 \downarrow _{\mathcal{D}} ^{\mathcal{C}} \in \text{add} (\mathcal{D}_0)$ is a Koszul $\mathcal{D}$-module. If $\mathcal{C}$ is furthermore standardly stratified, then it is a projective $\mathcal{C}_0$-module by Theorem 5.8. Therefore, all $\Omega^i( \mathcal{C}_0 )_i$ are projective $\mathcal{C}_0$-modules according to Lemma 2.14. Thus $\mathcal{C}_0$ is a Koszul $\mathcal{C}$-module by Theorem 5.13, and hence $\mathcal{C}$ is a Koszul category. This proves the first statement.

Now suppose that $\mathcal{C}$ is Koszul. Then $\mathcal{C}$ is a projective $\mathcal{C}_0$-module. If $M$ is a Koszul $\mathcal{C}$-module, $M \downarrow _{\mathcal{D}} ^{\mathcal{C}}$ is a Koszul $\mathcal{D}$-module by Theorem 5.12. Furthermore, $M$ is a projective $\mathcal{C}_0$-module by Corollary 2.18. Conversely, if $M$ is a projective $\mathcal{C}_0$-module and $M \downarrow _{\mathcal{D}} ^{\mathcal{C}}$ is a Koszul $\mathcal{D}$-module, then by Lemma 2.14 $\Omega^i (M)_i$ are projective $\mathcal{C}_0$-modules for all $i \geqslant 0$. By Theorem 5.13 $M$ is a Koszul $\mathcal{C}$-module.
\end{proof}

\section{Finite EI Categories}

When applying the generalized Koszul theory to a directed category $\mathcal{C}$ in the previous section, we take for granted that there is already a grading on the morphisms in $\mathcal{C}$ such that the degree 0 component is formed precisely by endomorphisms in $\mathcal{C}$. But in practice it is very hard to find such a grading for $\mathcal{C}$. Actually, we do not even know the existence of such gradings in general. In this section we will focus on finite EI categories, whose $k$-linearizations form a type of directed categories with combinatorial properties. These properties can be used to define a length grading on the set of morphisms and completely determine whether an arbitrary finite EI category can be graded by this length grading.

In this section we only consider skeletal and \textit{connected} finite EI categories $\mathcal{E}$, i.e., for every pair $x, y\in \Ob \mathcal{E}$, there is a chain of objects $x_0 =x, x_1, x_2, \ldots, x_n = y$ such that either $\mathcal{E} (x_i, x_{i+1}) \neq \emptyset$ or $\mathcal{E} (x_{i+1}, x_i) \neq \emptyset$ for every $0 \leqslant i \leqslant n-1$.

First we introduce some results from \cite{Li}. A morphism $\alpha$ in $\mathcal{E}$ is called \textit{unfactorizable} if $\alpha$ is not an automorphism, and whenever there is a decomposition $\alpha = \alpha_1 \circ \alpha_2$, either $\alpha_1$ or $\alpha_2$ is an automorphism. The composite morphism of an unfactorizable morphism and an automorphism is still unfactorizable. Therefore, all unfactorizable morphisms from an object $x$ to another object $y$ form an $(\Aut _{\mathcal{E}} (y), \Aut _{\mathcal{E}} (x))$ bi-set. Every non-isomorphism can be expressed as a composite of unfactorizable morphisms. This decomposition is not unique in general. We say a finite EI category $\mathcal{E}$ satisfies the \textit{Unique Factorization Property} (UFP) if whenever a non-isomorphism $\alpha$ has two decompositions into unfactorizable morphisms:
\begin{equation*}
\xymatrix{x=x_0 \ar[r]^{\alpha_1} & x_1 \ar[r]^{\alpha_2} & \ldots \ar[r]^{\alpha_m} &x_m=y \\
x=x_0 \ar[r]^{\beta_1} & y_1 \ar[r]^{\beta_2} & \ldots \ar[r]^{\beta_n} & y_n=y}
\end{equation*}
then $m=n$, $x_i = y_i$, and there are $ h_i \in \Aut _{\mathcal{E}} (x_i)$ such that the following diagram commutes, $1 \leqslant i \leqslant n-1$:
\begin{align*}
\xymatrix{ x_0 \ar[r]^{\alpha_1} \ar@{=}[d]^{id} & x_1 \ar[r]^{\alpha_2} \ar[d]^{h_1} & \ldots \ar[r]^{\alpha_{\ldots}} \ar[d]^{h_{\ldots}} & x_{n-1} \ar[r]^{\alpha_n} \ar[d]^{h_{n-1}} & x_n \ar@{=}[d]^{id} \\
x_0 \ar[r]^{\beta_1} & x_1 \ar[r]^{\beta_2} & \ldots \ar[r]^{\beta_{\ldots}} & x_{n-1} \ar[r]^{\beta_n} & x_n}
\end{align*}
Finite EI categories with this property are called \textit{finite free EI categories} by us. For every finite EI category $\mathcal{E}$ there is a unique (up to isomorphism) finite free EI category $\hat{\mathcal{E}}$ (called the \textit{free EI cover}) and a covering functor $\hat{F}: \hat{\mathcal{E}} \rightarrow \mathcal{E}$ such that $\hat{F}$ is the identity map restricted to objects, isomorphisms and unfactorizable morphisms. The functor $\hat{F}$ induces a surjective algebra homomorphism $\varphi: k \hat{\mathcal{E}} \rightarrow k\mathcal{E}$. Therefore, $k\mathcal{E} \cong k \hat{\mathcal{E}} /I$, where $I$ is the kernel of $\varphi$. We have the following description of $I$:

\begin{lemma}
The $k \hat{\mathcal{E}}$-ideal $I$ as a vector space is spanned by elements of the form $\hat{\alpha} - \hat{\beta}$, where $\hat{\alpha}$ and $\hat{\beta}$ are morphisms in $\hat{\mathcal{E}}$ with $\hat{F} (\hat{\alpha}) = \hat{F} (\hat{\beta})$.
\end{lemma}

\begin{proof}
Let $U$ be the vector space spanned by elements $\hat{\alpha} - \hat{\beta}$ such that $\hat{F} (\hat{\alpha}) = \hat{F} (\hat{\beta})$. Clearly, $U \subseteq I$ and we want to show the other inclusion. Let $x \in U$. By the definition of category algebras, $x$ can be expresses uniquely as $\sum_{i=1} ^n \lambda_i \alpha_i$ where $\alpha_i$ are pairwise different morphisms in $\hat {\mathcal{E}}$ and $\lambda_i \in k$. Then $\varphi (\sum _{i=1}^n \lambda_i \alpha_i) = \sum_{i=1}^n \lambda_i \hat{F} (\alpha_i) = 0$. Those $\hat{F} (\alpha_i)$ are probably not pairwise different in $\mathcal{E}$. By changing the indices if necessary, we can write the set $\{ \alpha_i \}_{i=1}^n$ as a disjoint union of $l$ subsets: $\{ \alpha_1, \ldots, \alpha_{s_1} \}$, $\{\alpha_{s_1+1}, \ldots, \alpha_{s_2} \}$ and so on, until $\{\alpha_{s_{l-1}+1}, \ldots, \alpha_{s_l} \}$ such that two morphisms have the same image under $\hat{F}$ if and only if they are in the same set.

Now we have:
\begin{align*}
\varphi(x) = (\lambda_1 + \ldots + \lambda_{s_1}) \hat{F} (\alpha_{s_1}) + \ldots + (\lambda _{s_{l-1}+1} + \ldots + \lambda_{s_l}) \hat{F} (\alpha_{s_l}) = 0.
\end{align*}
Therefore,
\begin{align*}
\lambda_1 + \ldots + \lambda_{s_1} = \ldots = \lambda _{s_{l-1}+1} + \ldots + \lambda_{s_l} = 0,
\end{align*}
and hence
\begin{align*}
x &= [\lambda_2 (\alpha_2 - \alpha_1) + \ldots + \lambda_{s_1}(\alpha_{s_1} - \alpha_1)] + \ldots \\
& + [\lambda _{s_{l-1}+2} (\alpha_{s_{l-1}+2} - \alpha_{s_{l-1}+1}) + \ldots + \lambda_{s_l} (\alpha_{s_l} - \alpha_{s_{l-1}+1})]
\end{align*}
is contained in $U$.
\end{proof}

If $\mathcal{E}$ is a finite free EI category, we can put a \textit{length grading} on its morphisms as follows: automorphisms and unfactorizable morphisms are given grades 0 and 1 respectively; if $\alpha$ is a factorizable morphism, then it can be expressed (probably not unique) as a composite $\alpha_n \alpha_{n-1} \ldots \alpha_2 \alpha_1$ with all $\alpha_i$ unfactorizable and we assign $\alpha$ grade $n$. This grading is well defined by the Unique Factorization Property of finite free EI categories. It is clear that this length grading cannot be applied to an arbitrary finite EI category. We say a finite EI category can be \textit{graded} if this length grading is well defined on it. The following proposition gives us criterions to determine whether an arbitrary finite EI category can be graded.

\begin{proposition}
Let $\mathcal{E}$ be a finite EI category. Then the following are equivalent:
\begin{enumerate}
\item $\mathcal{E}$ is a graded finite EI category.
\item For each factorizable morphism $\alpha$ in $\mathcal{E}$, whenever it has two factorizations $\alpha_1 \circ \ldots \circ \alpha_m$ and $\beta_1 \circ \ldots \circ \beta_n$ into unfactorizable morphisms, we have $m=n$.
\item Let $\hat{\mathcal{E}}$ be the free EI cover of $\mathcal{E}$ and $\hat{F}: \hat{\mathcal{E}} \rightarrow \mathcal{E}$ be the covering functor. If two morphisms $\hat{\alpha}$ and $\hat{\beta}$ in $\hat{\mathcal{E}}$ have the same image under $\hat{F}$, then they have the same length in $\hat{\mathcal{E}}$.
\end{enumerate}
\end{proposition}

\begin{proof}
It is easy to see that if condition (2) holds, our grading works for $\mathcal{E}$, and hence (1) is true. Otherwise, if a factorizable morphism $\alpha$ has two decompositions $\alpha_n \circ \ldots \circ \alpha_1$ and $\beta_m \circ \ldots \circ \beta_1$ with $m \neq n$, then $\alpha$ should be assigned a grade $n$ by the first decomposition, and a grade $m$ by the second decomposition. Thus our grading cannot be applied to $\mathcal{E}$. This proves the equivalence of (1) and (2).

Now let $\alpha$ be an arbitrary morphism in $\mathcal{E}$ which has two different decompositions $\alpha_n \circ \ldots \circ \alpha_1$ and $\beta_m \circ \ldots \circ \beta_1$ into unfactorizable morphisms. Since $\hat{\mathcal{E}}$ is the free EI cover of $\mathcal{E}$, these unfactorizable morphisms are also unfactorizable morphisms in $\hat{\mathcal{E}}$. Let $\hat{\alpha}$ and $\hat{\beta}$ be the composite morphisms of these $\alpha_i$'s and $\beta_i$'s in $\hat{\mathcal{E}}$ respectively. Thus $\hat{\alpha} - \hat{\beta}$ is contained in $U$ since they have the same image $\alpha$ under $\hat{F}$. If (3) is true, then $m = n$ since $\hat{\alpha}$ and $\hat{\beta}$ have lengths $m$ and $n$ respectively. Therefore (3) implies (2). We can check that (2) implies (3) in a similar way.
\end{proof}

The following two lemmas are from \cite{Li}.

\begin{lemma}
Let $\mathcal{E}$ be a finite free EI category and $\alpha: x \rightarrow y$ be an unfactorizable morphism. Define $H = \Aut _{\mathcal{E}} (y)$ and $H_0 = \Stab _H (\alpha)$. If $| H_0 |$ is invertible in $k$, then the cyclic module $k\mathcal{E} \alpha$ is projective.
\end{lemma}

\begin{proof}
This is Lemma 5.2 of \cite{Li}, where we assumed that the automorphism groups of all objects are invertible in $k$ but only used the fact that $| H_0 |$ is invertible in $k$. Here we give a sketch of the proof. Let $e = \frac{1} { | H_0 |} \sum_{h\in H_0} h$. Then $e$ is well defined since $| H_0 |$ is invertible in $k$, and is an idempotent in $k\mathcal{E}$. Now define a map $\varphi: k\mathcal{E}e \rightarrow k\mathcal{E} \alpha$ by sending $re$ to $r\alpha$ for $r \in k\mathcal{E}$. We can check that $\varphi$ is an $k\mathcal{E}$-module isomorphism. Thus $k\mathcal{E} \alpha$ is projective. See \cite{Li} for a detailed proof.
\end{proof}

\begin{lemma}
Let $\mathcal{E}$ be a finite free EI category and $\alpha: x \rightarrow y$ and $\alpha': x' \rightarrow y'$ be two distinct unfactorizable morphisms in $\mathcal{E}$. Then $k\mathcal{E} \alpha \cap k\mathcal{E} \alpha' = 0$ or $k\mathcal{E} \alpha = k\mathcal{E} \alpha'$.
\end{lemma}

\begin{proof}
This is Lemma 5.1 of \cite{Li}. We give a sketch of the proof. Notice that $k\mathcal{E} \alpha$ is spanned by all morphisms of the form $\beta \alpha$ where $\beta: y \rightarrow z$ is a morphism starting at $y$. Similarly, $k\mathcal{E} \alpha'$ is spanned by all morphisms of the form $\beta' \alpha'$ where $\beta': y' \rightarrow z'$ is a morphism starting at $y'$. If $x \neq x'$ or $y \neq y'$, then by the Unique Factorization Property of finite free EI categories we conclude that the set $\mathcal{E} \alpha \cap \mathcal{E} \alpha' = \emptyset$, and the conclusion follows. If $x =x'$ and $y= y'$, then the set $\mathcal{E} \alpha$ coincides with the set $\mathcal{E} \alpha'$ if and only if there is an automorphism $h \in \Aut _{\mathcal{E}} (y)$ such that $h\alpha = \alpha'$ again by the UFP. Otherwise, we must have $\mathcal{E} \alpha \cap \mathcal{E} \alpha' = \emptyset$. The conclusion follows from this observation.
\end{proof}

\begin{remark}
The reader can check that the conclusion of Lemma 6.3 is true for any non-isomorphisms $\alpha$ in $\mathcal{E}$ by using the UFP. Moreover, a direct check shows that it is also true for automorphisms. Similarly, we can also prove that Lemma 6.4 still holds if we assume that $\alpha$ and $\alpha'$ are two morphisms with the same target and source.
\end{remark}

Theorem 5.8 has a corresponding version for finite EI categories.

\begin{proposition}
Let $\mathcal{E}$ be a graded finite EI category. Then $k \mathcal{E}$ is a Koszul algebra if and only if $k\mathcal{E}$ is a quasi-Koszul algebra and $\mathcal{E}$ is a standardly stratified category (in a sense defined in \cite{Webb}) with respect to the canonical partial order on $\Ob \mathcal{E}$.
\end{proposition}

\begin{proof}
By the decompositions
\begin{equation*}
k\mathcal{E} _0 = \bigoplus _{x \in \Ob \mathcal{E}} k \Aut _{\mathcal{E}} (x), \quad k\mathcal{E} _i = \bigoplus _{x \neq y \in \Ob \mathcal{E}} k {\mathcal{E}} (x, y)_i, i>0
\end{equation*}
we conclude that all $k \mathcal{E} _i$ are projective $k \mathcal{E} _0$-modules if and only if $k \mathcal{E} (x, y) _i$ are projective $k \Aut _{\mathcal{E}} (y)$-modules for all $i \geqslant 0$, $x, y \in \Ob \mathcal{E}$. Notice that $k \mathcal{E} (x, y) _i$ is spanned by morphisms from $x$ to $y$ with length $i$, and these morphisms form several orbits under the action of $\Aut _{\mathcal{E}} (y)$. Suppose that there are $n$ distinct orbits and take a representative $\alpha_j$ from each orbit. Then we have a decomposition $k \mathcal{E} (x, y) _i \cong \bigoplus _{j=1}^n k \Aut _{\mathcal{E}} (y) \alpha_j$. Thus $k \mathcal{E} (x, y) _i$ is a projective $k \Aut _{\mathcal{E}} (y)$-module if and only if each $k \Aut _{\mathcal{E}} (y) \alpha_j$ is a projective $k \Aut _{\mathcal{E}} (y)$-module, and if and only if the stabilizer of $\alpha_j$ in $\Aut _{\mathcal{E}} (y)$ has an order invertible in $k$. This happens if and only if $\mathcal{E}$ is standardly stratified by Theorem 2.5 in \cite{Webb}. In conclusion, all $k \mathcal{E} _i$ are projective $k \mathcal{E} _0$-modules if and only if $\mathcal{E}$ is standardly stratified in a sense defined in \cite{Webb}.

Notice that $k \mathcal{E} _0$ is the direct sum of several group algebras, and hence is self-injective. If $k \mathcal{E}$ is Koszul, then it is quasi-Koszul by Theorem 2.17. Moreover, all $k \mathcal{E} _i$ are projective $k \mathcal{E} _0$-modules (see the last paragraph of Section 2). Therefore, $\mathcal{E}$ is standardly stratified.

Conversely, if $k\mathcal{E}$ is standardly stratified, then all $k \mathcal{E} _i$ are projective $k \mathcal{E} _0$-modules. By Corollary 2.18, $k \mathcal{E}_0$ is a Koszul $k \mathcal{E}$-module if $k\mathcal{E}$ is quasi-Koszul.
\end{proof}

In the first paragraph of the above proof we have showed that a graded finite EI category $\mathcal{E}$ is standardly stratified in the sense of \cite{Webb} if and only if its $k$-linearization as a graded directed category is standardly stratified in our sense. Actually this is still true for an arbitrary finite EI category by comparing Theorem 5.7 in this paper and Theorem 2.5 in \cite{Webb}. This is not surprising since the $k$-linearization of $\mathcal{E}$ is precisely the associated category of the algebra $k\mathcal{E}$.

\begin{proposition}
Let $\mathcal{E}$ be a finite EI category which might not be graded. Then $\mathcal{E}$ is standardly stratified if and only if $M = \bigoplus _{x \in \Ob \mathcal{E}} \Aut _{\mathcal{E}} (x)$ viewed as a  $k\mathcal{E}$-module has finite projective dimension.
\end{proposition}

\begin{proof}
Consider the $k$-linearization of $\mathcal{E}$, which is a directed category. By Theorem 5.8 and the remark we made in the paragraph before this proposition, we conclude that $\pd M < \infty$.

Conversely, suppose that $\mathcal{E}$ is not standardly stratified. Then there is a non-isomorphism $\gamma: t \rightarrow y$ such that the order of $H_{\gamma}$ is not invertible in the field $k$ by Theorem 2.5 in \cite{Webb}, where $H = \Aut _{\mathcal{E}} (y)$ and $H_{\gamma} = \Stab _H (\gamma)$. For this object $y$, define $S$ to be the set of objects $w$ such that there is a non-isomorphism $\beta: w \rightarrow y$ satisfying that $| H_{\beta}|$ is not invertible in $k$. This set $S$ is nonempty since $t \in S$. It is a poset equipped with the partial order inherited from the canonical partial order on $\Ob \mathcal{E}$. Take a fixed object $z$ which is maximal in this set and define $I_{ > z} = \{ x \in \Ob \mathcal{E} \mid x > z \}$.

By our definition, for an arbitrary object $x \in I_{>z}$ and a non-isomorphism $\alpha: x \rightarrow y$ (if it exists), the group $H _{\alpha} \leqslant H$ has an order invertible in $k$. Therefore, the $kH$-module $kH \alpha$ is projective. Since the value of $k \mathcal{E} 1_x$ on $y$ is 0 or is spanned by all non-isomorphisms from $x$ to $y$, and these non-isomorphisms form a disjoint union of $H$-orbits, we conclude that the value of $k\mathcal{E} 1_x$ on $y$ is a projective $kH$-module (notice that we always view 0 as a zero projective module). With the same reasoning, we know that the value of $k\mathcal{E} 1_z$ on $y$ is not a projective $kH$-module.

Consider the $k\mathcal{E}$-module $L = k \Aut _{\mathcal{E}} (z)$. We claim that $\pd L = \infty$. If this is true, then $\pd M = \infty$ since $L$ is a direct summand of $M$. We prove this claim by showing the following statement: for each $i \geqslant 1$, every projective cover of $\Omega^i(L)$ is supported on $I_{ > z}$; the value $\Omega^i(L) (y)$ of $\Omega(L)$ on $y$ is non-zero and is not a projective $kH$-module. Clearly, $\Omega(L)$ is spanned by all non-isomorphisms starting from $z$ and is supported on $I_{ > z}$; $\Omega(L) (y)$, spanned by all non-isomorphisms from $z$ to $y$, is non-zero. Moreover, $\Omega(L)(y)$ coincides with the value of $k\mathcal{E} 1_z$ on $y$ and is not a projective $kH$-module. Therefore our statement is true for $i=1$.

Suppose that this statement is true for $n$, and let $P$ be a projective cover of $\Omega^n (L)$. The exact sequence
\begin{equation*}
\xymatrix{ 0 \ar[r] & \Omega^{n+1} (L) \ar[r] & P \ar[r] & \Omega^n (L) \ar[r] & 0}
\end{equation*}
gives rise to an exact sequence
\begin{equation*}
\xymatrix{ 0 \ar[r] & \Omega^{n+1} (L) (y) \ar[r] & P(y) \ar[r] & \Omega^n (L)(y) \ar[r] & 0.}
\end{equation*}

Let us focus on the above sequences. Since $\Omega^n (L)$ is supported on $I_{ > z}$, so are $P$ and $\Omega^{n+1} (L)$. By the induction hypothesis $\Omega^n (L) (y) \neq 0$, Thus $P(y) \neq 0$. But $P$ is supported on $I_{ > z}$, so $P \in \text{add} (\bigoplus _{x \in I_{ > z}}  k\mathcal{E} 1_x)$. Notice that the value of each $k \mathcal{E} 1_x$ on $y$ is zero or a nontrivial projective $kH$-module. Therefore, $P(y)$ is a projective $kH$-module. Again by the induction hypothesis, $\Omega^n (L) (y)$ is not a projective $kH$-module, so $\Omega^n (L) (y) \ncong P(y)$, and $\Omega^{n+1} (L) (y)$ is non-zero. It cannot be a projective $kH$-module. Otherwise, $\Omega^{n+1} (L) (y)$ is also an injective $kH$-module and hence the above sequence splits, so $\Omega^n (L) (y)$ as a summand of $P(y)$ is a projective $kH$-module, too. But this contradicts the induction hypothesis.

We proved the induction hypothesis for $\Omega^{n+1} (L)$. Thus our statement and claim are proved. Consequently, $\pd M = \infty$.
\end{proof}

Now we can prove:

\begin{theorem}
Let $\mathcal{E}$ be a finite free EI category. Then the following are equivalent:
\begin{enumerate}
\item $\pd k\mathcal{E}_0 \leqslant 1$;
\item $k \mathcal{E}$ is a Koszul algebra;
\item $\mathcal{E}$ is standardly stratified;
\item $\pd k\mathcal{E}_0 < \infty$.
\end{enumerate}
\end{theorem}

\begin{proof}
It is clear that $\pd k\mathcal{E}_0 = 0$ if and only if $\mathcal{E}$ is a finite EI category with a single object since we only consider connected categories. In this situation, $k \mathcal{E} = k \mathcal{E}_0$, $J = 0$, and all statements are trivially true. Thus without loss of generality we suppose that $\pd k\mathcal{E}_0 \neq 0$.

Observe that $\pd k\mathcal{E} _0 =1$ if and only if $\Omega( k\mathcal{E} _0) = J = \bigoplus _{i \geqslant 1} k\mathcal{E}_i$ is projective. Since $J$ is spanned by all non-isomorphisms in $\mathcal{E}$ and each non-isomorphism can be written as a composition of unfactorizable morphisms, it is generated in degree 1. Thus $k \mathcal{E} _0$ is a Koszul $k\mathcal{E}$-module, and (1) implies (2). By Proposition 6.6, (2) implies (3). The statements (3) and (4) are equivalent by Proposition 6.7.

Now we prove that (3) implies (1). If $\mathcal{E}$ is standardly stratified, then for every morphism $\alpha : x \rightarrow y$ in $\mathcal{E}$, the order of $\Stab _H (\alpha)$ is invertible in $k$, where $H = \Aut _{\mathcal{E}} (y)$. By Lemma 6.4, $J$ is a direct sum of some $k\mathcal{E}$-modules $k\mathcal{E} \alpha_i$'s with each $\alpha_i$ unfactorizable. By Lemma 6.3, each $k\mathcal{E} \alpha_i$ is projective. Therefore, $J$ is also projective, i.e., $\pd k\mathcal{E}_0 = 1$.
\end{proof}

This theorem and Theorem 5.12 give us a way to construct Koszul algebras in the classical sense. Indeed, let $\mathcal{E}$ be a standardly stratified finite free EI category and define $\mathcal{D}$ to be the subcategory formed by removing all non-identity automorphisms. Then by Theorem 5.12 $k\mathcal{D}$ is a Koszul algebra in the classical sense since $k\mathcal{E}_0$ is a Koszul $k \mathcal{E}$-module in the generalized sense by the previous theorem.

Let us get more information about the projective resolutions of $k\mathcal{E} _0$ for arbitrary finite free EI categories. In general, $\Omega( k\mathcal{E} _0) \cong J = \bigoplus _{i \geqslant 1} k\mathcal{E}_i$ is not projective, but it is still a direct sum of some $k\mathcal{E}$-modules $k\mathcal{E} \alpha$'s with each $\alpha$ unfactorizable by Lemma 6.4. Thus the projective resolutions of $k\mathcal{E}_0$ is completely determined by the projective resolutions of those $k\mathcal{E} \alpha$'s.

\begin{lemma}
Let $\mathcal{E}$ be a finite free EI category and $\alpha: x \rightarrow y$ be an unfactorizable morphism. Grade the $k \mathcal{E}$-module $k\mathcal{E} \alpha$ by putting $\alpha$ in degree 1, namely, $(k \mathcal{E} \alpha)_1 = k \Aut _{\mathcal{E}} (y) \alpha$. Then $\Omega (k \mathcal{E} \alpha)$ is 0 or is generated in degree 1, and $\Omega (k \mathcal{E} \alpha)_1 = \Omega (k \mathcal{E} \alpha) (y)$, the value of $\Omega (k \mathcal{E} \alpha)$ on $y$.
\end{lemma}

\begin{proof}
Let $H = \Aut _{\mathcal{E}} (y)$ and $H_0 = \Stab_H (\alpha)$. If $| H_0|$ is invertible in $k$, then by Lemma 6.3, $k\mathcal{E} \alpha$ is a projective $k\mathcal{E}$-module, so $\Omega^i( k \mathcal{E} \alpha) =0$ for all $i \geqslant 1$, in particular $\Omega (k \mathcal{E} \alpha) =0$. The conclusion is trivially true. Thus we only need to deal with the case that $|H_0|$ is not invertible.

Consider the projective presentation
\begin{equation*}
\xymatrix {0 \ar[r] & N \ar[r] & k\mathcal{E} 1_y [1] \ar[r]^p & k \mathcal{E} \alpha \ar[r] & 0}
\end{equation*}
where $p$ maps $1_y$ to $\alpha$. Since $\Omega (k \mathcal{E} \alpha)$ is isomorphic to a direct summand of $N$, it is enough to show that $N$ is generated in degree 1, and $N_1 = N(y)$.

Notice that $k\mathcal{E} 1_y$ is spanned by all morphisms in $\mathcal{E}$ with source $y$, and $k\mathcal{E} \alpha$ is spanned by all morphisms in $\mathcal{E}$ of the form $\beta \alpha$ where $\beta$ is a morphism in $\mathcal{E}$ with source $y$. We claim that $N$ is spanned by vectors of the form $\beta_1 - \beta_2$ with $\beta_1 \alpha = \beta_2 \alpha$, where $\beta_1$ and $\beta_2$ are two morphisms with source $y$.

Clearly, every such difference is contained in $N$. Conversely, let $v \in N$. Then $v$ can be written as $\sum _{i=1} ^n \lambda_i \beta_i$ such that $\lambda_i \in k$ and $\beta_i$ are pairwise different morphisms with source $y$. By the definition of $p$, $\sum _{i=1} ^n \lambda_i \beta_i \alpha = 0$. Those $ \beta_i \alpha$ might not be pairwise different in $\mathcal{E}$. Now we apply the same technique used in the proof of Lemma 6.1. By changing the indices if necessary, we can group the same morphisms together and suppose that $ \beta_1 \alpha = \ldots = \beta_{s_1} \alpha$, $ \beta_{s_1+1} \alpha = \ldots = \beta_{s_2} \alpha$ and so on, until $ \beta_{s_{l-1}+1} \alpha = \ldots = \beta_{s_l} \alpha$.

We have:
\begin{align*}
p(v) = (\lambda_1 + \ldots + \lambda_{s_1}) \beta_{s_1} \alpha + \ldots + (\lambda _{s_{l-1}+1} + \ldots + \lambda_{s_l}) \beta_{s_l} \alpha = 0.
\end{align*}
Therefore,
\begin{align*}
\lambda_1 + \ldots + \lambda_{s_1} = \ldots = \lambda _{s_{l-1}+1} + \ldots + \lambda_{s_l} = 0,
\end{align*}
and hence
\begin{align*}
v &= [\lambda_2 (\beta_2 - \beta_1) + \ldots + \lambda_{s_1}(\beta_{s_1} - \beta_1)] + \ldots \\
& + [\lambda _{s_{l-1}+2} (\beta_{s_{l-1}+2} - \beta_{s_{l-1}+1}) + \ldots + \lambda_{s_l} (\beta_{s_l} - \beta_{s_{l-1}+1})].
\end{align*}
So $v$ can be written as a sum of these differences.

Now we can prove the lemma. Take an arbitrary object $z \in \Ob \mathcal{E}$ and consider the value $N (z)$. If it is 0, the conclusion holds trivially. Suppose that $N (z) \neq 0$. By the above description, $N (z)$ is spanned by vectors $\beta_1 - \beta_2$ such that $\beta_1, \beta_2$ are two morphisms from $y$ to $z$, and $\beta_1 \alpha = \beta_2 \alpha$. By the equivalent definition of UFP described in Remark 6.5, there is an automorphism $h \in \Aut _{\mathcal{E}} (y)$ such that $\beta_1 = \beta_2 h$ and $\alpha = h^{-1} \alpha$. Therefore $h\alpha = \alpha$, and $1-h \in N (y)$. Thus $\beta_1 - \beta_2 = \beta (1-h) \in k\mathcal{E} \cdot N (y)$. Since $z$ is taken to be an arbitrary object, $N$ is generated by $N (y)$, which is clearly equal to $N_1$.
\end{proof}

From this lemma we can get:

\begin{proposition}
Let $\mathcal{E}$ be a finite free EI category, then $\Ext _{k\mathcal{E}} ^2 (k\mathcal{E}_0, k\mathcal{E}_0) = 0$.
\end{proposition}

\begin{proof}
Since $\Omega (k \mathcal{E} _0) \cong J$, it is enough to show $\Ext _{k \mathcal{E}} ^1 (J, k \mathcal{E} _0) =0$. The conclusion holds trivially if $J$ is projective. Otherwise, since $J$ is the direct sum of some $k\mathcal{E} \alpha$'s with $\alpha$ unfactorizable, by the above lemma we know that $\Omega J$ is generated in degree 1.

Applying the functor $\Hom _{k\mathcal{E}} (-, k\mathcal{E} _0)$ to $0 \rightarrow \Omega J \rightarrow P \rightarrow J \rightarrow 0$ we get
\begin{equation*}
0 \rightarrow \Hom _{k\mathcal{E} } (J, k\mathcal{E}_0) \rightarrow \Hom _{k\mathcal{E} } (P, k\mathcal{E} _0) \rightarrow \Hom _{k\mathcal{E}} (\Omega J, k\mathcal{E}_0) \rightarrow \Ext _{k\mathcal{E}} ^1 (J, k\mathcal{E}_0) \rightarrow 0.
\end{equation*}
Since all modules are generated in degree 1, the sequence
\begin{equation}
0 \rightarrow \Hom _{k\mathcal{E} } (J, k\mathcal{E}_0) \rightarrow \Hom _{k\mathcal{E}} (P, k\mathcal{E} _0) \rightarrow \Hom _{k\mathcal{E}} (\Omega J, k\mathcal{E}_0)
\end{equation}
is isomorphic to the sequence
\begin{equation*}
0 \rightarrow \Hom _{k\mathcal{E} _0} (J_1, k\mathcal{E}_0) \rightarrow \Hom _{k\mathcal{E} _0} (P_1, k\mathcal{E} _0) \rightarrow \Hom _{k\mathcal{E} _0} ((\Omega J)_1, k\mathcal{E}_0)
\end{equation*}
obtained by applying the exact functor $\Hom _{k \mathcal{E}_0} (-, k\mathcal{E} _0)$ to the exact sequence $0 \rightarrow (\Omega J)_1 \rightarrow P_1 \rightarrow J_1 \rightarrow 0$. Thus the last map in sequence 6.1 is surjective, so $\Ext _{k \mathcal{E}} ^1 (J, k \mathcal{E} _0) =0$.
\end{proof}

The fact that $\Omega J$ is generated in degree 1 implies $\Ext^2 _{k\mathcal{E}} (k\mathcal{E}_0, k\mathcal{E}_0) =0$. Actually the converse statement is also true. Indeed, consider the exact sequence $0 \rightarrow \Omega J \rightarrow P \rightarrow J \rightarrow 0$. If $\Ext^2 _{k\mathcal{E}} (k\mathcal{E}_0, k\mathcal{E}_0) =0$, applying the exact functor $\Hom _{k\mathcal{E} _0} (-, k\mathcal{E}_0)$ we get the exact sequence
\begin{equation*}
0 \rightarrow \Hom _{k\mathcal{E} _0} (J, k\mathcal{E} _0) \rightarrow \Hom _{k\mathcal{E} _0} (P, k\mathcal{E} _0) \rightarrow \Hom _{k\mathcal{E} _0} (\Omega J, k\mathcal{E} _0) \rightarrow 0,
\end{equation*}
which is isomorphic to
\begin{equation*}
0 \rightarrow \Hom _{k\mathcal{E} _0} (J_1, k\mathcal{E} _0) \rightarrow \Hom _{k\mathcal{E} _0} (P_1, k\mathcal{E} _0) \rightarrow \Hom _{k\mathcal{E} _0} (\Omega J / J (\Omega J), k\mathcal{E} _0) \rightarrow 0
\end{equation*}
since both $J$ and $P$ are generated in degree $1$. Applying the functor $\Hom _{ k\mathcal{E} _0} (-, k \mathcal{E}_0 )$ again, we recover $0 \rightarrow \Omega J / J (\Omega J) \rightarrow P_1 \rightarrow J_1 \rightarrow  0$. Therefore, $\Omega J / J (\Omega J) \cong (\Omega J)_1$, so $\Omega J$ is generated in degree 1.

Finite free EI categories with quasi-Koszul category algebras have very special homological properties. For example:

\begin{proposition}
Let $\mathcal{E}$ be a finite free EI category. Then the following are equivalent:
\begin{enumerate}
\item $\Ext _{k\mathcal {E}} ^i (k \mathcal{E} _0, k\mathcal{E} _0) = 0$ for all $i \geqslant 2$;
\item for every unfactorizable morphism $\alpha: x \rightarrow y$ and $i \geqslant 0$, either $\Omega ^i (k \mathcal{E} \alpha) $ are all 0, or they are all generated in degree 1 (in which case it is generated by $\Omega ^i (k \mathcal{E} \alpha) (y)$);
\item $k \mathcal{E}$ is a quasi-Koszul algebra.
\end{enumerate}
\end{proposition}

\begin{proof}
If $k \mathcal{E}$ is a quasi-Koszul algebra, then
\begin{equation*}
\Ext _{k\mathcal {E}} ^i (k \mathcal{E} _0, k\mathcal{E} _0) = \Ext _{k\mathcal{E}} ^2 (k \mathcal{E} _0, k\mathcal{E} _0) \cdot \Ext _{k\mathcal{E}} ^{i-2} (k \mathcal{E} _0, k\mathcal{E} _0)
\end{equation*}
for every $i \geqslant 2$. But $\Ext _{k\mathcal{E}} ^2 (k \mathcal{E} _0, k\mathcal{E} _0) =0$ by Proposition 6.10, so (3) implies (1). Clearly, (1) implies (3).

Notice that $k \mathcal{E} \alpha$ is a isomorphic to a direct summand of $J \cong \Omega (k \mathcal{E}_0)$. Thus we only need to prove the equivalence of the following two statements:
\begin{enumerate}[(1')]
\item $\Ext _{k\mathcal {E}} ^i (J, k\mathcal{E} _0) = 0$ for every $i \geqslant 1$;
\item $\Omega ^i (J) =0 $ or is generated in degree 1 for every $i \geqslant 1$.
\end{enumerate}

Since the technique we use is similar to that in the proof of Proposition 6.10, we only give a sketched proof. In the case that $J$ is projective, i.e., $\mathcal{E}$ is standardly stratified, then (1') and (2') are trivially true, hence they are equivalent. Now suppose that $J$ is not projective. From the proof of Proposition 6.10 and the paragraph after it we conclude that $\Omega J$ is generated in degree 1 if and only if $\Ext _{k\mathcal{E}} ^1 (J, k\mathcal{E} _0) =0$. Replacing $J$ by $\Omega J$ (which is also generated in degree 1 either by the induction hypothesis or by the hypothesis $\Ext _{k\mathcal{E}} ^1 (J, k\mathcal{E} _0) =0$) and using the same technique, we get $\Omega ^2 (J)$ is generated in degree 1 if and only if $\Ext _{k\mathcal{E}} ^2 (J, k\mathcal{E} _0) =0$. The equivalence of (1') and (2'), and hence the equivalence of (1) and (2), come from induction.
\end{proof}

The reader may guess that the category algebra of a finite free EI category is always quasi-Koszul in our sense because of the following reasons: Finite free EI categories generalize finite groups and acyclic quivers, for which the associated algebras are all quasi-Koszul;  by Proposition 6.6 and Theorem 6.8, for an arbitrary finite EI category $\mathcal{E}$, $k\mathcal{E}$ is Koszul if $\mathcal{E}$ is standardly stratified and one of the following condition holds: $k \mathcal{E}$ is quasi-Koszul, or $\mathcal{E}$ is a finite free EI category; and we have proved that $\Ext _{k\mathcal{E}} ^2 (k\mathcal{E} _0, k\mathcal{E} _0) =0$ if $\mathcal{E}$ is a finite free EI category. Unfortunately, this conjecture is false, as shown by the following example.

\begin{example}
Let $\mathcal{E}$ be the following finite EI category where: $\Aut _{\mathcal{E}} (x) = \langle 1_x \rangle$, $\Aut _{\mathcal{E}} (z) = \langle 1_z \rangle$, $\Aut _{\mathcal{E}} (y) = \langle h \rangle$ is a group of order 2; $\mathcal{E} (x, y) = \{ \alpha \}$, $\mathcal{E} (y, z) = \{ \beta \}$ and $\mathcal{E} (x, z) = \{ \beta \alpha \}$. The reader can check that $\mathcal{E}$ is a finite free EI category and then the length grading can be applied on it. Let $k$ be an algebraically closed field with characteristic 2.

\begin{equation*}
\xymatrix { x \ar[r] ^{\alpha} & y \ar[r] ^{\beta} & z.}
\end{equation*}

The indecomposable direct summands of $k \mathcal{E}$ and $k \mathcal{E}_0$ are:
\begin{equation*}
P_x = \begin{matrix}   x_0 \\ y_1 \\ z_2 \end{matrix}, \qquad P_y = \begin{matrix} & y_0 & \\ y_0 & & z_1 \end{matrix}, \qquad P_z = z_0, \qquad k \mathcal{E} _0 \cong x_0 \oplus z_0 \oplus \begin{matrix} y_0 \\ y_0 \end{matrix}.
\end{equation*}
We use indices to mark the degrees of composition factors. The reader should bear in mind that the two simple modules $y$ appearing in $P_y$ have the same degree.

Take the summand $x_0$ of $k\mathcal{E} _0$. By computation, we get
\begin{equation*}
\Omega (x_0) = \begin{matrix}   y_1 \\ z_2 \end{matrix}, \qquad \Omega ^2 (x_0) = \begin{matrix} y_1 \end{matrix}, \qquad \Omega ^3 (x_0) = y_1 \oplus z_2.
\end{equation*}
Applying $\Hom _{k\mathcal{E}} (-, k\mathcal{E} _0)$ to the exact sequence
\begin{equation*}
\xymatrix{0 \ar[r] & \Omega^3 (x_0) \ar[r] & P_y [1] \ar[r] & \Omega^2 (x_0) \ar[r] &0}
\end{equation*}
we get $\Ext^3 _{k\mathcal{E}} (k\mathcal{E}_0, k\mathcal{E}_0 ) \neq 0$. Consequently, $k \mathcal{E}$ is not a quasi-Koszul algebra in our sense by the previous proposition.
\end{example}

We aim to characterize finite free EI categories with quasi-Koszul category algebras. For this goal, we make the following definition:

\begin{definition}
Let $\mathcal{E}$ be a finite EI category. An object $x \in \Ob \mathcal{E}$ is called left regular if for every morphism $\alpha$ with target $x$, the stabilizer of $\alpha$ in $\Aut _{\mathcal{E}} (x)$ has an order invertible in $k$. Similarly, $x$ is called right regular if for every morphism $\beta$ with source $x$, the stabilizer of $\beta$ in $\Aut _{\mathcal{E}} (x)$ has an order invertible in $k$.
\end{definition}

\begin{remark}
We make some comments for this definition.
\begin{enumerate}
\item If $x \in \Ob \mathcal{E}$ is maximal, i.e., there is no non-isomorphisms with source $x$, then $x$ is right regular by convention; similarly, if $x$ is minimal, then it is trivially left regular.
\item The category $\mathcal{E}$ is standardly stratified if and only if every object $x \in \Ob \mathcal{E}$ is left regular; similarly, $\mathcal{E} ^{\textnormal{op}}$ is standardly stratified if and only if every object $x \in \Ob \mathcal{E}$ is right regular
\item If $\mathcal{E}$ is a finite free EI category and $x \in \Ob \mathcal{E}$. Then $x$ is left regular if and only if for every $\alpha$ with target $x$, the $k\mathcal{E}$-module $k \mathcal{E} \alpha$ is a left projective $k \mathcal{E}$-module. Similarly, $x$ is right regular if and only if for every $\beta$ with source $x$, the right $k\mathcal{E}$-module $\beta (k \mathcal{E})$ is a right projective $k \mathcal{E}$-module.
\end{enumerate}
\end{remark}

\begin{lemma}
Let $\mathcal{E}$ be a finite free EI category and $\beta: x \rightarrow y$ be a morphism with $x \in \Ob \mathcal{E}$ right regular. Then there exists some idempotent $e$ in $k \mathcal{E}$ such that $\beta (k \mathcal{E}) \cong e (k \mathcal{E})$ as right $k \mathcal{E}$-modules by sending $e$ to $\beta$. In particular, $\beta (k G) \alpha \cong e (kG) \alpha$ as vector spaces for every morphism $\alpha$ with target $x$, where $G = \Aut _{\mathcal{E}} (x)$.
\end{lemma}

\begin{proof}
Let $G_0 = \Stab _G (\alpha)$ and $e = \sum _{g \in G_0} g / | G_0 |$. This is well defined since $x$ is right regular. Then we can prove $\beta (k \mathcal{E}) \cong e (k \mathcal{E})$ as right $k\mathcal{E}$-modules in a way similar to the proof of Lemma 6.3. The isomorphism is given by sending $er$ to $\beta r$ for $r \in k\mathcal{E}$. Since the image of $e (kG) \alpha \subseteq k\mathcal{E}$ is exactly $\beta (k \mathcal{E}) \alpha$, we deduce that $e (kG) \alpha \cong \beta (kG) \alpha$ as vector spaces.
\end{proof}

Using these concepts, we can get a sufficient condition for the category algebra of a finite free EI category to be quasi-Koszul.

\begin{theorem}
Let $\mathcal{E}$ be a finite free EI category such that every object $x \in \Ob \mathcal{E}$ is either left regular or right regular. Then $k \mathcal{E}$ is quasi-Koszul.
\end{theorem}

\begin{proof}
By the second statement of Proposition 6.11, it is enough to show that for each unfactorizable $\alpha: x \rightarrow y$ and every $i \geqslant 1$, $\Omega ^i (k \mathcal{E} \alpha)$ is 0 or generated by $\Omega ^i (k \mathcal{E} \alpha) (y)$. Let $H = \Aut _{\mathcal{E}} (y)$ and $H_0 = \Stab _H (\alpha)$. If $| H_0 |$ is invertible in $k$, then $k \mathcal{E} \alpha$ is a projective $k\mathcal{E}$-module, and the conclusion follows. So we only need to deal with the case that the order of $H_0$ is not invertible in $k$.

By Lemma 6.9, $\Omega (k \mathcal{E} \alpha)$ is generated in degree 1, or equivalently, generated by its value $\Omega (k \mathcal{E} \alpha) (y) = 1_y \Omega (k \mathcal{E} \alpha)$ on $y$. Now suppose that $\Omega ^i (k \mathcal{E} \alpha)$ is also generated in degree 1, or equivalently, generated by its value $\Omega ^i (k \mathcal{E} \alpha) (y) = 1_y \Omega ^i (k \mathcal{E} \alpha)$ on $y$, where $i \geqslant 1$. We claim that $\Omega ^{i+1} (k \mathcal{E} \alpha)$ is generated by $\Omega ^{i+1} (k \mathcal{E} \alpha) (y)$, which is clearly equal to $\Omega ^{i+1} (k \mathcal{E} \alpha)_1$. If this is true, then conclusion follows from Proposition 6.11.

Take an arbitrary object $z \in \Ob \mathcal{E}$ such that $\mathcal{E} (y, z) \neq \emptyset$. (In the case $\mathcal{E} (y, z) = \emptyset$, $\Omega ^s (k \mathcal{E} \alpha) (z) =0$ for $s \geqslant 0$, and the claim is trivially true.) The morphisms in $\mathcal{E} (y, z)$ form a disjoint union of orbits under the right action of $H$. By taking a representative $\beta_i$ from each orbit we have $\mathcal{E} (y, z) = \bigsqcup _{i=1} ^n \beta_i H$. Since $|H _0|$ is not invertible, $y$ is not left regular. By the assumption, $y$ must be right regular. Therefore, by the previous lemma, for each representative morphism $\beta_s$, $1 \leqslant s \leqslant n$, there exist some idempotent $e_i$ such that $\beta_s (k \mathcal{E}) \cong e_s (k \mathcal{E})$ as right projective $k \mathcal{E}$-modules, and $\beta_s (k \mathcal{E}) \alpha \cong e_s (k \mathcal{E}) \alpha$ as vector spaces.

Consider the exact sequence
\begin{equation*}
\xymatrix{ 0 \ar[r] & \Omega ^{i+1} (k \mathcal{E} \alpha) \ar[r] & P^i \ar[r] & \Omega^i (k \mathcal{E} \alpha) \ar[r] & 0,}
\end{equation*}
where we assume inductively that $\Omega ^i (k \mathcal{E} \alpha)$ is generated in degree 1, or equivalently generated by its value on $y$. Thus $P^i \in \text{add} (k \mathcal{E} 1_y [1])$. Observe that the segment of a minimal projective resolution of the $k\mathcal{E}$-module $k \mathcal{E} \alpha$
\begin{equation*}
\xymatrix{ P^{i+1} \ar[r] & P^i \ar[r] & \ldots \ar[r] & P^0 \ar[r] & k\mathcal{E} \alpha \ar[r] & 0}
\end{equation*}
induces a minimal projective resolution of the $kH$-module $kH \alpha$:
\begin{equation*}
\xymatrix{ P^{i+1} (y) \ar[r] & P^i (y) \ar[r] & \ldots \ar[r] & P^0 (y) \ar[r] & kH \alpha \ar[r] & 0.}
\end{equation*}
Thus $\Omega ^j (k\mathcal{E} \alpha)_1 = \Omega^j (k\mathcal{E} \alpha) (y) = \Omega ^j _{kH} (kH \alpha)$ for $1 \leqslant j \leqslant i+1$.

Applying the exact functor $\Hom _{k \mathcal{E}} (k \mathcal{E} 1_y, -)$ to the exact sequence
\begin{equation}
\xymatrix {0 \ar[r] & \Omega ^{i+1} (k \mathcal{E} \alpha) \ar[r] & P^i \ar[r] & \Omega^i (k \mathcal{E} \alpha) \ar[r] & 0},
\end{equation}
we get an exact sequence
\begin{equation*}
\xymatrix {0 \ar[r] & \Omega ^{i+1} (k \mathcal{E} \alpha) (y) \ar[r] & 1_y P^i \ar[r] & \Omega^i (k \mathcal{E} \alpha) (y) \ar[r] & 0},
\end{equation*}
which can be identified with
\begin{equation*}
\xymatrix {0 \ar[r] & \Omega ^{i+1} _{kH} (kH \alpha) \ar[r] & P^i(y) \ar[r] & \Omega^i _{kH} (kH \alpha) \ar[r] & 0}.
\end{equation*}
Applying the exact functor $\Hom_{kH} (\bigoplus _{s=1}^n kHe_s, -)$ to the above sequence, we have another exact sequence
\begin{equation}
0 \rightarrow \bigoplus _{s=1}^n e_s \Omega ^{i+1} _{kH} (kH \alpha) \rightarrow \bigoplus _{s=1}^n e_s P^i(y) \rightarrow \bigoplus _{s=1}^n e_s \Omega^i _{kH} (kH \alpha) \rightarrow 0.
\end{equation}

Since $\Omega^i (k \mathcal{E} \alpha)$ is generated by $\Omega ^i (k \mathcal{E} \alpha) (y) = \Omega^i _{kH} (kH \alpha)$ by the induction hypothesis, the value of $\Omega ^i (k\mathcal{E} \alpha )$ on $z$ is $\sum _{s=1}^n \beta_s \cdot \Omega_{kH} ^i (kH \alpha)$ (this is well defined as $\Omega ^i _{kH} (kH \alpha) \subseteq (kH) ^{\oplus m}$ for some $m \geqslant 0$). We check that this sum is actually direct by the UFP of $\mathcal{E}$. In conclusion,
\begin{equation}
\Omega ^i (k\mathcal{E} \alpha ) (z) = \bigoplus _{s=1}^n \beta_s \cdot \Omega_{kH} ^i (kH \alpha) \cong \bigoplus _{s=1}^n e_s \Omega_{kH} ^i (kH \alpha).
\end{equation}
Similarly, the value of $P^i$ on $z$ is
\begin{equation}
P^i (z) = \bigoplus _{s=1}^n \beta_s \cdot P^i(y) \cong \bigoplus _{s=1}^n e_s P^i(y).
\end{equation}

Restricted to $z$, sequence 6.2 gives rise to
\begin{equation}
\xymatrix {0 \ar[r] & \Omega ^{i+1} (k \mathcal{E} \alpha) (z) \ar[r] & P^i(z) \ar[r] & \Omega^i (k \mathcal{E} \alpha)(z) \ar[r] & 0}.
\end{equation}

On one hand, $\bigoplus _{s=1}^n \beta_s \Omega ^{i+1} _{kH} (kH \alpha) \subseteq \Omega ^{i+1} (k \mathcal{E} \alpha) (z)$. On the other hand, we have:
\begin{align*}
& \dim_k \bigoplus _{s=1}^n \beta_s \Omega ^{i+1} _{kH} (kH \alpha) = \dim_k \bigoplus _{s=1}^n e_s \Omega ^{i+1} _{kH} (kH \alpha) \quad \text{ by Lemma 6.15} \\
& = \dim_k \bigoplus _{s=1}^n e_s P^i(y) - \dim_k \bigoplus _{s=1}^n e_s \Omega^i _{kH} (kH \alpha) \quad \text{ by sequence 6.3} \\
& = \dim_k P^i (z) - \dim_k \Omega ^i (k\mathcal{E} \alpha ) (z) \quad \text{ by identities 6.4 and 6.5} \\
& = \dim_k \Omega ^{i+1} (k \mathcal{E} \alpha) (z) \quad \text{ by sequence 6.6.}
\end{align*}
Therefore, $\Omega ^{i+1} (k\mathcal{E} \alpha ) (z) = \bigoplus _{s=1}^n \beta_s \Omega ^{i+1} _{kH} (kH \alpha) = \bigoplus _{s=1}^n \beta_s \Omega ^{i+1} _{k \mathcal{E}} (k \mathcal{E} \alpha) (y)$ since $\Omega ^{i+1} _{kH} (kH \alpha) = \Omega ^{i+1} _{k \mathcal{E}} (k \mathcal{E} \alpha) (y)$. That is, the value of $\Omega ^{i+1} (k \mathcal{E} \alpha)$ on $z$ is generated by $\Omega ^{i+1} (k \mathcal{E} \alpha) (y)$. Since $z$ is arbitrary, our claim holds, and the conclusion follows from induction.
\end{proof}

\section{Standardly Stratified Algebras with Linear Standard Modules}

Theorem 5.10 tells us that the Yoneda category $E(\mathcal{C}_0)$ of a directed Koszul category $\mathcal{C}$ is still a directed Koszul category, so is standardly stratified as well. Moreover, the homological dual functor $E$ interchanges standard modules and indecomposable projective modules. Let $A$ be a Koszul algebra which is standardly stratified with respect to a poset of orthogonal primitive idempotents $(\{ e_{\lambda} \} _{\lambda \in \Lambda}, \leqslant)$. We may ask a similar question: is the Koszul dual algebra $\Gamma = \Ext _A ^{\ast} (A_0, A_0)$ standardly stratified with respect to $(\{ e_{\lambda} \} _{\lambda \in \Lambda}, \leqslant)$ (or $(\{ e_{\lambda} \} _{\lambda \in \Lambda}, \leqslant ^{\textnormal{op}})$) as well? (Here we identify the primitive idempotents of $A$ and $\Gamma$ in the following way: let $e$ be a primitive idempotent of $A$. Then it is also a primitive idempotent of $A_0$. Therefore, $A_0 e$ is a projective $A_0$-module, and $\Ext _A ^{\ast} (A_0e, A_0)$ is an indecomposable summand of $\Gamma$. This summand corresponds to a primitive idempotent of $\Gamma$, which we still denote by $e$.) This question has been studied in \cite{Agoston1,Agoston2,Drozd,Mazorchuk2,Mazorchuk3}. However, in all these papers $A_0$ is supposed to be a semisimple algebra. By modifying the technique used in \cite{Agoston1}, we get a sufficient condition for the Yoneda algebra $\Gamma$ to be standardly stratified with respect to the opposite order.

Throughout this section $A$ is a graded finite-dimensional basic $k$-algebra with $A_0$ self-injective. We choose a complete set of orthogonal primitive idempotents $\{ e_{\lambda} \} _{\lambda \in \Lambda}$ and let $\leqslant$ be a partial order on this set.

\begin{theorem}
If $A$ is standardly stratified with respect to $\leqslant$ such that all standard modules are concentrated in degree 0 and Koszul. Then $A_0 \cong \Delta$ and $\Gamma = \Ext _A ^{\ast} (A_0, A_0)$ is standardly stratified with respect to the poset $(\{ e_{\lambda} \} _{\lambda \in \Lambda}, \leqslant ^{\textnormal{op}})$, where $\Delta$ is the direct sum of all standard modules.
\end{theorem}

We show the first statement since it is relatively easier and leave the proof of the second statement to the end of this section. Let $\Delta_{\lambda}$ be a standard module with graded projective cover $P_{\lambda}$. Since $\Delta_{\lambda}$ is a Koszul $A$-module and concentrated in degree 0, $\Delta_{\lambda} = (\Delta_{\lambda})_0$ is a projective $A_0$-module by Corollary 2.5. The surjection $P_{\lambda} \rightarrow \Delta_{\lambda}$ induces a surjection $(P_{\lambda})_0 \rightarrow \Delta_{\lambda}$, so $\Delta_{\lambda}$ is a summand of $(P_{\lambda})_0$, and hence is isomorphic to $(P_{\lambda})_0$ since $(P _{\lambda}) _0$ is indecomposable. Put all these standard modules together we find $A_0 \cong \Delta$.

Take a minimal element $\mu \in \Lambda$ and let $e = e_{\mu}$. Let $\Lambda_1 = \Lambda \setminus \{ \mu \}$ and $\epsilon = \sum _{\lambda \in \Lambda_1} e_{\lambda}$. Viewed as an idempotent of $\Gamma$, $e$ is maximal with respect to $\leqslant ^{\textnormal{op}}$. The basic idea to prove the second statement is to show that $\Gamma e \Gamma$ is a projective $\Gamma$-module and the quotient algebra $\Gamma / \Gamma e \Gamma$ is standardly stratified with respect to the poset $(\{e_{\lambda} \} _{\lambda \in \Lambda_1}, \leqslant ^{\textnormal{op}})$. Then the conclusion follows from induction.

We collect a list of preliminary results in the following lemmas, where the algebra $A$ is the same as in Theorem 7.1 if we do not specify it particularly.

\begin{lemma}
The algebra $\epsilon A \epsilon$ is standardly stratified with respect to the poset $( \{e_{\lambda} \} _{\lambda \in \Lambda_1}, \leqslant) $ and has standard modules $\epsilon A_0 e_{\lambda}$, $\lambda \in \Lambda_1$, which are all concentrated in degree 0. Moreover,
\begin{equation*}
(\epsilon A \epsilon)_0 = \bigoplus _{\lambda \in \Lambda_1} \epsilon A_0 e_{\lambda} = \epsilon A_0 \epsilon = \epsilon A_0
\end{equation*}
is a self-injective algebra. If $M$ is a Koszul $A$-module, then $\epsilon M$ is a Koszul $\epsilon A \epsilon$-module. In particular, all standard modules of $\epsilon A \epsilon$ are Koszul.
\end{lemma}

\begin{proof}
The algebra $\epsilon A \epsilon$ has projective modules $\epsilon A e_{\lambda}$, $\lambda \in \Lambda_1$. Notice that each $A e_{\lambda}$ has a $\Delta$-filtration and the standard module $A_0 e \cong \Delta_{\mu}$ cannot appear in the filtration since $e$ is a minimal primitive idempotent. We conclude that $\epsilon A \epsilon$ has standard modules $\epsilon \Delta_{\lambda} \cong \epsilon A_0 e_{\lambda}$, and $\epsilon A \epsilon$ has a filtration formed by $\epsilon A_0 e_{\lambda}$, $\lambda \in \Lambda_1$. This proves the first statement.

Clearly,
\begin{equation*}
(\epsilon A \epsilon)_0 = \epsilon A_0 \epsilon = \bigoplus _{\lambda \in \Lambda_1} \epsilon A_0 e_{\lambda}.
\end{equation*}
We claim $\epsilon A_0 e =0$, which implies $\epsilon A_0 = \epsilon A_0 \epsilon + \epsilon A_0 e = \epsilon A_0 \epsilon$. Indeed,
\begin{equation*}
\epsilon A_0 e \cong \Hom _{A_0} (A_0 \epsilon, A_0 e) \cong \bigoplus _{\lambda \in \Lambda_1} \Hom_A (\Delta_{\lambda}, \Delta_{\mu}) =0
\end{equation*}
Since $A$ is standardly stratified and $\mu$ is minimal in $\Lambda$.

By Proposition 2.5 on page 35 of \cite{Auslander}, the exact functor $F = \Hom _{A_0} (A_0 \epsilon, -)$ gives an equivalence between a subcategory $\mathcal{M}$ of $A_0$-mod and the category $\epsilon A_0 \epsilon$-mod. Since all projective (injective as well) modules $A_0 e_{\lambda}$ with $\lambda \in \Lambda_1$ are contained in $\mathcal{M}$, and $F$ sends these projective (injective, resp) modules to projective (injective, resp) modules of $\epsilon A_0 \epsilon$. Consequently, $\epsilon A_0 \epsilon = \epsilon A_0$ is self-injective.

Let $M$ be a Koszul $A$-module and
\begin{equation*}
\xymatrix{ \ldots \ar[r] & P^1 \ar[r] & P_0 \ar[r] & M \ar[r] & 0}
\end{equation*}
be a Koszul projective resolution of $M$. That is, each $P^i$ is generated in degree $i$. Applying the exact functor $\Hom_A (A \epsilon, -)$ we get a Koszul projective resolution of $\epsilon M$ as follows
\begin{equation*}
\xymatrix{ \ldots \ar[r] & \epsilon P^1 \ar[r] & \epsilon P^0 \ar[r] & \epsilon M \ar[r] & 0.}
\end{equation*}
Thus $\epsilon M$ is a Koszul $\epsilon A \epsilon$-module. Since the standard modules of $\epsilon A \epsilon$ are indecomposable summands of $\epsilon A_0$, and $A_0 = \Delta$ is a Koszul $A$-module, we conclude that every standard module of $\epsilon A \epsilon$ is a Koszul $\epsilon A \epsilon$-module as well.
\end{proof}

\begin{lemma}
Let $M$ be a Koszul $A$-module up to a degree shift. Then
\begin{enumerate}
\item $M$ has a $\Delta$-filtration.
\item $\Ext_A ^i (M, A_0 e) =0$ for $i > 0$.
\item $[\Omega^i (M): \Delta_{\mu}] = 0$ for each $i \geqslant 1$, where $[\Omega^i (M): \Delta_{\mu}]$ is the number of $\Delta$-filtration factors of $M$ isomorphic to $\Delta_{\mu}$.
\item $M = A \epsilon M$ if and only if $[M: \Delta_{\mu}] =0$.
\end{enumerate}
\end{lemma}

\begin{proof}
Without loss of generality we assume that $M$ is Koszul. By Corollary 2.18 $M$ is a projective $A_0$-module, and hence has a $\Delta$-filtration since $A_0 = \Delta$.

To prove the second statement, it is enough to show $\Ext_A ^1 (M, A_0 e) =0$. Indeed, for $i > 1$, by Corollary 2.5 and using Lemma 2.10 recursively we have
\begin{equation*}
\Ext _A ^i (M, A_0 e) \cong \Ext _A ^{i-1} (\Omega M, A_0 e) \cong \ldots \cong \Ext _A ^1 (\Omega^{i-1} (M), A_0 e).
\end{equation*}
Since $M$ is Koszul, $\Omega^{i-1} (M)$ is also Koszul up to a degree shift. Therefore we can replace $M$ by $\Omega ^{i-1} (M)$ and use induction.

The projective presentation $0 \rightarrow \Omega M \rightarrow P \rightarrow M \rightarrow 0$ gives a surjective map $\Hom_A (\Omega M, A_0 e) \rightarrow \Ext^1 _A (M, A_0 e)$. The syzygy $\Omega M$ is Koszul up to a degree shift and hence has a $\Delta$-filtration. Since $e$ is minimal and the algebra $A$ is standardly stratified with respect to the poset $(\{ e_{\lambda} \} _{\lambda \in \Lambda}, \leqslant)$, every $\Delta$-filtration of $JP = \bigoplus _{i \geqslant 1} P_i$ has no factors isomorphic to $A_0 e$. Therefore, the $\Delta$-filtration of $\Omega M \subseteq \bigoplus _{i \geqslant 1} P_i$ has no factors isomorphic to $A_0 e$ either. But $\Hom_A (\Delta_{\lambda}, \Delta_{\mu}) = 0$ for $\lambda \neq \mu$. Thus $\Hom_A (\Omega M, \Delta_{\mu}) = 0$, so $\Ext^1 _A (M, A_0 e) = 0$. This proves the second statement.

To prove (3), it suffices to show $[\Omega M : \Delta_{\mu}] = 0$ and the conclusion comes from induction. But this fact has been established in last paragraph.

Now we prove (4). Notice that $A \epsilon M$ is the trace of $A \epsilon$ in $M$. If $[M : \Delta_{\mu}] = 0$, then in particular $M_0 \in \text{add} (A_0 \epsilon)$, and $M$ is in the trace of $A \epsilon$ since it is generated in degree 0. Conversely, if $M = A \epsilon M$, then $M$ is in the trace of $A \epsilon$, i.e., $M$ is a quotient module of some $ (A \epsilon) ^{\oplus m}$. Since $[A \epsilon : \Delta_{\mu}] =0$, we deduce that $[M : \Delta_{\mu}] =0$.
\end{proof}

We define an operator $\Pi$ on $A$-gmod as follows:
\begin{equation*}
\Pi(M) = \left\{
\begin{array}{rl}
A \epsilon M & \text{if } M \neq A \epsilon M,\\
\Omega M & \text{if } M = A \epsilon M.
\end{array} \right.
\end{equation*}

\begin{lemma}
Suppose that $M$ and $A \epsilon M$ are Koszul $A$-modules up to a common degree shift. Then
\begin{enumerate}
\item For all $i \geqslant 1$,
\begin{equation*}
\Pi^i(M) = \left\{
\begin{array}{rl}
\Omega^i(M) & \text{if } M = A \epsilon M,\\
\Omega^{i-1} (A \epsilon M) & \text{if } M \neq A \epsilon M.
\end{array} \right.
\end{equation*}
\item There is some $l \in \mathbb{Z}$ such that $\Pi^l (M) = 0$.
\end{enumerate}
\end{lemma}

\begin{proof}
If $M = A \epsilon M$, then $\Pi(M) = \Omega M$. By the previous lemma, $[\Omega M : \Delta_{\mu}] =0$, so $\Omega M = A \epsilon \Omega M$ and $\Pi (\Omega M) = \Omega^2 (M)$. Using induction we get $\Pi^i (M) = \Omega^i (M)$.

If $M \neq A \epsilon M$, then $\Pi(M) = A \epsilon M$. Clearly, $A \epsilon (A \epsilon M) = A \epsilon M$, so $\Pi^2(M) = \Omega (A \epsilon M)$ by the definition. Applying statements (3) and (4) of the previous lemma recursively, we get $\Pi^i (M) = \Omega^{i-1} (A \epsilon M)$.

To prove (2), it suffices to show that there is a number $l \in \mathbb{Z}$ such that $\Pi^l (\Pi(M)) = \Omega^l (\Pi(M)) = 0$. Since $\Pi(M) = \Omega M$ or $\Pi(M) = A \epsilon M$, both of which are Koszul up to a degree shift, $\Pi (M)$ has a $\Delta$-filtration by the above lemma. Since $A$ is standardly stratified and finite-dimensional, $A_0 \cong \Delta$ has finite projective dimension and every $\Delta$-filtration of $\Pi (M)$ is of finite length. Thus the projective dimension of $\Pi (M)$ is finite.
\end{proof}

Recall $\Gamma = \Ext _A ^{\ast} (A_0, A_0)$ and the indecomposable summands of $\Gamma$ are precisely $\Ext _A ^{\ast} (A_0 e_{\lambda}, A_0), \lambda \in \Lambda$.

\begin{lemma}
Suppose that $M$ and $A \epsilon M$ are Koszul $A$-modules up to a common degree shift. Then the trace $\Gamma e \Ext _A ^{\ast} (M, A_0)$ of $\Gamma e$ in $\Ext _A ^{\ast} (M, A_0)$ is a projective $\Gamma$-module. Moreover,
\begin{equation*}
\dim_k \Ext ^{\ast} _A (M, A_0) - \dim_k \Gamma e \Ext ^{\ast} _A (M, A_0) = \dim_k \Ext ^{\ast} _{\epsilon A \epsilon} (\epsilon M, \epsilon A_0).
\end{equation*}
\end{lemma}

\begin{proof}
Without loss of generality we suppose that $M$ is generated in degree 0 and prove this lemma by induction on the least number $l$ such that $\Pi^l (M) = 0$, which always exists by Lemma 7.4. If $l = 0$, then $M =0$ and the conclusion holds trivially. Now suppose that it holds for $l \leqslant n-1$ and let $M$ be a Koszul $A$-module for which this least number is $n$. There are two cases.

(1). If $M \neq A \epsilon M$, then $\Pi(M) = A \epsilon M$. This happens if and only if $[M: \Delta_{\mu}] \neq 0$. Since $e = e_{\mu}$ is a minimal primitive idempotent, those factors isomorphic to $\Delta_{\mu} = A_0 e$ can only appear as direct summands of $M_0$ by the second statement of Lemma 7.3. Therefore we get an exact sequence $0 \rightarrow A \epsilon M \rightarrow M \rightarrow (A_0 e) ^{\oplus a} \rightarrow 0$, where all terms are Koszul by our assumption. By Proposition 2.11, this sequence gives rise to:
\begin{equation*}
0 \rightarrow \Ext ^{\ast} _A ((A_0 e) ^{\oplus a}, A_0) \rightarrow \Ext ^{\ast} _A (M, A_0) \rightarrow \Ext ^{\ast} _A (A \epsilon M, A_0) \rightarrow 0,
\end{equation*}
that is:
\begin{equation}
\xymatrix{ 0 \ar[r] & (\Gamma e) ^{\oplus a} \ar[r] & \Ext ^{\ast} _A (M, A_0) \ar[r]^p & \Ext ^{\ast} _A (A \epsilon M, A_0) \ar[r] & 0.}
\end{equation}
By taking the trace of $\Gamma e$ in those terms appearing in this exact sequence, we get another exact sequence:
\begin{equation}
\xymatrix{ 0 \ar[r] & (\Gamma e) ^{\oplus a} \ar[r] & \Gamma e \Ext ^{\ast} _A (M, A_0) \ar[r] & \Gamma e \Ext ^{\ast} _A (A \epsilon M, A_0) \ar[r] & 0.}
\end{equation}
Notice that $A \epsilon M$ is Koszul, $A \epsilon A \epsilon M = A \epsilon M$, and $\Pi^{n-1} (A \epsilon M) = \Pi^n (M) =0$. By the induction hypothesis, $\Gamma e \Ext ^{\ast} _A (A \epsilon M, A_0)$ is a projective $\Gamma$-module. Then the above sequence splits and $\Gamma e \Ext ^{\ast} _A (M, A_0)$ is a projective $\Gamma$-module as well.

Comparing sequence 7.1 to sequence 7.2, we get:
\begin{align*}
& \dim_k \Ext ^{\ast} _A (M, A_0) - \dim_k \Gamma e \Ext ^{\ast} _A (M, A_0) \\
& = \dim_k \Ext ^{\ast} _A (A \epsilon M, A_0) - \dim_k \Gamma e \Ext ^{\ast} _A (A \epsilon M, A_0) \\
& = \dim_k \Ext ^{\ast} _{\epsilon A \epsilon} (\epsilon A \epsilon M, \epsilon A_0)
\end{align*}
by the induction hypothesis. But $\epsilon A \epsilon M = \epsilon M$. Thus the conclusion is true for $M$ by induction.\\

(2). If $M = A \epsilon M$, i.e., $[M: \Delta_{\mu}] = 0$, then $\Pi(M) = \Omega M$. Let $P$ be a graded projective cover of $M$. Then $P \in \text{add} (A \epsilon)$ and $P_0 \in \text{add} (A_0 \epsilon)$. Since $M$ is Koszul, $M_0 \cong P_0$, so $\Ext _A ^i (M, A_0) = \Ext _A ^{i-1} (\Omega M, A_0)$ for $i>0$ by Lemma 2.10, and
\begin{equation*}
\Hom_A (M, A_0) \cong \Hom_{A_0} (M_0, A_0) \cong \Hom _{A_0} (P_0, A_0) \cong \Hom_A (P_0, A_0).
\end{equation*}
Therefore, the following exact sequence:
\begin{equation*}
\xymatrix { 0 \ar[r] & \bigoplus _{i \geqslant 1} \Ext _A ^i (M, A_0) \ar[r] & \Ext _A^{\ast} (M, A_0) \ar[r] & \Hom_A (M, A_0) \ar[r] & 0}
\end{equation*}
is isomorphic to
\begin{equation}
\xymatrix { 0 \ar[r] & \Ext _A ^{\ast} (\Omega M, A_0) \ar[r] & \Ext _A^{\ast} (M, A_0) \ar[r] & \Hom_A (P_0, A_0) \ar[r] & 0. }
\end{equation}
Applying the exact functor $\Hom _{\Gamma} (\Gamma e, -)$ to the above sequence, we get
\begin{equation*}
\xymatrix { 0 \rightarrow e \Ext _A ^{\ast} (\Omega M, A_0) \rightarrow e \Ext _A^{\ast} (M, A_0) \rightarrow e \Hom_A (P_0, A_0) \rightarrow 0.}
\end{equation*}
Notice that $e \Hom_A (P_0, A_0) = 0$ since $A_0 e = \Delta_{\mu}$, $P_0 \in \text{add} (A_0 \epsilon)$ and $[A_0 \epsilon : \Delta_{\mu}] = 0$. Thus $e \Ext _A ^{\ast} (M, A_0) \cong e \Ext _A ^{\ast} (\Omega M, A_0) [-1]$ (as graded $\Gamma$-modules) and
\begin{equation}
\Gamma e \Ext _A ^{\ast} (M, A_0) \cong \Gamma e \Ext _A ^{\ast} (\Omega M, A_0)
\end{equation}
Since $[\Omega M : \Delta_{\mu} ] = 0$, by Lemma 7.3, $A \epsilon \Omega M = \Omega M$. So $\Gamma e \Ext _A ^{\ast} (\Omega M, A_0)$ is a projective $\Gamma$-module by the induction hypothesis on $\Pi(M) = \Omega M$, and $\Gamma e \Ext _A ^{\ast} (M, A_0)$ is projective as well.

The exact sequence $ 0 \rightarrow \Omega M \rightarrow P \rightarrow M \rightarrow 0$ gives the exact sequence $0 \rightarrow \epsilon \Omega M \rightarrow \epsilon P \rightarrow \epsilon M \rightarrow 0$. Clearly, $\epsilon P$ is a projective cover of $\epsilon M$, so $\epsilon \Omega M = \Omega (\epsilon M)$. Since $\epsilon M$ is a Koszul $\epsilon A \epsilon$-module by Lemma 7.2, we have an exact sequence similar to (7.3):
\begin{equation}
0 \rightarrow \Ext _{\epsilon A \epsilon} ^{\ast} (\Omega(\epsilon M), \epsilon A_0) \rightarrow \Ext _{\epsilon A \epsilon} ^{\ast} (\epsilon M, \epsilon A_0) \rightarrow \Hom _{\epsilon A \epsilon} (\epsilon P_0, \epsilon A_0) \rightarrow 0.
\end{equation}
Therefore,
\begin{align*}
& \dim_k \Ext ^{\ast} _A (M, A_0) - \dim_k \Gamma e \Ext ^{\ast} _A (M, A_0) \quad \text{by (7.3) and (7.4)} \\
& = \dim_k \Ext ^{\ast} _A (\Omega M, A_0) + \dim_k \Hom_A (P_0, A_0) - \dim_k \Gamma e \Ext ^{\ast} _A (\Omega M, A_0) \\
&= \dim_k \Ext ^{\ast} _{\epsilon A \epsilon} (\epsilon \Omega M, \epsilon A_0) + \dim_k \Hom_A (P_0, A_0) \quad \text{by induction on } \Omega M \\
& = \dim_k \Ext ^{\ast} _A (\epsilon M, \epsilon A_0) + \dim_k \Hom_A (P_0, A_0) - \dim_k \Hom _{\epsilon A \epsilon} (\epsilon P_0, \epsilon A_0),
\end{align*}
where the last identity comes from (7.5).

We establish the identity $\dim_k \Hom_A (P_0, A_0) = \dim_k \Hom_{\epsilon A \epsilon} (\epsilon P_0, \epsilon A_0)$ and finish the proof by induction. Take an arbitrary indecomposable summand of $P_0$, which is isomorphic to a certain $A_0 e_{\lambda}$. Since $A \epsilon M = M$, $P_0 \cong M_0$ has no summands isomorphic to $A_0 e$. Therefore $\lambda \in \Lambda_1$ and $e_{\lambda} \epsilon = e_{\lambda}$. Now
\begin{align*}
\Hom_A (A_0 e_{\lambda}, A_0) & \cong \Hom_{A_0} (A_0 e_{\lambda}, A_0) \cong e_{\lambda} A_0 = e_{\lambda} \epsilon A_0\\
& \cong \Hom_{\epsilon A_0 \epsilon} (\epsilon A_0 e_{\lambda}, \epsilon A_0) \cong \Hom _{\epsilon A \epsilon} (\epsilon A_0 e_{\lambda}, \epsilon A_0).
\end{align*}
Since $P_0$ is a direct sum of these summands, the identity holds.
\end{proof}

Now we can prove the second statement of Theorem 7.1.

\begin{proof}
We use induction on the size of the poset $\Lambda$. If this number is 1, then $A$ and $\Gamma = \Lambda ^{\textnormal{op}}$ both are local algebra. Clearly $\Gamma$ is standardly stratified. Suppose that the conclusion is true for posets with sizes at most $m-1$ and let $\Lambda$ be an poset with $m$ elements. Take $e$ be a minimal idempotent and define $\epsilon, \Lambda_1$ as before.

Let $M = A_0$ in the previous lemma. We conclude that $\Gamma e \Gamma$ is a projective $\Gamma$-module. By Lemma 7.2 $\epsilon A \epsilon$ satisfies all the conditions in the theorem and $\Lambda_1$ has only $m-1$ elements. Therefore, by the induction hypothesis, $\Gamma' = \Ext ^{\ast} _{\epsilon A \epsilon} (\epsilon A_0, \epsilon A_0)$ is standardly stratified with respect to the poset $(\{ e_{\lambda} \} _{\lambda \in \Lambda_1}, \leqslant ^{\textnormal{op}} )$. Thus it is enough to show $\Gamma / \Gamma e \Gamma \cong \Gamma'$.

There is an algebra homomorphism $\varphi: \Gamma \rightarrow \Gamma'$ induced by the exact functor $F = \Hom_A (A \epsilon, -)$ in the following way: $\varphi$ sends an $n$-fold exact sequence
\begin{equation*}
\xymatrix {0\ar[r] & A_0 \ar[r] & M^n \ar[r] & \ldots \ar[r] & M^1 \ar[r] & A_0 \ar[r] & 0}
\end{equation*}
representing an element $g \in \Gamma_n$ to an $n$-fold exact sequence
\begin{equation*}
\xymatrix {0 \ar[r] & \epsilon A_0 \ar[r] & \epsilon M^n \ar[r] & \ldots \ar[r] & \epsilon M^1 \ar[r] & \epsilon A_0 \ar[r] & 0}
\end{equation*}
representing an element $g' \in \Gamma'_n$ for all $n \geqslant 0$. Every element $x \in (\Gamma e)_n = \Ext _A^n (A_0e, A_0)$ can be represented by an exact sequence:
\begin{equation*}
\xymatrix {0\ar[r] & A_0 \ar[r] & M^n \ar[r] & \ldots \ar[r] & M^1 \ar[r] & A_0e \ar[r] & 0.}
\end{equation*}
Since $\Hom_A (A\epsilon, -)$ sends $A_0e$ to $\Hom_A (A\epsilon, A_0e) =0$, $\varphi$ maps every $(\Gamma e)_n$ and hence $\Gamma e$ to 0. Thus the ideal $\Gamma e \Gamma$ generated by $\Gamma e$ is also sent to 0, and $\varphi$ gives rise to an algebra homomorphism $\bar{\varphi}: \Gamma/ \Gamma e \Gamma \rightarrow \Gamma'$.

Clearly, $\varphi$ maps $\Gamma_0 = \End_{A_0} (A_0)$ onto $\End _{\epsilon A \epsilon} (\epsilon A_0) = \epsilon \Gamma_0 \epsilon$. Moreover, we have
\begin{equation*}
\Ext _A^1 (A_0, A_0) \cong \Hom_A (\bigoplus _{i \geqslant 1} A_i, A_0) \cong \Hom _{A_0} (A_1, A_0) \cong \Hom _{\epsilon A_0 \epsilon} (\epsilon A_1, \epsilon A_0),
\end{equation*}
where the last isomorphism is induced by the functor $F$ (see the last paragraph of the proof of Lemma 7.5). But
\begin{equation*}
\Hom _{\epsilon A_0 \epsilon} (\epsilon A_1, \epsilon A_0) \cong \Hom _{\epsilon A \epsilon} (\bigoplus _{i \geqslant 1} \epsilon A_i, \epsilon A_0) \cong \Ext _{\epsilon A \epsilon} ^1 (\epsilon A_0, \epsilon A_0).
\end{equation*}
Thus, the homomorphism $\varphi$ induced by $F$ maps $\Ext _A^1 (A_0, A_0)$ to $\Ext _{\epsilon A \epsilon} ^1 (\epsilon A_0, \epsilon A_0)$ surjectively. Since $A_0$ ($\epsilon A_0$, resp.) is a Koszul $A$-module ($\epsilon A \epsilon$-module, resp.), both $\Gamma = \Ext _A  ^{\ast} (A_0, A_0)$ and $\Gamma' = \Ext _{\epsilon A \epsilon} ^{\ast} (\epsilon A_0 \epsilon, \epsilon A_0 \epsilon)$ are generated in degree 0 and degree 1 as algebras. Therefore, the map $\varphi: \Gamma \rightarrow \Gamma'$ is a surjective algebra homomorphism, so $\bar{\varphi}: \Gamma / \Gamma e \Gamma \rightarrow \Gamma'$ is surjective as well.

Let $M = A_0$ in the previous lemma. We get dim$_k \Gamma' = \dim_k \Gamma - \dim_k \Gamma e \Gamma$. Therefore, as a surjective homomorphism between two $k$-algebras with the same dimension, $\bar{\varphi}$  must be an isomorphism. This completes the proof.
\end{proof}

There is a similar conclusion in the case that $A$ is a \textit{quasi-hereditary} algebra, the definition for which can be found in \cite{Dlab}.

\begin{corollary}
Let $A$ be the same as in the previous theorem. If $A$ is quasi-hereditary with respect to the poset $(\{ e_{\lambda } \}_{\lambda \in \Lambda}, \leqslant)$, then $\Gamma = \Ext _A ^{\ast} (A_0, A_0)$ is also quasi-hereditary with respect to the poset $(\{ e_{\lambda } \}_{\lambda \in \Lambda}, \leqslant^{\textnormal{op}})$.
\end{corollary}

\begin{proof}
We already proved that $\Gamma$ is standardly stratified with respect to the poset $(\{ e_{\lambda } \}_{\lambda \in \Lambda}, \leqslant^{\textnormal{op}})$. It suffices to check that the ideal $\Gamma e \Gamma$ is a hereditary ideal, i.e, the endomorphism algebra of $\Gamma e$ is one-dimensional.

Since $e$ is maximal as an idempotent of $\Gamma$, the standard $\Gamma$-module corresponding to $e$ is exactly $\Gamma e$. Thus
\begin{equation*}
\Hom _{\Gamma} (\Gamma e, \Gamma e) \cong e \Gamma e \cong \Ext _A ^{\ast} (A_0e, A_0e) = \Hom_A (A_0e, A_0e)
\end{equation*}
where the last identity follows from (2) of Lemma 7.3. Since $A_0e$ is a standard module of the quasi-hereditary algebra $A$, $\Hom_A (A_0e, A_0e) \cong k$. The conclusion follows from induction.
\end{proof}

\end{document}